\documentclass[10pt, a4paper]{article}

\usepackage{graphicx}
\usepackage{color}
\usepackage{amsmath}
\usepackage{amssymb}
\usepackage{amscd}
\usepackage{bbm}
\usepackage{tikz}
\usepackage{hyperref}
\hypersetup{
    bookmarks=true,         
    colorlinks=true,       
    linkcolor=orange,          
    citecolor=green,        
    filecolor=magenta,      
    urlcolor=cyan           
}

\usepackage[utf8]{inputenc}

\usepackage{marginnote}

\usepackage{amsthm}
\usepackage{bbm} 

\newtheorem{Theorem}{Theorem}[section]
\newtheorem{Lemma}[Theorem]{Lemma}
\newtheorem{Corollary}[Theorem]{Corollary}
\newtheorem{Assumption}[Theorem]{Assumption}
\newtheorem{Proposition}[Theorem]{Proposition}
\newtheorem{Definition}[Theorem]{Definition}

\newtheorem{design-crit}[Theorem]{Design Criterion}

\theoremstyle{remark}

\newtheorem{Remark}{Remark}

\newcommand{\inr}[1]{\bigl< #1 \bigr>}
\newcommand{\biginr}[1]{\left< #1 \right>}

\newcommand{\bignorm}[1]{\left\|#1\right\|}%
\newcommand{\norm}[1]{\|#1\|}%
\newcommand{\Tr}[1]{{\mathrm{Tr}\left(#1\right)}}%

\newcommand\ove[2]{\overset{#2}{#1}\;}

\newcommand{\vertiii}[1]{{\left\vert\kern-0.25ex\left\vert\kern-0.25ex\left\vert #1 
    \right\vert\kern-0.25ex\right\vert\kern-0.25ex\right\vert}}

\newcommand{\beginproof}{{\bf Proof. {\hspace{0.2cm}}}}
\def \endproof
{{\mbox{}\nolinebreak\hfill\rule{2mm}{2mm}\par\medbreak}}

\DeclareMathOperator*{\argmax}{argmax}
\DeclareMathOperator*{\argmin}{argmin}

\def\ds1{\textrm{1\kern-0.25emI}} 

\newcommand \E{\mathbb{E}}
\newcommand \R{\mathbb{R}}

\newcommand \cB{{\cal B}}
\newcommand \cC{{\cal C}}

\newcommand \cE{{\cal E}}

\newcommand \cI{{\cal I}}

\newcommand \cL{{\cal L}}

\newcommand \cO{{\cal O}}
\newcommand \cP{{\cal P}}

\newcommand \cS{{\cal S}}

\newcommand \cX{{\cal X}}
\newcommand \cY{{\cal Y}}
\newcommand \cZ{{\cal Z}}

\newcommand \bC{{\mathbb C}}

\newcommand \bE{{\mathbb E}}

\newcommand \bN{{\mathbb N}}

\newcommand \bP{{\mathbb P}}

\newcommand \bR{{\mathbb R}}

\newcommand{\bV}{{\boldsymbol{V}}}

\newcommand{\simplequote}[1]{`#1'}
\newcommand{\doublequote}[1]{``#1''}

\newcommand{\Perm}{\mathfrak{S}}

\usepackage{comment}
\usepackage{subcaption}
\usepackage{caption}
\usepackage[skip=-5pt, font={small}]{subcaption}
\usepackage{algorithmic}
\usepackage{mathabx}
\usepackage[nolabel, final]{showlabels}

\usepackage{multirow,bigdelim}
\usepackage{array}     

\usepackage{tikz}

\usepackage{bbold}
\usepackage{tikz-cd}
\usepackage{multirow}

\begin{document}

\title{Learning with a linear loss function. Excess risk and estimation bounds for ERM, minmax MOM and their regularized versions. Applications to robustness in sparse PCA.}
\author{Guillaume Lecué and Lucie Neirac  \\
ESSEC, Business school and CREST, ENSAE, IPParis.}


\maketitle

\begin{abstract}
Motivated by several examples, we consider a general framework of learning with linear loss functions. In this context, we provide excess risk and estimation  bounds that hold with large probability for four estimators: ERM, minmax MOM and their regularized versions. 
These general bounds are applied for the problem of robustness in sparse PCA. In particular, we improve the state of the art result for this this problems, obtain results under weak moment assumptions as well as for adversarial contaminated data.
\end{abstract}

{\bf{keywords:}} SDP relaxation, empirical processes, robustness, heavy-tailed, adversarial contamination, high-dimensional statistics.


\section{Introduction}\label{sec:intro}

Community detection, phase recovery, signed clustering, angular group synchronization, \textsc{Max-Cut}, sparse PCA, and the sparse single index model are all classical topics in machine learning and statistics. At first glance, they are pretty different problems with different types of data and different goals. However, they can all be written in such a way that a common analysis of various estimators introduced for these problems can be analyzed the same way. All these problems can be recast in the classical machine learning framework of risk minimization \cite{MR1719582}. It is therefore possible to leverage the vast literature related to risk minimization to derive excess risk and estimation bounds as well as algorithms for all the problems cited above as well as many other ones. It appears that the general framework that can encapsulate all these problems relies in fact on a simple loss function, maybe the simplest one: the linear loss function. Indeed, this observation is the baseline of \cite{MR4329809}: several estimators introduced recently in some of the problems cited at the beginning are in fact empirical risk minimizers (ERM) for linear loss functions. They can therefore be analyzed using all the machinery \cite{MR1719582,BoucheronLugosiMassart,Kol11} developed during the last fifty years for ERM in this very specific framework of the linear loss function. 

General excess risk and estimation bounds have therefore been obtained in \cite{MR4329809} for ERM using a linear loss function. State-of-the art techniques like localization, homogeneity argument, local curvature and complexity fixed points equation have been used in \cite{MR4329809} to obtain these general bounds that have then been applied in Community detection, phase recovery, signed clustering, angular group synchronization and \textsc{Max-Cut}. This new perspective allowed us to obtain new results or recover older one with a new proof technique but most importantly it showed that a common analysis of several problems that looks a priori very different can be performed. 

The aim of the article is to push forward the analysis of statistical procedures based on linear loss functions and to show that this viewpoint allows to deal with the problem of structural risk minimization and of robustness\footnote{in all this article, robustness means robust to data contamination and to heavy-tailed data.} in all the problems cited above and in many other ones (some of them are given below). Indeed, \cite{MR4329809} only deals with ERM procedures. However, some problems rely on some structure such as sparsity and other are facing the problem of robustness. For these issues, ERM is not the right answer and these two problems call for other procedures such as regularized ERM (for structural risk minimization) or the recently introduced minmax MOM estimator \cite{MR4102681} (for robustness issues). It is therefore the first contribution of this paper to derive general statistical bounds for regularized ERM, minmax MOM and its regularized version in the framework of linear loss functions. As an illustration these bounds are applied to the problem of sparse PCA. Using our viewpoint, we improve state-of the art results for this problem (improvement on the rates and the deviation for less stringent assumptions) as well as getting robust (to heavy-tailed data and to adversarial contamination) versions of these results thanks to the minmax MOM approach. Another aim of this article is to show that the linear loss function appears in many problems and so we provide a list of problems that can be recast in this framework. But first, we explain how linear loss function appear only recently, even though they are simpler than many other loss functions previously used in machine learning such as the quadratic or the logistic loss functions.    

Statistics, machine learning and optimization got closer during the last twenty years and gave birth in part to \textit{data sciences}. One consequence of these connections is that nowadays statistical estimators and machine learning procedures should be computable on a laptop in a reasonable amount of time and should not be purely theoretical objects. This viewpoint shed some lights on algorithms from the statistical perspective and may now be seen as statistical procedures that can receive a statistical analysis such as satisfying excess risk bounds. For instance, statistical properties of some gradient-descent based algorithms and  SDP relaxation procedures have been obtained during the last twenty years. In particular, the SDP relaxation has proved to be very successful first in optmization and nowadays in statistics for many graph related issues such as community detection. From our perspective, SDP relaxation has been at the origin of many examples of ERMs based on a linear loss function. 

Semidefinite programming (SDP) as a mathematical concept was introduced in the late 1980s and early 1990s. The foundations of SDP were laid down by researchers such as Yurii Nesterov, Arkadi Nemirovski, and others (\cite{NEMIROVSKII198521}, \cite{Boyd1997}, \cite{Nesterov1998SemidefiniteRA}), who extended the ideas of linear programming to semidefinite matrices, allowing for the optimization of linear functions subject to semidefinite constraints. The theoretical development and algorithms for solving SDP problems gained significant attention during this period, leading to its establishment as a fundamental optimization framework within the mathematical community.

The growing interest in SDPs in recent years is due to several compelling factors. One of the main factors is its broad applicability, as it can address a wide variety of complex problems arising in various mathematical contexts, including graph theory \cite{Gaar_2022,10.1007/978-3-642-01929-6_8}, combinatorial optimization \cite{gutekunst2019semidefinite}, signal processing \cite{Luo2006AnIT}, quantum information \cite{Wang_2016}, for the {K}oml\'{o}s conjecture \cite{MR3945254} or in integer programming \cite{MR3500324}. Its potency lies in its ability to efficiently handle non-convex and combinatorial optimization challenges by approximating them with convex semidefinite constraints. At the same time, the development of efficient algorithms for solving SDP problems, such as interior-point methods \cite{Helmberg1996} and first-order methods \cite{Monteiro2003}, has significantly improved the feasibility of tackling large-scale SDPs, thereby broadening the range of possibilities for applying SDP to real-world problems. 

From our point-of-view SDP relaxation provide many examples of machine learning procedures such as ERM or RERM (regularized ERM) based on a linear loss function. We are now providing some of these examples and later we will dive deeper into the example of sparse PCA.

\paragraph{Notations.} Throughout this paper, we use uppercase letters for matrix and lowercase letters for vectors. For a matrix $A\in\bR^{N\times P}$, we note $A\geq 0$ to indicate that $A_{ij}\geq 0$ for any $(i,j)\in\{1, \ldots, N\}\times\{1, \ldots, P\}$ and $A\succeq 0$ to say that $A$ is positive semidefinite. For $A$ and $B\in\bR^{N\times P}$, we define their Frobenius inner product as $\inr{A, B} := \mathrm{Trace}(B^\top A)$, and we write $A\circ B$ for their element-wise product. If $x$ is a vector in $\bC^d$ then $|x|$ denotes the vector in $\bR^d$ made of the modules of the coordinates of $x$. We denote $[N]=\{1, \ldots, N\}$.

\paragraph{Community detection.}  SDPs have been used to handle the problem of community detection on graphs in \cite{guedon2016community} or \cite{MR3901009} under the Stochastic Block Model assumption, which is as follows. We consider a set of vertices $V=\{1, \cdots, d\}$, and assume it is partitioned into $K$ communities $\cC_1, \cdots, \cC_K$ of arbitrary sizes $|\cC_1|, \cdots, |\cC_K|$. For any pair of nodes $i,j\in V$, we denote by $i \sim j$ when $i$ and $j$ belong to the same community, and by $i\not\sim j$ if $i$ and $j$ do not belong to the same community. For each pair $(i, j)$ of nodes from $V$, we draw an edge between $i$ and $j$ with a fixed probability $p_{ij}$ independently from the other edges. We assume that there exist numbers $p$ and $q$ satisfying $0<q<p<1$, such that $p_{ij} > p $ if $i\sim j$ and $i\neq j$, $p_{ij} = 1$ if $i = j$ and $p_{ij} < q$ otherwise. We denote by $A=(A_{i,j})_{1\leq i,j,\leq d}$ the observed symmetric adjacency matrix, such that, for all $1\leq i\leq j\leq d$,  $A_{ij}$ is distributed according to a Bernoulli of parameter $p_{ij}$. The c{}ommunity structure of such a graph is captured by the membership matrix $\bar{Z} \in \R^{d\times d}$, defined by $\bar{Z}_{ij}=1$ if $i\sim j$, and $\bar{Z}_{ij}=0$ otherwise. The objective is to reconstruct $\bar{Z}$ from the observation $A$. Lemma 7.1 of \cite{guedon2016community} shows that the membership matrix $\bar Z$ is given by the following oracle:
\begin{align*}
Z^*\in \argmax_{Z\in\cC} \inr{\bE[A], Z} ,\qquad \cC := \left\{Z \in \R^{d\times d}, Z \succeq 0, Z \geq 0, {\mathrm diag}(Z)\preceq I_d, \sum_{i, j=1}^d Z_{ij} \leq \lambda\right\} 
\end{align*}
where $\lambda = \sum_{i, j=1}^d \bar{Z}_{ij} = \sum_{k=1}^K |\mathcal{C}_k|^2$ denotes the number of nonzero elements in the membership matrix $\bar{Z}$. Since only the A matrix is observed, the authors consider the following estimator for $Z^*$:
\begin{align*}
\hat Z\in \argmax_{Z\in\cC} \inr{A, Z}. 
\end{align*}
This estimator is therefore obtained as the solution of an ERM with the linear loss function $Z\to\ell_Z(A):=-\inr{A, Z}$, constructed from a single observation of the random matrix $A$.

\paragraph{Variable clustering.} SDP estimators have been used in \cite{bunea2018model} to solve the variable clustering problem. The problem is that of grouping into clusters similar components of a vector $X\in\bR^d$, that is to find a partition $G = \{G_1, \ldots, G_K\}$ of $\{1, \ldots, d\}$ that separates the components of $X$. To that end, the authors observe $N$ independant copies $X_1, \ldots, X_N$ of $X$ and place themselves in the case where the covariance matrix $\Sigma$ of $X$ follows a block model. To describe this model, we need to define the membership matrix $Q \in \bR^{p\times K}$ associated with a partition $G$ as $Q_{ak} = \mathbbm{1}_{\left\{a \in G_k\right\}}$. Then, $\Sigma$ is said to follow an exact $G$-block covariance model when it decomposes as $\Sigma = QCQ^\top + \Gamma$, where $C$ is a symmetric $K\times K$ matrix and $\Gamma$ is a diagonal $d\times d$ matrix. For a given partition $G$, we also introduce its corresponding membership matrix $Z^*\in\bR^{d\times d}$ defined by $Z^*_{ij} = |G_k|^{-1}\mathbbm{1}_{\left\{i \mbox{ and } j \mbox{ belong to the same group } G_k\right\}}$. There is a one-to-one correspondence between partitions $G$ and their corresponding membership matrices, so that looking for $G$ is equivalent to looking for $Z^*$. Using the $K$-means algorithm and a relaxation of it given in \cite{PengWei2007}, the authors show that the best partition for the $X_i$'s can be estimated with the one corresponding to the following membership matrix:
\begin{align*}
\hat Z\in \argmax_{Z\in\cC} \inr{A, Z} ,\qquad \cC := \left\{Z\in\bR^{d\times d}: Z\succeq 0, Z \geq 0 , \sum_{j}Z_{ij} = 1 \forall i, \mathrm{Tr}\left(Z\right) = K\right\} 
\end{align*}
where $A:=\frac{1}{N}\sum_{i = 1}^NX_iX_i^\top$ is the empirical covariance of the $X_i$'s. In the noiseless case, we would have $Z^*\in \argmax_{Z\in\cC} \inr{\bE[A], Z}$. The estimator $\hat Z$ can therefore be seen as an ERM with the linear loss function $Z\to\ell_Z(A):=-\inr{A, Z}$, constructed from the observation of $A$.

\paragraph{Angular synchronization.}

The angular synchronization problem consists of estimating $d$ unknown angles $\theta_1, \cdots, \theta_d$ (up to a global shift angle) given a noisy subset of their pairwise offsets $\delta_{ij} = \theta_i - \theta_j$. This problem is investigated in \cite{ThightnessAngularSynchr2016}. The authors consider that they observe $d(d-1)/2$ measurements of the following form:
\[
a_{ij} = e^{\iota\delta_{ij}} + \epsilon_{ij}, \quad \mbox{ for } 1\leq i < j \leq d.
\]  
They assume the $(\epsilon_{ij})_{i<j}$'s to be $i.i.d$ complex Gaussian variables. The problem can be rewritten under the following form:
\[
A = X\bar X^\top + \sigma W
\]
with $X\in\bC^d$ defined by $X_i = e^{\iota\theta_{i}}$, $W$ being a complex Wigner matrix and $\sigma > 0$ being the variance of the noise. The aim is then to reconstruct the vector $x^*=(e^{\iota \theta_i})_{i=1}^d$, whose maximum likelihood estimator is, up to a global rotation of its coordinates, the unique solution of the following maximization problem:
\[
\argmax_{x \in \cE} \left\{ \bar{x}^\top \; \bE A \; x \right\} \mbox{ where } \cE := \left\{x \in \mathbb{C}^d:|x_i| = 1 \mbox{ for all } i = 1, \ldots,d \right\}
\]
By noticing that $ \cE = \{ Z \in \mathbb{H}_n: Z \succeq 0, \mathrm{diag}(Z) = \mathbbm{1}_d, \mathrm{rank}(Z) = 1 \}$, they lead to the following SDP formulation of the problem, after removing the rank constraint:
\begin{align}\label{eq:angular_sync_oracle}
Z^* \in \argmin_{Z\in\cC} \left(-\inr{\bE[A], Z}\right) \mbox{ where } \cC := \{ Z \in \mathbb{H}_n: Z \succeq 0, \mathrm{diag}(Z) = \mathbbm{1}_d\}.
\end{align}
They show that in this setting, $x^*$ can be obtain from $Z^*$ as its leading unit-length eigen vector. Since $\bE[A]$ is not known and only observed through $A$, $\hat Z \in \argmin_{Z\in\cC} \left(-\inr{A, Z}\right)$ is a natural estimator for $Z^*$. This is therefore another example of an ERM estimator based on the observation of the matrix $A$ and the linear loss function $Z\to \ell_Z(A) = -\inr{A, Z}$.

\paragraph{\textsc{Max-Cut}.} 

In \cite{hong2021optimal}, the authors propose an SDP estimator to handle the \textsc{Max-Cut} problem. The \textsc{Max-Cut} problem is a classical graph theory problem, which consists of taking a graph with vertices $V := \{1, \ldots, d\}$ and edges $E\subset V\times V$ and finding a partition $S\cup\bar{S}=V$ of vertices such that the number of edges connecting a vertex in $S$ to a vertex in $\bar{S}$ is maximal among all possible partitions. Most of the time, we observe only $A\in\{0, 1\}^{d\times d}$ a noisy or partial version of the adjacency matrix of the graph. Hence, the true adjacency matrix of the graph is not observed but it is usually assumed to be equal to the expectation $\bE A$ of the observed one $A$. Hence, $A$ is considered as our data and from this data, we wish to find an optimal partition $S^*$ of the original graph. Choosing a partition $S$ being equivalent to choosing $x\in\{-1, 1\}^N$, it is shown in \cite{goemans1995improved} via a lifting argument that an optimal partition is a first eigenvector of a solution to the following optimization problem:
\[
Z^*\in\argmin_{Z\in\bR^{d\times d}}\left(\inr{\bE[A], Z}: Z\succeq 0, Z_{ii} = 1~\forall i, \mathrm{rank}(Z)=1\right).
\]
 Then, using an SDP relaxation by removing the rank constraint, we recover the classical \textsc{Max-Cut} SDP relaxation procedures introduced by Goemans and Williamson. The ERM counterpart based on the data $A$ is
\begin{equation*}
 \hat Z\in\argmin_{Z\in\cC}\inr{A, Z} \mbox{ for } \cC:=\left\{Z\in\bR^{d\times d}, Z\succeq 0, Z_{ii} = 1~\forall i\right\}
 \end{equation*}It is indeed an ERM procedure based on the observation of $A$ and the linear loss function $Z\to\ell_Z(A):=\inr{A, Z}$ over a convex set.

 \paragraph{Phase recovery.} 

 The former problem is close to the one of phase recovery, which aims at recovering a vector $x\in \bC^d$ from the noisy observation of the amplitude of $N$ random linear measurements: $X = |Bx| \in\bR^N$, with $B\in\bC^{N\times d}$ a random matrix. In \cite{waldspurger2013phase}, the authors use a strategy that involves separating phase from amplitude and optimizing only the values of the phase variables. In the noiseless case, they write $x = B^+\mathrm{diag(X)u}$, where $u\in\bC^N$ is a phase vector and $B^+\in\bC^{N\times N}$ is the pseudo-inverse of $B$. In this format, they show that finding $x\in\bC^d$ such that $|Bx| = X$ is equivalent to solving the following problem:
\begin{align*}
   z^* \in \argmin_{ z\in\cE}  \inr{\bE[A], z\bar{z}^\top} \mbox{ where } \cE := \left\{z\in\bC^N: \,|z_i| = 1\,,\forall i \in [N]\right\} 
\end{align*}
and $A := \left(XX^\top\right) \circ \left(I_N - BB^+\right)$. Writing $Z = z\bar{z}^\top$, this problem is equivalent to the following one:
\begin{align*}
   \min \left( \inr{\bE[A], Z}: Z\succeq 0, Z_{ii} = 1 \forall i , \mathrm{rank}(Z) = 1\right) 
\end{align*}
which may be relaxed by dropping the rank constraint:
\[
Z^*\in \argmin_{Z\in\cC} \inr{\bE[A], Z} \mbox{ for } \cC := \left\{Z\in\bR^{N\times N}: Z\succeq 0, Z_{ii} = 1~\forall i \in [N]\right\}. 
\] The optimal value of $z^*$ is then obtained as the first eigenvector of the oracle $Z^*$. An estimator of $Z^*$ from the observation of $A$ is then $ \hat Z \in \argmin_{ Z\in\cC} \inr{A, Z}$ which is a SDP optimization problem that we see as an ERM with the linear loss function $Z\to\ell_Z(A):=\inr{A, Z}$.

\paragraph{Distance metric learning} 

SDP estimators can also be used in learning distance metrics, as it is done in \cite{Xing2002}. Learning distances is particularly important, as the choice of a metric that is correctly adapted to the input space is crucial to the acuity of many learning algorithms, especially in clustering, where it is essential to take deep account of the relationships between the data. Let's consider a set of points $\left(X_i\right)_{i=1, \ldots, N}\in\bR^d$ that we observe partially or with noise. Now, consider the task of learning a distance metric of the form
\[
d_Z(X, Y) = \sqrt{\mathrm{Tr((X-Y)(X-Y)^\top Z)}},
\]where $Z\succeq 0$ is positive semidefinite. We note that, since one has $\mathrm{Tr((X-Y)(X-Y)^\top Z)} = \norm{Z^{1/2}(X-Y)}_2^2$, learning such a distance metric amounts to finding a rescaling of data that replaces each point $X$ with $Z^{1/2}X$ and applying the standard Euclidean metric to the rescaled data. Now, assume that we want the $X_i$'s to be as close as possible to each other for this metric. This leads us to solve the problem $\min_{Z\succeq 0}\sum_{i,j=1}^Nd_Z(X_i, X_j)^2$. However,  this last problem is trivially solved by $Z=0$ hence, we may add some constraints: we suppose to know $M$ points $(Y_i)_{i=1, \ldots,M}$, distinct from the $X_i$'s, for which we want $\sum_{i, j = 1}^Md_Z(Y_i, Y_j)\geq 1$ to be satisfied. This prevent the situation where $d_Z$ collapses the dataset into a single point. Let us then define $A := \sum_{i, j = 1}^N\left(X_i-X_j\right)\left(X_i-X_j\right)^\top$. In the noiseless case, the matrix $Z^*$ we are looking for can then be taken as a solution to the following problem:
\begin{align*}
    Z^*\in\argmin_{Z\in\cC} \inr{\bE[A], Z} \mbox{ where } \cC := \left\{Z\in\bR^{d\times d} : Z\succeq 0, \sum_{i, j=1}^M\inr{\left(Y_i-Y_j\right)\left(Y_i-Y_j\right)^\top, Z   }^{1/2}\geq 1\right\}.
\end{align*}
One can show that the set $\cC$ is convex (see Appendix \ref{app:distance_metric_learning_convexity}). In practice, the observation $A$ is a noisy version of $\bE[A]$, so we just replace $\bE[A]$ with $A$ to get an estimator of $Z^*$: $\hat Z\in\argmin_{Z\in\cC} \inr{A, Z}$ which is again an ERM estimator with the linear loss function $Z\to\ell_Z(A) = \inr{A, Z}$, constructed from an observation of the random matrix $A$. 

\paragraph{Noisy optimal transport.}

Let $\cX = (x_1, \ldots, x_N)$ and  $\cY=(y_1, \ldots, y_N)$ be two clouds of points in $\bR^d$. The quadratic optimal transport problem (or quadratic assignment problem) is defined by the $W_2$-Wasserstein distance
\begin{equation}\label{eqn:def:W2}
    W_2^2(\cX,\cY) = \min_{\tau \in \Perm_N} \sum_{i=1}^N \|x_i - y_{\tau(i)}\|^2
\end{equation}where $\Perm_N$ is the set of all permutations of $[N]$. Finding a solution to \eqref{eqn:def:W2} is a standard problem in optimal transport that can be lifted to the matrix problem 
\begin{equation*}
Z^*\in\argmin_{Z\in \cC} \sum_{i,j}\norm{x_i-y_j}_2^2 P_{ij}
\end{equation*} and $\cC$ is the set of all $N\times N$ bi-stochastic matrices (i.e. of matrices with non-negative entries summing to one along rows and columns). Indeed, if $\tau^*$ denotes an optimal solution to \eqref{eqn:def:W2} then for all $i\in[N]$, $Z^*_{i \tau^*(i)}=1$ and $Z^*_{i j}=0$ when $j\neq \tau^*(i)$. 

Let us now assume that we do not observe exactly the points in $\cX$ and $\cY$ but we only have access to a noisy version of these points: for all $i\in[N]$, $X_i=x_i+\sigma G_i$ and $Y_i=y_i+\sigma G_i^\prime$ where $\sigma\geq0$ and $(G_i,G_i^\prime)_{i=1}^N$ are $2N$ i.i.d. standard mean zero random vectors in $\bR^d$. The quadratic assignment problem for this two noisy cloud of points is a solution to the problem 
\begin{equation*}
\hat Z\in\argmin_{Z\in \cC} \inr{A, Z} \mbox{ where } A = (\norm{X_i-Y_j}_2^2)_{1\leq i, j\leq N}
\end{equation*} and it can be shown that in the free noise case, we have $Z^*\in\argmin_{Z\in \cC} \inr{\bE A, Z}$. The noisy quadratic OT problem is to identify a sharp phase transition that is a $\sigma^*$ such that 1) if $\sigma<\sigma^*$ then with high probability $\hat Z = Z^*$ and 2) for all $\sigma > \sigma^*$, with probability larger than $1/2$,  $\hat Z \neq Z^*$. Once again, one may looked at $\hat Z$ as an ERM for a linear loss function.

\paragraph{The sparse single index model.}
For this last example, we consider a semi-parametric model where an output $Y\in \bR$ is generated from an input $X\in \bR^d$, via a \simplequote{link} function in the following way:
\[
Y = f\left(\inr{X, \beta^*}\right) + \epsilon
\]where $\beta^*\in\bR^d$ is assumed to be a $k$-sparse unit vector, $f:\bR\to\bR$ is an unknown univariate measurable function and $\epsilon$ is a noise that is generally assumed to be independent of the input. The entries of $X$ are assumed to be $i.i.d$ with a given density $p_0$. The joint density of $X$ is then $p = \otimes_{j = 1}^d p_0$ with respect to the Lebesgue measure. We define a univariate score function $s:x\in\bR\to\bR$ by $s(x)=-p_0'(x)/p_0(x)$, defined for $p_0$-almost all $x\in\bR$ and the first and second score functions associated with $p$ are defined for $p$-almost all $x=(x_j)_{j=1}^d$ by
\[
S(x) = \left(s(x_j)\right)_{1\leq j \leq d} \in\bR^d \mbox{ and } T(x) = S(x)S(x)^\top - \mathrm{diag}\left(\left(s^\prime(x_j)\right)_{1\leq j \leq d}\right).
\]

Unlike the previous examples, the dimension $d$ may be larger than $N$ however, the target index $\beta^*$ is assumed to be $k$-sparse with $k<N$. We therefore fall into the realm of structural learning. The work of \cite{yang2018steins} focuses on this problem where it is proved that $\beta^*$ can be obtained as the leading eigenvector of 
\[
Z^* \in \argmin_{Z\in\cC} \left(- \inr{\bE[A], Z}\right) \mbox{ where } \cC:=\left\{0\preceq W\preceq I_d, ~\mathrm{Tr}(W)=1\right\}
\]
and $A := YT(X)$. Regularized ERM promotes the sparsity structure via a $\ell_1$-regularization. The oracle $Z^*$ can then be estimated as follows:
\[
\hat Z \in \argmin_{Z\in\cC}\left( - \inr{A, Z} + \lambda\norm{Z}_1\right)
\]
which takes the form of a regularized ERM estimator based on the observation $A$, the linear loss function $Z\to \ell_Z(A) = -\inr{A, Z}$ and a $\ell_1$ regularization.

\paragraph{Goal of the paper.} The list of examples provided above indicates that there is a real interest in the general study of linear loss functions in machine learning (we will provide later one more examples in structural learning for which we will provide a complete statistical analysis). Our aim is to propose such a unified methodology to obtain statistical properties of classical machine learning procedures based on  linear loss functions such as the SDP procedures introduced above that we are now looking as ERM procedures constructed with a linear loss function. We continue the work begun in \cite{MR4329809} and go further here by presenting three other estimators that address the two problems of structural risk minimization and robustness. Our machine learning viewpoint allows to introduce new procedures (addressing the previously mentioned two issues) as well as study their statistical properties. 

\paragraph{Framework.} Our general framework is as follows. Let $H$ be a Hilbert space. Let $A$ be a random vector in $H$ that we observe and $\cC\subset H$ be a constraint set (most of the time it will be a  convex set). We suppose to be interested in an object which is the solution to the \simplequote{oracle} optimization problem
\begin{align}\label{eq:def_oracle}
    Z^*\in \argmax_{Z\in\cC}\inr{\bE[A], Z}.
\end{align}
 In some cases, $Z^*$ is not our direct object of interest, but knowing about it enables us to achieve our objective (for instance, by retrieving one of its first eigen-vector). We then propose several estimators for the estimation of the oracle $Z^*$, among which we will choose depending on the presence or not of some particular structure and on the quality of the data (presence or not of corrupted data). 

The first estimator we propose is the one studied in \cite{MR4329809} and is the standard ERM estimator built on the random matrix $A$ but for the (non standard) linear loss function, that is $Z\to\ell_Z(A) = -\inr{A, Z}$:
\begin{align}\label{eq:def_natural_estimator}
    \hat{Z}\in \argmax_{Z\in\cC}\inr{A, Z}.
\end{align}
Then, we turn to two classical machine learning and statistics problems: structured learning and robustness. Leveraging on our view point (i.e. all the previous procedures are all ERMs), we attack the structural learning problem by proposing a regularized version of this ERM estimator by adding a regularization function to the objective function in \eqref{eq:def_natural_estimator}. Afterwards, we turn to the robustness problem and introduce an estimator based on the median of means (MOM) principle, which has been introduced in  \cite{MR4102681} and that is called the minmax MOM. This latter estimator addresses the problem of robustness and can be constructed whatever the loss function is and in particular it fits our linear loss function setup. We show that the resulting estimators  are robust to data contamination as well as to heavy-tailed data. As for ERMs, we present a classical and a regularized version of the minmax MOM estimator in this setup.

For each of those estimators we are able to propose statistical guarantees when $\bE[A]$ is only partially observed through $A$. In particular, our approach leads to new non-asymptotic rates of convergence or exact reconstruction properties for a wide range of estimators that fall within our framework. Then, in order to show the versatility of our approach, we apply these general bounds to the sparse PCA problem. Using our approach we are able to handle the this classical statistical problem using our general excess risk and estimation bounds. As a result we improve the state-of-the art results in sparse PCA as well as introduce new procedures with statistical optimal guarantees that solve the problems of robust structural learning for this problem. Efficient robust gradient descent based  algorithms may easily be derived from these procedures as in \cite{MR4102681}; we will however not dive deeper into the algorithmic consequences of our approach.


\section{General excess risk on estimation bounds for ERM, minmax MOM estimators and their regularized versions}\label{sec:general_bounds}
In this section, we provide high probability excess risk and estimation bounds satisfied by four procedures (ERM, minmax MOM and their regularized versions) in the setup introduced above, that is for the linear loss function. The results for ERM are taken from \cite{MR4329809} and are recalled here for completeness and because it presents an \simplequote{easy} setup for the introduction of two key tools: local complexity fixed points and local curvature equations. The proofs of all the results are postponed to Section~\ref{sec:proofs}. They use state-of-the art machinery such as localization, homogeneity argument, local curvature and fixed point complexity parameters. 

In particular, there are several ways to localize around the oracle depending on the metric used; it can be either the excess risk itself or a natural local curvature metric, denoted later by the $G$ function or the standard $L_2$ metric with respect to the probability measure of the data. Depending on the metric, this defines different local curvatures and different fixed points. For each type of localization, we state a statistical result. We therefore obtain various bounds for each of the four estimators in this section. Hence, this section  provides a complete description of the results one can obtain for these estimators in the setup of linear loss functions and for any regularization norm. We will apply these results in the sparse PCA framework later to show how these general bounds can be applied in a concrete example.

\subsection{General framework}\label{subsec:General_framework}

Throughout this section, we place ourselves in the classical context considered in machine learning and provide its relation with the setup from the Introduction section, in particular, we provide for each example the random matrix $A$ appearing in (\ref{eq:def_oracle}) and (\ref{eq:def_natural_estimator}).

Let $H$ be a Hilbert space and $X$ be a random vector with values in $H$ distributed according to a distribution $P$. For any function $g:H \rightarrow \bR $ for which it makes sense, we denote by $Pg := \bE_{X\sim P}[g(X)]$ the expectation of the $g$ function  under the distribution $P$. For each $p\geq 1$, we denote by $\norm{g}_{L_p} = (P[|g|^p])^\frac{1}{p}$ its $L_p(P)$-norm. Let $\cC$ be a subset of $H$. For all $Z$ in $H$, the loss function of $Z$ is the \textit{linear loss function}, $\ell_Z  : X\in H \rightarrow -\inr{X, Z}$ (it is an alignment measure, which quantifies the error made when estimating $Z$ with $X$). As usual in machine learning, we are interested in the best element in $H$ that minimizes the risk (i.e. the expectation of the loss function) over $\cC$, i.e. we want to estimate/learn/infer/test
\begin{align}\label{eq:general_oracle}
Z^* \in \argmin_{Z\in \cC} P \ell_Z.
\end{align}Sometimes $Z^*$ is called the oracle because it is a quantity we would like to know but we usually cannot have a direct access to it because the distribution $P$ of $X$ is not known to the Statistician and so is the risk function $Z\to P\ell_Z$. However, we have access to a sample distributed according to $P$. This sample / dataset is denoted by $\{X_i:i\in[N]\}$ where $N\in\bN$ is called the sample size. From a mathematical point of view  $(X_i)_{i\in[N]}$ is a family of i.i.d. random variables distributed according to $P$ -- in the section below concerning median-of-means estimators we will relax this assumption and consider a situation where a fraction of the dataset may have been corrupted by an adversary, in that case the $X_i$'s are not anymore assumed to be i.i.d..

The setup we just introduced is pretty much the same as in the Introductory section. We just have to identify the random matrix $A$ for each particular examples. Since, the 'linear loss function' setup is not standard in machine learning, we provide the  connection between $A$ and the $X_i$'s for each example:
\begin{itemize}
    \item in community detection, $N=1$ and $A=X_1$ is the adjacency matrix of the observed graph;
    \item in variable clustering, $A:=\frac{1}{N}\sum_{i = 1}^NX_iX_i^\top$ is the empirical covariance of the observed variables $X_i$'s;
    \item in angular synchronization, $A = \left(e^{\iota\delta_{ij}} + \epsilon_{ij}\right)_{1\leq i < j \leq d}$ is made of the noisy measurements of the pairwise offsets;
    \item in the \textsc{Max-Cut} problem, $A$ is the adjacency matrix of the observed graph;
    \item in phase recovery, $A := \left(XX^\top\right) \circ \left(I_N - BB^+\right)$, where $X$ is the vector of the $N$ observed measurements and $B$ is the measurement matrix;
   \item in distance metric learning, $A := \sum_{i, j = 1}^N\left(X_i-X_j\right)\left(X_i-X_j\right)^\top$ where the $X_i$'s are the observed data from which we want to learn the metric; 
       \item in noisy optimal transport, $A = (\norm{X_i-Y_j}_2^2)_{1\leq i, j\leq N}$, where $\left\{X_1, \ldots, X_N\right\}$ and $\left\{Y_1, \ldots, Y_N\right\}$ are the two sets of observed vectors that we wish to transport one over the other;
 
        \item in the sparse single index model, $A = \frac{1}{N}\sum_{i=1}^N Y_iT(X_i)$, where for any $i\in[N]$, $Y_i = f(\inr{X_i, \beta^*}) + \epsilon_i$ is the noisy output associated to the input $X_i$ via the link function $f$, and $T(X_i)\in\bR^{d\times d}$ is the second order score matrix of $X$.
\end{itemize}

\begin{Remark}
Most of the problems introduced in Section \ref{sec:intro} are presented as maximization problems, whereas ERM is a minimization problem. Given the linearity of the loss function, there are several ways to write the maximization problem into a minimization one: one may take the opposite of the linear loss function, or replace $A$ with $-A$, or $\cC$ with $-\cC$. Here, we consider the loss function $\ell_Z :A \to -\inr{A, Z}$, i.e. we take the opposite of the loss function, which is still a linear one.
\end{Remark}

Moving back to the \doublequote{learning with a linear loss function} introduced at the beginning of this section, we want to estimate/learn the oracle $Z^*$ from the data $(X_i)_{i\in[N]}$. Let $\hat{Z}$ be an estimator  constructed with these data. The quality of prediction of $\hat Z$ is measured via the excess risk $P\cL_{\hat Z}$ where $Z\in \cC\to \cL_Z := \ell_{Z} - \ell_{Z^*}$ is called the excess loss. The quality of estimation of $\hat Z$ is measured by the error rate $\norm{\hat{Z}-Z^*}_{L_2}^2$, where $L_2$ is taken with respect to the $P$ distribution.

There are many ways to construct estimators in the machine learning context considered here. We will see four of them below. The most classical one is the empirical risk minimization procedure \cite{MR1719582} introduced in the next section. Before moving to the construction of estimators, we say a word about the set $\cC$. In all examples introduced in Section~\ref{sec:intro}, $\cC$ is a convex set because of algorithmic considerations. For our theoretical purpose, we will however need a weaker assumption given now: the star-shapped property.

\begin{Definition}\label{def:star-shaped}
We say that a set $\cC$ is star-shaped in $Z^*$ when for all $Z\in \cC$, the segment $[Z,Z^*]$ is in $\cC$.
\end{Definition}

In all our results we will assume $\cC$ to be star-shaped in $Z^*$. This property is satisfied in all examples introduced in Section~\ref{sec:intro} because a convex set is star-shaped in any of its elements. 

\subsection{The ERM estimator and its regularized version: definition and general bounds}
In this section, we consider the \simplequote{$i.i.d$ setup} introduced in the previous section and consider the standard ERM estimator and its regularized version for which we provide high probability excess risk and estimation bounds. The bounds for the ERM are taken from \cite{MR4329809}. We reproduce them here because they introduce key quantities (localization, local curvature and complexity fixed points) in an 'easy' setup and they will appear in the study of the three other estimators in a more convoluted way.

\subsubsection{ERM for the linear loss function}\label{sec:ERM}
For any loss function and in particular for the linear one considered here $\ell_Z:X\in H\to -\inr{X,Z}$, defined for all $Z\in\cC$, the ERM is 
\begin{align*}{}
\hat{Z} \in \argmin_{Z\in \cC} P_N\ell_Z \mbox{ where } P_N \ell_Z = \frac{1}{N}\sum_{i=1}^N \ell_Z(X_i) = \frac{1}{N}\sum_{i=1}^N \inr{-X_i, Z}. 
\end{align*}The ERM is the natural empirical version of the oracle $Z^*$ since $P \ell_Z$ appearing in the definition of $Z^*$ in \eqref{eq:general_oracle} has been replaced by its empirical counterpart $P_N \ell_Z$. When there is only one observation, ie $N=1$, for instance in the community detection problem, we simply have $P_N \ell_Z = P_1\ell_Z = -\inr{X_1, Z} = -\inr{A, Z}$. 

The study of the statistical properties of ERM estimators goes back to \cite{MR0474638} and has been at the heart of many researches since then \cite{Kol11,MR2182250}. The results recalled below are for the special case of the linear loss function and are taken from \cite{MR4329809}. They are however based on nowadays classical concepts in machine learning.

A key quantity driving the rate of convergence of the ERM is a local complexity fixed point parameter. This kind of parameter carries all the statistical complexity of the problem. It can however be hard to compute (see for instance \cite{MR4329809} or  Section~\ref{sec:processus_sto} below), since it requires to control with large probability the supremum of an empirical processes indexed by a "localized classes", i.e. of the set $\cC$ intersected with a neighborhood of the oracle. We now define such a complexity fixed point related to the problem we are considering here. 

  \begin{Definition}
    \label{def:fixed_point}[Complexity fixed point parameter] 
    Let $0<\Delta<1$. The fixed point complexity parameter at deviation $1-\Delta$ is
    \begin{equation}\label{eq:fixed_point}
    r^*(\Delta) = \inf\left(r>0 : \bP\left[\sup_{Z\in\cC: P\cL_Z \leq r}(P_N-P)\cL_Z\leq \frac{r}{2}\right]\geq 1-\Delta\right).
    \end{equation}
  \end{Definition}

In what follows, we give some statistical properties of the ERM $\hat{Z}$ build from this complexity parameter. They are taken from  \cite{MR4329809} even though they have been obtained for the special case $N=1$ and $X_1=A$ they can be easily extended to the setup considered here for a general sample size $N$.

 \begin{Theorem}[Theorem~1 in \cite{MR4329809}]\label{theo:main_ERM}
    We assume that the constraint $\cC$ is star-shaped in $Z^*$. Then, for all $0<\Delta<1$, with probability at least $1-\Delta$, it holds true that $P\cL_{\hat{Z}}\leq r^*(\Delta)$.
  \end{Theorem}

From Theorem \ref{theo:main_ERM}, we get that one way to grab some information on the ERM is to get an upper bound for the complexity fixed point $r^*(\Delta)$. To that end, one needs to understand the shape of the sets $\cC\cap\{Z:P\cL_{Z}\leq r\}$ for $r>0$. This task may however be hard because of the shape of the neighborhoods $\{Z:P\cL_{Z}\leq r\}$ given by the excess risk. In that case, it has been shown \cite{chinot2018statistical} that one can leverage on a \textit{local curvature}  of the excess risk to introduce easier to compute fixed points. We are now introducing the complexity fixed point associated with this other localization and then the notion of local curvature. In what follows, $G$ is some function from $H$ to $\bR$. 
 \begin{Definition}
    \label{def:fixed_point_Gloc}[Complexity fixed point parameter with G-localization]
    Let $0<\Delta<1$. The fixed point complexity parameter with respect to the $G$-localization at deviation $1-\Delta$ is
    \begin{equation}\label{eq:fixed_point}
    r_G^*(\Delta) = \inf\left(r>0 : \bP\left[\sup_{Z\in\cC: G(Z-Z^*)\leq r}(P_1-P)\cL_Z\leq \frac{r}{2}\right]\geq 1-\Delta\right).
    \end{equation}
  \end{Definition}

 The difference between $r^*$ and $r_G^*$ lies in the fact that the local subsets are not defined with the same proximity function: $r^*$ used the excess risk function for localization whereas $r_G^*$ uses the $G$ function. The latter $G$ function should play the role of a simple description of the curvature of the excess risk around the oracle as it is granted in the following assumption.  
 \begin{Assumption}\label{ass:curvature_G}
    For all $Z\in\cC$, if $P\cL_Z\leq r^*_G(\Delta)$ then $P\cL_Z\geq G(Z^*-Z)$.
 \end{Assumption} 

 There are examples where one can show a curvature of the excess risk over the entire set $\cC$ - this is for instance the case in the sparse PCA example below (see Lemma~\ref{coro:curvature_excess_risk} below). In that case, we speak about a \textit{global} curvature. What shows the following result is that we only need a \textit{local} curvature of the excess risk around $Z^*$ to hold in order to get statistical bounds for the ERM $\hat Z$.

 \begin{Theorem}[Corollary~1 in  \cite{MR4329809}]\label{coro:main_coro_ERM}
    We assume that the constraint $\cC$ is star-shaped in $Z^*$ and that the ``local curvature'' Assumption~\ref{ass:curvature_G} holds  for some $0<\Delta<1$. With probability at least $1-\Delta$, it holds true that
    \begin{equation*}
    r^*_G(\Delta)\geq P\cL_{\hat Z}\geq G(Z^*-\hat Z).
    \end{equation*}
 \end{Theorem}

 Finally a third and final estimation bound is given in the following for cases where Assumption \ref{ass:curvature_G} is hard to verify. They are situations where the shape of the local subsets $\cC\cap\{Z:P\cL_Z\leq r\}$ is hard to understand. In that case, we can simplify this assumption by considering neighborhoods with respect to the $G$ function. 

\begin{Assumption}\label{ass:curvature_G_2}
For all $Z\in \cC$, if $G(Z^*-Z)\leq r_G^*(\Delta)$, then $P\cL_Z\geq G(Z^*-Z)$.
\end{Assumption}
The following result establishes that, under Assumption~\ref{ass:curvature_G_2}, $\hat Z$ is a good estimate of $Z^*$ with respect to the $G$ function, but no guarantee on the excess risk is obtained.
\begin{Theorem}[Theorem~2 in \cite{MR4329809}]\label{theo:main_2}
    We assume that the constraint $\cC$ is star-shaped in $Z^*$ and that the ``local curvature'' Assumption~\ref{ass:curvature_G_2} holds  for some $0<\Delta<1$. We assume that the $G$ function is continuous, $G(0)=0$ and $G(\lambda (Z^*-Z))\leq \lambda G(Z^*-Z)$ for any $\lambda\in[0,1]$ and $Z\in \cC$. Then, with probability at least $1-\Delta$, it holds true that $G(Z^*-\hat Z)\leq r^*_G(\Delta)$.
  \end{Theorem}

We refer the reader to \cite{Lucie1} for the application of these results in community detection, signed clustering, angular group synchronization (for both multiplicative and additive models) and the \textsc{Max-Cut} problem. All these problems share the feature that the oracle $Z^*$ does not have some special structure onto which one can leverage to improve the rates of convergence. They are however situations such as in sparse PCA or in the sparse single index model where the target has a structure that can be beneficial in order to improve statistical performance. In such cases, one may consider some regularization procedures like in the following section. 

\subsubsection{Regularized ERM for the  linear loss}\label{sec:RERM}
We focus here on structural learning in which targets/oracles have a structure (such as sparsity, low rank or regularity) onto which the statistician can leverage to construct more statistically efficient estimators. The typical approach to this problem is to regularize the  ERM in order to force the estimator toward the desired structure.  

We place ourselves in the framework defined above in Section~\ref{subsec:General_framework} except that we need here a regularization function, i.e. a function that favors some structure. In this work, we consider a general norm defined at least on the span of $\cC$ and denoted by $\norm{\cdot}$. Typical examples are the $\ell_1$ norm and the trace-norm used in high-dimensional statistics to induce sparsity or low-rank. When $Z^*$ has some structure a natural way to force an estimator toward $Z^*$ is by adding a mutliple of this norm. This yields to the regularized ERM, later called RERM:
\begin{align}\label{def:RERM_estimator}
\hat{Z}^{\mathrm{RERM}} \in \argmin_{Z\in \cC} \left( P_N \ell_Z + \lambda \norm{Z}\right)
\end{align}
where $\lambda>0$ is called the regularization parameter and has the role to  make a trade-off between the data adequation term $P_N \ell_Z$ and the regularization term $\norm{Z}$.

As for the ERM, convergence rates achieved by the RERM $\hat{Z}^{\mathrm{RERM}}$ are driven by a local complexity fixed point parameter. However, the regularization norm appears in this type of parameter: it is now the set $\cC$ intersected with balls with respect to $\norm{\cdot}$ centered at $Z^*$ (and for some radius) that are ``localized'' by some neighborhood of $Z^*$. Somehow the model in structural learning is of the form $\cC\cap \{Z:\norm{Z-Z^*}\leq r\}$.  As in the ERM case, one may consider two different ways to construct localization: either via the excess risk or via a local curvature $G$ function.  However, to avoid a lengthy presentation,  we  focus only on the latter one, i.e. on a localization via a local curvature $G$ function because it is this result that we will use for the our application later in sparse PCA. In what follows, we  consider a function $G:H\rightarrow \bR$, which characterizes the curvature of the objective function, i.e. the risk, $Z\in H\rightarrow P \ell_Z$ around its minimizer $Z^*$. 

\begin{Definition}\label{def:RERM_Complexity_G_function}
For parameter $A>0$, radius $\rho > 0$ and deviation parameter $\delta \in (0, 1)$, we define the complexity fixed point for the structural learning with a linear loss function by
\begin{align*}
r^*_{\mathrm{RERM, G}}(A, \rho, \delta) = \inf\left(r > 0: \bP\left(\sup_{Z\in \cC : \norm{Z-Z^*}\leq \rho, G(Z-Z^*)\leq r}|(P-P_N)\cL_Z| \leq \frac{r}{3A}\right)\geq 1-\delta\right),
\end{align*}where we recall that for all  $Z\in\cC$, $\cL_Z=\ell_Z-\ell_{Z^*}$ is the excess loss function of $Z$.
\end{Definition}

After introducing the fixed point $r^*_{\mathrm{RERM, G}}(A, \rho, \delta)$, we are now in a position to introduce the $G$ function. As we already mentioned above, the $G$ function describes the curvature of the excess risk locally around the oracle. 

\begin{Assumption}\label{ass:A-spars}
We assume there exist $A>0$, $\rho^*>0$ and $\delta \in (0, 1)$ such that, for all $ Z \in \cC$ satisfying $G(Z-Z^*) = r^*_{\mathrm{RERM, G}}(A, \rho^*, \delta)$ and $\norm{Z-Z^*}\leq \rho^*$, then $AP\cL_Z\geq G(Z-Z^*) $.
\end{Assumption}

We now leverage on the structure inducing property of the regularization norm and explain what features must the radius $\rho^*$ appearing in Assumption~\ref{ass:A-spars} have in relation to this property. We will use the assumption below, that is adapted from the one in \cite{lecue2017regularization}, to get the statistical bounds satisfied by the RERM estimator $\hat Z^{RERM}$. The idea is that the regularization norm $\norm{.}$ is expected to promote some structure by having a large subdifferential at elements in $H$ having this structure. First, let us recall what the subdifferential of $\norm{.}$ at a point $Z$ is:
\begin{align*}
(\partial\norm{.})_Z := \left\{\Phi \in H:\norm{Z+h}-\norm{Z}\geq \inr{\Phi, h} \mbox{ for all } h \in H\right\}.
\end{align*}Elements in $(\partial\norm{.})_Z$ are called the \textit{subgradients} of $\norm{\cdot}$ in $Z$.  What matters in structural learning to get fast rates is that $Z^*$ is close to an element with a structure induced by the regularization norm. Therefore we consider the set of all subgradients of $\norm{\cdot}$ of points close to $Z^*$: 
\begin{equation*}
\mbox{for any } \rho>0 \, : \quad \Gamma_{Z^*}(\rho) = \underset{Z\in Z^*+\frac{\rho}{20}B}{\bigcup}(\partial\norm{.})_Z
\end{equation*}where $B$ is the unit ball of $\norm{.}$. We expect $\Gamma_{Z^*}(\rho)$ to be a large subset of the unit dual sphere (or dual ball, when $0\in Z^*+(\rho/20)B$) of $\norm{.}$ when $Z^*$ is structured or close to a structured element in $H$, for the notion of structure associated with $\norm{.}$. This intuition is formalized in the following definition. 

\begin{Definition}\label{def:RERM_A-sparsity} For $A>0$, $\rho > 0$ and $\delta\in(0, 1)$ we define:
\begin{align*}
H_{\rho, A} := \left\{Z\in \cC: \norm{Z-Z^*} = \rho ~~ and ~~ G(Z-Z^*)\leq r_{\mathrm{RERM}, G}^*(A, \rho, \delta)\right\}
\end{align*}
and
\begin{align*}
\Delta(\rho, A) := \inf_{Z\in H_{\rho, A} } \sup_{\Phi \in \Gamma_{Z^*}(\rho)} \inr{\Phi, Z-Z^*}.
\end{align*}
We say that $\rho > 0$ satisfies the \textbf{$A$-sparsity equation} when $\Delta(\rho, A)\geq (4/5)\rho$.
\end{Definition}
Note that it is always true that $\Delta(\rho, A)\leq\rho$ -- because $\norm{Z-Z^*}=\rho$ and $\Phi$ is a subgradient of $\norm{\cdot}$ -- hence, a radius $\rho$ satisfying the $A$-sparsity equation is somehow extremal up to the absolute constant $4/5$ (the analysis works for any other absolute constant, there is nothing special with $4/5$). It means that $\Gamma_{Z^*}(\rho)$ is almost as big as the unit dual sphere (or ball) of $\norm{\cdot}$. \\ 

All the material introduced above (complexity fixed points, local curvatures and the sparsity equation) are the corner stones of our statistical  analysis of RERMs. Once introduced,  we are in a position to state our main result on RERM estimators for linear loss functions and a general regularization norm. 

\begin{Theorem}\label{theo:main-penalized}
Let $\delta\in (0, 1)$. Assume that the constraint set $\cC$ is star-shaped in $Z^*$. Consider a continuous function $G:H\rightarrow \bR$ such that $G(0)=0$. Suppose the existence of $A>0$ and $\rho^*>0$ such that that Assumption \ref{ass:A-spars} holds and $\rho^*>0$ satisfies the $A$-sparsity equation from Definition \ref{def:RERM_A-sparsity}. Define the function $r^*(.):= r^*_{\mathrm{RERM, G}}(A, ., \delta)$ and assume that
\begin{align}\label{hyp:lambda} 
\frac{10}{21A}\frac{r^*(\rho^*)}{\rho^*} < \lambda < \frac{2}{3A}\frac{r^*(\rho^*)}{\rho^*}.
\end{align}
Then, with probability at least $1-\delta$, the following bounds hold for the RERM estimator defined in (\ref{def:RERM_estimator}):
\begin{align*}
\norm{\hat{Z}^{\mathrm{RERM}}-Z^*} \leq \rho^* ~~,~~ G(\hat{Z}^{\mathrm{RERM}}-Z^*) \leq r^*(\rho^*) ~~and ~~P\cL_{\hat{Z}^{\mathrm{RERM}}} \leq \frac{r^*(\rho^*)}{A}.
\end{align*}
\end{Theorem}

We note that in the case where $G$ is the risk function $Z\rightarrow P\ell_Z$ - that is when the excess risk is used for localization, because, by linearity $G(Z-Z^*)=P\ell_{Z-Z^*}=P \cL_Z$ - Assumption \ref{ass:A-spars} is trivially verified with $A=1$, and as a consequence Theorem \ref{theo:main-penalized} applies.

\subsection{Median of Means estimators: definitions and general bounds}
In this section, we move to the construction and the statistical analysis of another family of estimators introduced in \cite{MR4102681} whose aims are to solve robustness issues related to adversarial contamination of the dataset as well as heavy-tailed data. We are interested here in the case where our data could be contaminated by possible outliers generated by an adversary and the inliers data may be heavy-tailed. Even though the framework seems not in favor of statisticians because the dataset is of poor quality, we still want to achieve the same statistical performance as if there was no outliers and light-tailed (such as sub-gaussian) data.   It is known that the classical ERM or RERM approaches from the previous section do not perform well in general on this type of dataset and that is the reason why we move to MOM estimators.

The statistical framework considered in this section cannot be the ideal i.i.d. setup considered in the previous section that fits well for ERM and RERM. Indeed, the i.i.d. framework do not allow for adversarial corruption. That is why we consider the following setup in this section. 

\begin{Assumption}\label{ass:adv_contamination}[Adversarial contamination setup]
Let $N$ i.i.d. random vectors $(\widetilde{X}_i)_{i=1}^N$ in $H$. These vectors are first given to an adversary who is allowed to modify up to $|\cO|$ of them. This modification does not have to follow any rule and is unknown to the statistician. This leads to the modified dataset $\{X_1, \ldots, X_N\}$ that the adversary gives to the statistician. Hence, the dataset at hands $\{X_1, \ldots, X_N\}$ is said to be \simplequote{adversarially} contaminated. It can be partitioned into two groups: the modified data $(X_i)_{i\in\cO}$, which can be seen as outliers and the \simplequote{good data}, or inliers, $(X_i)_{i\in\cI}$ such that for any $i\in\cI$, $X_i=\widetilde{X}_i$. Of course, the statistician does not know which data has been modified or not so that the partition $\cO\cup\cI = \{1, \ldots, N\}$ is unknown to the statistician.
\end{Assumption}

\begin{Remark}
Since there are two types of data considered in Assumption~\ref{ass:adv_contamination} (the 'good' $\tilde X_i$s and the corrupted ones $X_i$s), we need to be clear on the objects we will be using later: the risk function and its associated oracle are the one associated with the 'good' data:
\begin{equation*}
 Z\in\cC \to P\ell_Z= \bE \inr{-\tilde X, Z} \mbox{ and } Z^*\in\argmin_{Z\in\cC} P\ell_Z
 \end{equation*} where $\tilde X$ has the same probability distribution as $\tilde X_1, \ldots, \tilde X_N$. It is also the same for the $L_2$-norm: for all $Z\in H$, $\norm{Z}_{L_2} = \sqrt{\bE \inr{\tilde X, Z}^2}$. Note that the $L_2$-norm is in general different from the original Hilbert norm defining $H$, which is denoted by $\norm{\cdot}_2$.
\end{Remark}

The adversarial contamination setup addresses several questions in statistics regarding the rates of convergence, the probability deviations and the number of outliers. Many approaches have been introduced to answer these questions \cite{MR2488795}. There was an important renewal of this topic during the last ten years \cite{catoni_challenging_2012,MR3631028}. The approach we use in this section is based on the median-of-means principle \cite{MR702836,MR855970}: $[N]$ is partitioned into $K$ equal-size groups $B_1, \ldots, B_K$ (w.l.o.g. $K$ is assumed to divide $N$, otherwise we only have to remove some data). For any function $g:H\rightarrow \bR $ and $k\in [K]$ we define $P_{B_k}g=(K/N)\sum_{i\in B_k}g(X_i)$, the empirical mean of $g$ over $B_k$. Then, we define $\mathrm{MOM}_k(g)$ as the median of these $K$ empirical means:
\begin{align*}
\mathrm{MOM}_K(g) := \mathrm{Med}(P_{B_1}g, \ldots, P_{B_K}g).
\end{align*}This data partition scheme is at the heart of our approach to answer the robustness issues. It is used as a building block in the minmax MOM estimator. We recall its construction and provide its statistical properties in the remaining of this section as well as for its regularized version for the robust structural learning problem.

\subsubsection{The minmax MOM estimator for the linear loss function.}\label{sec:MOM}
To solve the robustness to adversarial corruption as well as to heavy-tailed data, one can use a systematic approach called the minimax MOM estimator in \cite{MR4102681}. It works whenever a loss function exists and a robust gradient descent algorithm may also be constructed out of it (see \cite{MR4102681} for more details). When the dataset has been splited into $K$ equal size blocks, it takes the following form:
\begin{align}\label{def:MOM-estimator}
\hat{Z}^{\mathrm{MOM}}_K \in \argmin_{Z\in \cC} \sup_{Z^{\prime}\in \cC}\mathrm{MOM}_K (\ell_{Z}-\ell_{Z^{\prime}})
\end{align}and can therefore be used in the particular case studied here of the linear loss function $x\to\ell_Z(x) = -\inr{Z,x}$. From our theoretical perspective, the aim of the minmax MOM estimator $\hat{Z}^{\mathrm{MOM}}_K $ is to achieve the rates of convergence for the same deviation probabilty in the contaminated and heavy-tailed setup as in the ideal i.i.d. setup with light-tailed data, as long as the number of outliers is not too large. It is the aim of the next section to prove such  statistical bounds. As for the ERM case, rates of convergence are given by local complexity fixed points  that depends on the choice of localization. Below, we consider three different ways to localize: either via the $L_2(P)$-norm, or via the excess risk or via some general curvature function $G$.

\paragraph{MOM estimator with excess-risk localization.} As previously for ERMs, the convergence rate of the minmax MOM estimator is driven by a local complexity fixed point parameters. In this section, we  consider the case where the excess risk is simple enough so that it can serve as a localization. In that case, there is no need to identify the curvature of the excess risk locally around $Z^*$ since the excess risk describes it by itself. There is therefore no curvature assumption. In the next two paragraphs the picture will be different.

\begin{Definition}\label{def:MOM-fixed-point-excess-risk}
Let $\sigma_1, \ldots, \sigma_N$ be $N$ independent Rademacher variables which are independent of the $\tilde X_i$'s. For $\gamma>0$, we define:
\begin{align*}
r^*_{\mathrm{MOM}, \mathrm{ER}}(\gamma) := \mathrm{inf}\left\{r>0~:~\mathrm{max}\left(\frac{\mathrm{E}(r)}{\gamma},\sqrt{12800}V_K(r)
\right)\leq r^2\right\} 
\end{align*}
where, for all $r>0$,
\begin{align*}
\mathrm{E}(r) := \bE\left[\underset{Z\in\cC:P\cL_Z\leq r^2}{\mathrm{sup}}\left|\frac{1}{N}\sum_{i=1}^N\sigma_i\cL_Z(\tilde X_i) 
\right|\right] 
\mbox{ and }
V_K(r) := \sqrt{\frac{K}{N}}\underset{Z\in\cC:P\cL_Z\leq r^2}{\mathrm{sup}} \sqrt{\bV ar(\cL_Z(\tilde X))}.
\end{align*}
\end{Definition}

In the case of excess risk localization, there is no need for other tools than the fixed point $r^*_{\mathrm{MOM}, \mathrm{ER}}(\gamma)$ to describe the rate of convergence of the minmax MOM. This is what shows the following result.

\begin{Theorem}\label{theo:main_MOM_loc_excess_risk}We consider the adversarial contamination  setup of Assumption~\ref{ass:adv_contamination}.
We assume that the constraint set $\cC$ is star-shaped in $Z^*$. Let $\gamma=1/6400$ and consider $K$, a divisor of $N$ such that $K\geq 100|\cO|$. Then, it holds true that with probability at least $1-\exp(- 72K/625)$, $P\cL_{\hat{Z}^{\mathrm{MOM}}_K} \leq r^*_{\mathrm{MOM}, \mathrm{ER}}(\gamma)^2$. 
\end{Theorem}

Compared to the fixed point from Definition~\ref{def:fixed_point} describing the rate of convergence of the ERM, we note that the one from Definition~\ref{def:MOM-fixed-point-excess-risk} uses a local Rademacher complexity, denoted by $E(r)$ , and a variance term, denoted by $V_K(r)$; there is no need to upper bound with high probability the supremum of an empirical process. For minmax MOM estimators, the task of computing fixed point complexity parameters is therefore easier. Moreover, as one can see in Theorem~\ref{theo:main_MOM_loc_excess_risk}, the convergence rate is obtained with an exponentially large probability even though no strong concentration property is assumed; only the existence of a second moment (so that the variance term $V_K(r)$ exists) is required. This shows the robustness to heavy-tail data of minmax MOM estimators for the linear loss function as well as its robustness with respect to adversarial contamination since it is proved in the setup of Assumption~\ref{ass:adv_contamination}. However, the computation of the complexity term $E(r)$ may require more moments than just $2$ in order to recover a Gaussian regime, i.e. a rate achieved when the data have a light subgaussian tail.

\paragraph{MOM estimator with $L_2$-localization.} In this section, we  consider the case where the behaviour / curvature of the excess risk locally around the oracle $Z^*$ is well described by the $L_2$-norm to the square. This is the situation when a margin assumption $ A P\cL_Z\geq \norm{Z-Z^*}_{L_2}^2, \forall Z\in\cC$ holds, i.e. with a margin parameter equal to $2$ \cite{MR1765618}. In that case, one needs to modify the definition of the complexity fixed point parameter by using a $L_2$-localization. 

\begin{Definition}\label{def:complexity_param_MOM}
Let $\sigma_1, \ldots, \sigma_N$ be independent Rademacher variables which are independent of the $\tilde X_i$'s. For $\gamma>0$, we define
\begin{align*}
r^*_{\mathrm{MOM}, L_2}(\gamma) := \inf \left(r>0 :\bE\left[\sup_{Z\in \cC: ~\norm{Z-Z^*}_{L_2}\leq r} \left|\frac{1}{N}\sum_{i=1}^N\sigma_i \cL_{Z}(\tilde X_i)\right|\right] \leq \gamma r^2 \right)
\end{align*}where we recall that $\norm{Z}_{L_2} = \sqrt{\bE \inr{\tilde X,Z}^2}$ for all $Z\in H$.
\end{Definition}


As we said above, we use the $L_2$-norm in the localization to define the fixed point $r^*_{\mathrm{MOM}, L_2}(\gamma)$ when it describes  the curvature of the excess risk around $Z^*$. We now formalize this property in the next assumption.

\begin{Assumption}\label{ass:curvature_MOM_L2}
There exists $A>0$ such that for any $Z\in \cC$, if $\norm{Z-Z^*}_{L_2}^2\leq C_{K,A} $, then $\norm{Z-Z^*}_{L_2}^2\leq AP\cL_Z$, where $C_{K,A}:=\max\left(r^*_{\mathrm{MOM}, L_2}(\gamma)^2, \gamma^{-1}A^2(K/N)\right)$ for $\gamma = 1/3200$.
\end{Assumption}


Looking at Assumption~\ref{ass:curvature_MOM_L2}, this may be surprising to have a quadratic term $\norm{Z-Z^*}_{L_2}^2$ describing a linear term $P\cL_Z = \inr{\bE \tilde X, Z^*-Z}$. However, one may see that the local curvature of the excess risk from Assumption~\ref{ass:curvature_MOM_L2} holds only for  $Z$ in $\cC$ not in $H$.  Thanks to the two tools introduced above (a local complexity fixed point and a curvature assumption), we are now ready to state our main result on the minmax MOM estimator in the adversarial contamination setup for a  $L_2$-localization.


\begin{Theorem}\label{theo:main_MOM_unloc} We consider the adversarial contamination setup of Assumption~\ref{ass:adv_contamination}. We assume that the constraint set $\cC$ is star-shaped in $Z^*$. Let $\gamma = 1/3200$. Assume the existence of $0<A<1$ such that Assumption \ref{ass:curvature_MOM_L2} holds. Let $K$ be a divisor of $N$ such that $K\geq 100|\cO|$. Then, it holds true that with probability at least $1-\exp(-  72K/625)$: 
\begin{align*}
P\cL_{\hat{Z}^{\mathrm{MOM}}_K} \leq \frac{C_{K,A}}{A} \mbox{ and } \norm{\hat{Z}^{\mathrm{MOM}}_K-Z^*}_{L_2}^2\leq C_{K,A}.
\end{align*}
\end{Theorem}

Theorem~\ref{theo:main_MOM_unloc} can be used under a margin assumption with a margin parameter equal to $2$. It can be extended to margin parameter other than $2$. However, one may be interested in other situations where the local curvature of the excess risk is not described by the square of the $L_2$ norm but for instance by the square of the native Hilbert norm of $H$ - as it will be the case for the sparse PCA problem. In the next paragraph, we provide a statistical bound for the minmax MOM estimator for a local curvature of the excess risk described by a general $G$ function.

\paragraph{MOM estimator with $G$ localization.} In this final paragraph regarding the minmax MOM estimator, we consider a general $G$ function describing locally the excess risk around $Z^*$ and derive statistical bounds when this function is used for localization.  When applied to the particular cases of the excess risk or the $L_2$ norm to the square, we recover the last two results. However, other $G$ functions may be considered, for instance,  if the calculation of $r^*_{\mathrm{MOM}, \mathrm{ER}}(\gamma)$ is too hard or if $L_2$-norm to the square does not describe well enough the excess risk. We need first to define a complexity fixed point for a localization w.r.t. a general $G$ function. Unlike in the previous section dealing with the $L_2$ to the square localization and as in the last but one section dealing with a excess risk localization, there is a variance term in this fixed point equation.

\begin{Definition}\label{def:local-fixed-point-MOM}
Let $\sigma_1, \ldots, \sigma_N$ be $N$ independent Rademacher variables which are independent of the $\tilde X_i$'s. For $G:H\rightarrow \bR$ and $\gamma>0$, we define:
\begin{align*}
r^*_{\mathrm{MOM, G}}(\gamma) := \mathrm{inf}\left\{r>0~:~\mathrm{max}\left(\frac{\mathrm{E}_G(r)}{\gamma},\sqrt{12800}V_{K,G}(r)
\right)\leq r^2\right\} 
\end{align*}
where, for all $r>0$,
\begin{align*}
\mathrm{E}_G(r) := \bE\left[\underset{Z\in\cC:G(Z-Z^*)\leq r^2}{\mathrm{sup}}\left|\frac{1}{N}\sum_{i=1}^N\sigma_i\cL_Z(\tilde X_i)
\right|\right] 
\mbox{ and }
V_{K,G}(r) := \sqrt{\frac{K}{N}}\underset{Z\in\cC:G(Z-Z^*)\leq r^2}{\mathrm{sup}} \sqrt{\bV ar(\cL_Z(\tilde X))}.
\end{align*}
\end{Definition}

The function $G$  characterizes the curvature of the excess risk $Z\in\cC\rightarrow P\cL_Z = \inr{\bE X, Z^*-Z}$ locally around its minimizer $Z^*$. This is formalized in the following assumption. 

\begin{Assumption}\label{ass:curvature-MOM}
There exist $A>0$ and $\gamma>0$ such that for all $Z\in \cC$, if $G(Z-Z^*) \leq (r^*_{\mathrm{MOM, G}}(\gamma))^2$, then $AP\cL_Z \geq G(Z-Z^*)$.
\end{Assumption}

The difference between $r^*_{\mathrm{MOM, ER}}$ and $r^*_{\mathrm{MOM, G}}$ is that the local subsets are not defined using the same proximity function to the oracle $Z^*$. The main advantage in finding a curvature function $G$ satisfying Assumption~\ref{ass:curvature-MOM} is that $r^*_{\mathrm{MOM, G}}$ may be easier to compute than $r^*_{\mathrm{MOM, ER}}$, since the shape of a neighborhood defined by $G$ may be easier to understand than the one defined by the excess risk. However, one always has $r^*_{\mathrm{MOM, ER}}\leq r^*_{\mathrm{ERM, G}}$ since there is no better way to describe the excess risk than the excess risk itself. We now obtain statistical bounds satisfied by the minmax MOM estimator \eqref{def:MOM-estimator} under this local curvature assumption.

\begin{Theorem}\label{coro:main-MOM-loc} We consider the adversarial contamination setup of Assumption~\ref{ass:adv_contamination}.
We assume that the constraint set $\cC$ is star-shaped in $Z^*$. We consider a continuous function $G:H\rightarrow \bR$ be a continuous function. Let $\gamma = 1/6400$. We assume the existence of $0<A<2$ such that the local curvature Assumption \ref{ass:curvature-MOM} holds for those values of $\gamma$ and $G$. Then, with probability at least $1-\exp(-72K/625)$ it holds true that:
\begin{align*}
P\cL_{\hat{Z}^{\mathrm{MOM}}_K} \leq \frac{1}{2}r^*_{\mathrm{MOM, G}}(\gamma)^2 \quad \mbox{ and } \quad G(Z^*-\hat{Z}^{\mathrm{MOM}}_K) \leq r^*_{\mathrm{MOM, G}}(\gamma)^2.
\end{align*}
\end{Theorem}

Theorem~\ref{coro:main-MOM-loc} may be applied in the examples introduced from Section~\ref{sec:intro} if one is willing to handle robustness issues for these (none structured) learning problems. If one wants to handle the robustness issues in structural learning then one may consider regularized versions of the minmax MOM estimator.

\subsubsection{Regularized minmax MOM estimators for the linear loss function}\label{sec:RMOM}
We are now considering the setup of structural learning that allows for high-dimensional statistics, i.e. when the dimension of the parameter to estimate $Z^*$ is larger than the number of observations. In that case, some structure is usually assumed to be satisfied by $Z^*$ and should be taken into account for the construction of estimators. On top of that, we consider a setup where the data may have been corrupted by some outliers and the inliers may be heavy-tailed. We therefore have to face several issues related to robustness and high-dimensions that we propose to solve using a regularized version of the minmax MOM estimator introduced in Section~\ref{sec:MOM}:
\begin{align}\label{def:RMOM-estimator}
\hat{Z}^{\mathrm{RMOM}}_{K, \lambda} \in \argmin_{Z\in \cC} \sup_{Z^\prime \in \cC} \left(\mathrm{MOM}_K (\ell_Z - \ell_{Z^\prime}) + \lambda(\norm{Z} - \norm{Z^\prime})\right)
\end{align}
where $\lambda >0$ is some regularization parameter and $\norm{\cdot}$ is a norm inducing some structure. In the following sections, we provide statistical guarantees for this estimator. As in the previous sections, the convergence rates depend on local complexity fixed points, local curvature properties of the excess risk and of the 'structure inducing power' of the regularization norm $\norm{\cdot}$. As previously, the choice of the localization function plays a key role in the definition of all these concepts. We therefore consider three paragraphs depending on the localization function used: it can either be the excess risk, the $L_2$-norm or some general function $G$. 

\paragraph{RMOM estimator with excess-risk localization.}As in the previous section, we start with the excess risk localization. 

\begin{Definition}\label{def:G-er-fixed-point-RMOM}
Let $\sigma_1, \ldots, \sigma_N$ be independent Rademacher variables which are independent of the $\tilde X_i$'s. For $\gamma>0$ and $\rho >0$, we define:
\begin{align*}
r^*_{\mathrm{RMOM, ER}}(\gamma, \rho) := \mathrm{inf}\left\{r>0~:~\mathrm{max}\left(\frac{\mathrm{E}(r, \rho)}{\gamma},400\sqrt{2}V_{K}(r, \rho)
\right)\leq r^2\right\} 
\end{align*}
where, for all $\rho,r>0$ and $\cC_{\rho, r} = \left\{Z\in\cC:\norm{Z-Z^*}\leq\rho,P\cL_Z\leq r^2\right\}$, 
\begin{align*}
\mathrm{E}(r, \rho) := \bE\left[\sup_{Z\in\cC_{\rho, r}}\left|\frac{1}{N}\sum_{i=1}^N\sigma_i\cL_Z(\tilde X_i)
\right|\right] \mbox{ and }
V_{K}(r, \rho) := \sqrt{\frac{K}{N}} \sup_{Z\in\cC_{\rho, r}}\sqrt{\bV ar(\cL_Z(\tilde X))}.
\end{align*}
\end{Definition}
The sparsity equation introduced for the study of the RERM in Definition~\ref{def:RERM_A-sparsity} has to be slightly modified according to this new definition of the complexity parameter. 

\begin{Definition}\label{def:sparsity_RMOM_er} For $\gamma>0$ and $\rho>0$, let
$
\bar{H}_{\rho} := \left\{Z\in\cC:\norm{Z-Z^*}=\rho\mbox{ and }P\cL_Z\leq r^*_{\mathrm{RMOM, ER}}(\gamma, \rho)^2\right\}
$
and
$
\bar{\Delta}(\rho) := \inf_{Z\in \bar{H}_{\rho}} \sup_{\Phi \in \Gamma_{Z^*}(\rho)} \inr{\Phi, Z-Z^*}.
$
We say that $\rho$ satisfies the \textbf{sparsity equation} if $\bar{\Delta}(\rho)\geq4\rho / 5$.  
\end{Definition}

We are now ready to state our main statistical result satisfied by the regularized minmax MOM  estimator for the linear loss function and for an excess-risk localization. 

\begin{Theorem}\label{theo:main_RERM_er}
We consider the adversarial contamination setup of Assumption~\ref{ass:adv_contamination}. Let $K\in[N]$ be such that $K\geq 100|\cO|$. Let $\rho^*>0$  satisfying the sparsity equation from Definition~\ref{def:sparsity_RMOM_er}. Let $\gamma = 1/3200$ and take $\lambda=(11/(40\rho^*))r^*_{\mathrm{RMOM, ER}}(\gamma, 2\rho^*)$ as regularization parameter. Then, with probability at least $1-2\exp(-72K/625)$,
\begin{align*}
P\cL_{\hat{Z}^{\mathrm{RMOM}}_{K, \lambda}} \leq r^*_{\mathrm{RMOM, ER}}(\gamma, 2\rho^*)^2 \mbox{ and } \norm{\hat{Z}^{\mathrm{RMOM}}_{K, \lambda}-Z^*}\leq 2\rho^*.
\end{align*}
\end{Theorem}
Note that one may replace $r^*_{\mathrm{RMOM, ER}}(\gamma, 2\rho^*)$ by any real number $r^*$ larger than $r^*_{\mathrm{RMOM, ER}}(\gamma, 2\rho^*)$. This observation is particularly useful since we usually only know how to upper bound local complexity fixed points such as $r^*_{\mathrm{RMOM, ER}}(\gamma, 2\rho^*)$  and that we use it to define $\lambda$, the regularization parameter.

\paragraph{RMOM estimator with $L_2$ localization.} In this section, we look at the case where the $L_2$-norm to the square is used to describe the local curvature of the excess risk. As we mentioned above, it is the case when the margin assumption with margin parameter equals to $2$ holds. We define below the appropriate complexity fixed point parameter, the local curvature assumption and the associated sparsity equation. 

\begin{Definition}\label{def:complexity_RMOM}
Let $(\sigma_i)_{i\leq N}$ be independent Rademacher variables independent of the $X_i$'s. For $\rho>0$ and $\gamma > 0$, we define:
\begin{align*}
r^*_{\mathrm{RMOM}, L_2}(\gamma, \rho) := \inf \left(r>0 :\bE\left[\sup_{Z\in \cC: \norm{Z-Z^*}\leq \rho, \norm{Z-Z^*}_{L_2}\leq r} \left|\frac{1}{N}\sum_{i=1}^N\sigma_i \cL_{Z}(\tilde X_i)\right|\right] \leq \gamma r^2 \right).
\end{align*}
\end{Definition}


We turn now to the sparsity equation that is used to construct the radius $\rho^*$ which defines the model $\cC\cap (Z^*+\rho^* B)$ where both $Z^*$ and $\hat Z^{RMOM}_{K,\lambda}$ lie (with high probability).

\begin{Definition}\label{def:A-spar-MOM}
For $\gamma$, $\rho$ and $A>0$, let:
\begin{align*}
C_{K}(\gamma, \rho, A) := \max\left(320000A^2\frac{K}{N}, r^*_{\mathrm{RMOM}, L_2}(\gamma, \rho)^2\right), 
\end{align*}
\begin{align*}
\widetilde{H}_{\rho, A}:= \left\{Z\in \cC: \norm{Z-Z^*}=\rho \mbox{ and } \norm{Z-Z^*}_{L_2}\leq \sqrt{C_{K}(\gamma, \rho, A)}\right\}
\mbox{ and } 
\widetilde{\Delta}(\rho, A) := \inf_{Z\in \widetilde{H}_{\rho, A}} \sup_{\Phi \in \Gamma_{Z^*}(\rho)} \inr{\Phi, Z-Z^*}.
\end{align*}
A real number $\rho >0$ is said to satisfy the \textbf{$A$-sparsity equation} if $\widetilde{\Delta}(\rho, A) \geq 4\rho / 5$. 
\end{Definition}


The next definition is the formal way to say that the $L_2$-norm to the square can be used to describe the curvature of the excess risk closed to the oracle.

\begin{Assumption}\label{ass:curvature_RMOM_L2}
There exists $A$, $\gamma$ and $\rho^*>0$ such that $\rho^*$ satisfies the $A$-sparsity equation from Definition \ref{def:A-spar-MOM} and for both $b\in\{1, 2\}$ and all $Z\in \cC$, if $\norm{Z-Z^*}_{L_2}^2= C_{K}(\gamma, b\rho^*, A) $ and $\norm{Z-Z^*}\leq b\rho^*$, then $\norm{Z-Z^*}_{L_2}^2\leq AP\cL_Z$.
\end{Assumption}


After introducing the three key concepts in structural learning: local complexity fixed point, local curvature assumption and the sparsity equation, we can now state our excess risk and estimation (w.r.t. to both $L_2$ and the regularization norm) bounds.

\begin{Theorem}\label{theo:main-RMOM-unl}We consider the adversarial contamination setup of Assumption~\ref{ass:adv_contamination}. 
Let $K$ be a divisor of $N$ and assume that $K\geq 100|\cO|$.  Grant Assumption \ref{ass:curvature_RMOM_L2} for some $A\in(0, 1]$, $\gamma=1/32000$ and $\rho^*$ that satisfies the $A$-sparsity equation from Definition \ref{def:A-spar-MOM}. Define $\lambda = (11/(40\rho^*)) C_K(\gamma, 2\rho^*, A)$. Then it holds true that with probability at least $1-2\exp(-72K/625)$:
\begin{align*}
\norm{\hat{Z}^{\mathrm{RMOM}}_{K, \lambda}-Z^*}\leq 2\rho^* \quad ,\quad P\cL_{\hat{Z}^{\mathrm{RMOM}}_{K, \lambda}}\leq \frac{93}{100} r^*_{\mathrm{RMOM, L_2}}(\gamma, 2\rho^*)^2 \quad\mbox{ and }\quad \norm{\hat{Z}^{\mathrm{RMOM}}_{K, \lambda}-Z^*}_{L_2}^2\leq r^*_{\mathrm{RMOM, L_2}}(\gamma, 2\rho^*)^2.
\end{align*}
\end{Theorem}

Again the same result as the one of Theorem~\ref{theo:main-RMOM-unl} holds if one replaces $r^*_{\mathrm{RMOM, L_2}}$ by any upper bound on $r^*_{\mathrm{RMOM, L_2}}$. 

\paragraph{RMOM estimator with $G$ localization.} Finally, we consider  a function $G:H\rightarrow \bR$ that is expected to describe well the local curvature of the excess risk and that is used to define all the subsequent localization. An example of such a $G$ function is given in the sparse PCA case studied later. Indeed, in Lemma~\ref{coro:curvature_excess_risk} below, we will use $Z\in\bR^{d\times d}\to G(Z)=\norm{Z}_2^2$ as a localization function (we recall that $\norm{\cdot}_2$ is the canonical norm over $H$; it is in general different from the $L_2$ one that was used above for localization). We are now introducing a complexity fixed point that uses the $G$ function for localization. 

\begin{Definition}\label{def:G-loc-fixed-point-RMOM}
Let $\sigma_1, \ldots, \sigma_N$ be independent Rademacher variables independent of the $\tilde X_i$'s. For $G:H\rightarrow \bR$ and $A, \gamma$ and $\rho>0$, we define:
\begin{align*}
r^*_{\mathrm{RMOM, G}}(\gamma, \rho) := \mathrm{inf}\left\{r>0~:~\mathrm{max}\left(\frac{\mathrm{E}_G(r, \rho)}{\gamma},400\sqrt{2}V_{K,G}(r, \rho)
\right)\leq r^2\right\} 
\end{align*}
where, for all $r,\rho>0$,
\begin{align*}
\mathrm{E}_G(r, \rho) := \bE\left[\underset{Z\in\cC_{\rho, r}}{\mathrm{sup}}\left|\frac{1}{N}\sum_{i=1}^N\sigma_i\cL_Z(\tilde X_i)
\right|\right]
\mbox{ and }
V_{K,G}(r, \rho) := \sqrt{\frac{K}{N}} \underset{Z\in\cC_{\rho, r}}{\mathrm{sup}} \sqrt{\bV ar(\cL_Z(\tilde X))},
\end{align*}
with $\cC_{\rho, r} = \{Z\in\cC:\norm{Z-Z^*}\leq\rho,G(Z-Z^*)\leq r^2\}$
\end{Definition}
An example of computation of an upper bound of  the local complexity fixed point $r^*_{\mathrm{RMOM, G}}(\gamma, \rho)$ is provided in the sparse PCA example in Lemma~\ref{Lemma:r2-MOM} below. The final ingredient to derive the rate of convergence is the radius $\rho$ that needs to satisfies a sparsity equation.

\begin{Definition}\label{def:sparsity_RMOM_G-loc}
For all $\gamma$ and $\rho>0$, consider
$
\bar{H}_{\rho} := \left\{Z\in\cC:\norm{Z-Z^*}=\rho\mbox{ and }G(Z-Z^*)\leq r^*_{\mathrm{RMOM, G}}(\gamma, \rho)^2\right\}
$
and
$
\bar{\Delta}(\rho) := \inf_{Z\in \bar{H}_{\rho}} \sup_{\Phi \in \Gamma_{Z^*}(\rho)} \inr{\Phi, Z-Z^*}.
$
We say that $\rho$ satisfies the \textbf{sparsity equation} if $\bar{\Delta}(\rho)\geq4\rho / 5$. 
\end{Definition}

 Finally, we write the assumption saying that the $G$ function is indeed appropriate to describe the excess risk locally around $Z^*$. 
\begin{Assumption}\label{ass:curvature_RMOM_G-loc}
There exists $A>0$, $\gamma>0$ and $\rho^*>0$ such that $\rho^*$ satisfies the spartsity equation from Definition \ref{def:sparsity_RMOM_G-loc} and for both $b\in\{1, 2\}$ and  all $Z\in\cC$, if $G(Z-Z^*)=r^*_{\mathrm{RMOM, G}}(\gamma^*, b\rho^*)^2$ and $\norm{Z-Z^*}\leq b\rho^*$, then $AP\cL_Z\geq G(Z-Z^*)$. 
\end{Assumption}

We are now ready to state the following result on the statistical properties of the regularized minimax MOM in the context of robust structural learning with a linear loss function and for a general $G$ function describing the local curvature of the excess risk.

\begin{Theorem}\label{theo:main_RMOM_G-loc}
We consider the adversarial contamination setup of Assumption~\ref{ass:adv_contamination}. Let $G:H\rightarrow \bR$ be a continuous function such that $G(0)=0$ and for all $\alpha\geq1$ and $Z\in\cC, G(\alpha(Z-Z^*))\geq \alpha G(Z-Z^*)$. Let $K\in[N]$ be such that $K\geq 100|\cO|$. Grant Assumption \ref{ass:curvature_RMOM_G-loc} for some $A\in(0, 1]$, $\gamma=1/32000$ and $\rho^*$ that satisfies the sparsity equation from Definition~\ref{def:sparsity_RMOM_G-loc}. Define $\lambda=(11/(40\rho^*))r^*_{\mathrm{RMOM, G}}(\gamma, 2\rho^*) $. Then with probability at least $1-2\exp(-72K/625)$, it holds true that:
\begin{align*}
\norm{\hat{Z}^{\mathrm{RMOM}}_{K, \lambda}-Z^*}\leq 2\rho^* \quad , \quad P\cL_{\hat{Z}^{\mathrm{RMOM}}_{K, \lambda}} \leq \frac{93}{100}r^*_{\mathrm{RMOM, G}}(\gamma, 2\rho^*)^2 \quad\mbox{ and }\quad G(\hat{Z}^{\mathrm{RMOM}}_{K, \lambda}-Z^*) \leq r^*_{\mathrm{RMOM, G}}(\gamma, 2\rho^*)^2.
\end{align*}
\end{Theorem}

In the sparse PCA example, Theorem~\ref{theo:main_RMOM_G-loc} will be applied for the study of a $\ell_1$-regularized minmax MOM estimator. However, applying Theorem~\ref{theo:main_RMOM_G-loc}  requires several intermediate results such as proving that $Z\to G(Z)=\norm{Z}_2^2$ can be used as a local curvature of the excess risk, find a $\rho^*$ satisfying the sparsity equation of Definition~\ref{def:sparsity_RMOM_G-loc} and compute an upper bound for the local complexity fixed point $r^*_{\mathrm{RMOM, G}}(\gamma, \rho)$. For the last task, one  needs to handle the variance term $V_{K,G}$ as well as the complexity term $E_G(r,\rho)$. For the latter, we need to find an upper bound on the expected supremum of a Rademacher process over the interpolation body $\cC_{\rho, r} = \{Z\in\cC:\norm{Z-Z^*}\leq\rho,G(Z-Z^*)\leq r^2\}$. This step is usually the hardest one and requires some techniques from empirical process theory that we are now developing in the next section.


\section{Two examples of computation of local complexity fixed points}\label{sec:processus_sto}
In this section, we present concentration and in expectation results for two specific interpolation norms of the difference between the covariance matrix and its empirical version. These results are typical results that we use to compute local complexity fixed points like the ones used in the previous section. Indeed, in order to use any of the general statistical bounds presented in Section~\ref{sec:general_bounds}, we have to compute local complexity fixed points. We provide two such examples in this section that will be useful for the next  section on the sparse PCA problem. Note that the bounds presented here hold under weak moment assumptions (i.e. roughly speaking $\log(d)$ moments are enough) and may be of independent interest.

In this section, we use the following notations:  $X_1, \ldots, X_N$ are i.i.d. centered random vectors in $\bR^d$ and we denote by $\Sigma$ their covariance matrix, i.e. $\E X_1^\top X_1 = \Sigma$. The entries of $\Sigma$ are denoted by $\Sigma_{pq}$ i.e. $\bE X_{1p}X_{1q}=\Sigma_{pq}$ for all $p,q\in[d]$ where $X_1=(X_{1j})_{j=1}^d$. We denote the empirical covariance matrix by $\hat\Sigma_N = (1/N)\sum_{i=1}^N X_i^\top X_i$ and its entries by $\hat\Sigma_{pq}$, $p,q\in[d]$. The aim of this section is to provide large deviation and in expectation upper bounds for the norm of $\Sigma-\hat \Sigma_N$ for two norms defined by interpolation bodies. The proofs of the two Theorems \ref{theo:main_processus_sto} and \ref{theo:main_processus_sto_SLOPE} below are postponed to Section \ref{proofs:processus_sto}.

\subsection{Control of $\norm{\Sigma-\hat \Sigma_N}$ for a $B_2/B_1$ interpolation norm.}In order to upper bound the deviation of the empirical covariance matrix $\hat \Sigma_N$ around $\Sigma$ w.r.t. some norm we need to assume some concentration properties on the $X_i$'s. We therefore consider such an assumption now.
\begin{Assumption}\label{assum:weak_moment}
There exists $w\geq0$ and $t\geq1$ such that the following holds:
for all $p,q\in[d]$ and all $2\leq r\leq 2\log(ed/k)+t$ we have $\norm{X_{1p}X_{1q}-\bE(X_{1p}X_{1q})}_{L_r}\leq w^2 r$.
\end{Assumption}

In other words, Assumption~\ref{assum:weak_moment} is a growth condition on the first $2\log(ed/k)+t$ moments of the products $X_{1p}X_{1q}$ of the coordinates of $X_1$. This growth condition is the one exhibited by sub-exponential (i.e. $\psi_1$) variables. This is, for instance, the case of a product of two sub-gaussian (i.e. $\psi_2$) variables because $\norm{UV}_{\psi_1}\leq \norm{U}_{\psi_2}\norm{V}_{\psi_2}$ and the $r$-th moment of a $\psi_\alpha$ variable growths like $r^{1/\alpha}$ (see Chapter~1 in \cite{MR3113826} for more details). Assumption~\ref{assum:weak_moment} does not require the existence of any moment beyond  the $(2\log(ed/k)+t)$-th moment and is therefore called a weak moment assumption: Assumption~\ref{assum:weak_moment} essentially assumes the existence of $\log(ed/k)$ subgaussian moments on the coordinates of the data. We will see below that this assumption is enough to get estimation result for the first $k$-sparse principal component in deviation with an improved rate of convergence of order 
\begin{equation}\label{eq:improved_rate}
\sqrt{\frac{k^2 \log(ed/k)}{N}}.
\end{equation}

Let $k\in[d]$. We denote by $\norm{\cdot}$ the following interpolation pseudo-norm onto $\bR^{d\times d}$ defined by 
\begin{equation}\label{eq:norm_interpol_mat}
\norm{A}= \sup\left(\inr{A,Z}:Z\in k B_1\cap B_2\right).
\end{equation}

\begin{Theorem}\label{theo:main_processus_sto}There exists an absolute constant $c_0$ such that the following holds.
Grant Assumption~\ref{assum:weak_moment} for some $w$ and $t\geq1$ and assume that $N\geq 2\log(ed/k)+t$. With probability at least $1-\exp(-t)$,
\begin{equation*}
	\norm{\hat \Sigma_N -\Sigma}\leq c_0 w^2 \sqrt{\frac{k^2(\log(ed/k)+t)}{N}}.
	\end{equation*}	
Moreover, if  $N\geq 2\log(ed/k)+1$, it holds true that $\bE\left[\norm{\hat \Sigma_N -\Sigma}\right] \leq c_0 w^2 \sqrt{6k^2\log(ed/k)/N}$.
\end{Theorem}	

\begin{Remark}
Classical estimation result require the number of observations to be larger than $s\log(ed/s)$ where $s$ is the sparsity of signal to be reconstructed. Here, we observe in Theorem~\ref{theo:main_processus_sto} that $N$ is only asked to be larger than $\log(ed/k)$ so it is a much weaker assumption than in the classical high-dimensional setup. The rational behind this phenomenon is that we do not have to lower bound a quadratic process since our loss function is linear. It is usually isomorphic or just lower bounds results on a quadratic processes that require $N$ to be larger than the sparsity up to a log factor. We don't have such a quadratic process to lower bound in our 'linear loss function' framework.
\end{Remark}

\subsection{Control  of $\norm{\Sigma-\hat \Sigma_N}$  for a $B_2$/SLOPE interpolation norm.}\label{subsec:SLOPE_interpolation}
As in the last section, we need some assumption on the existence of moments on the coordinates of $X_1$. We consider  such an assumption now.
\begin{Assumption}\label{assum:weak_moment_slope}
There exists $w\geq0$ and $t\geq3$ such that the following holds.
For all $p,q\in[d]$ and all $2\leq r\leq \log(ed^2)+t$ we have $\norm{X_{1p}X_{1q}-\bE(X_{1p}X_{1q})}_{L_r}\leq w^2 r$.
\end{Assumption}

Our aim is to analyze the statistical properties of a SLOPE regularization for the sparse PCA problem and to show that the optimal rate \eqref{eq:improved_rate} can be achieved by a unique regularization method which does not require the a priori knowledge of the sparsity parameter $k$. To that end we introduce the SLOPE regularization norm of a $d\times d$ matrix $A$
\begin{equation*}
\norm{A}_{SLOPE} = \sum_{p,q=1}^db_{pq}A_{(p,q)}^*
\end{equation*}where $\textbf{b}:=(b_{pq}:p,q\in[d])$ are decreasing  weights for some lexicographical order over $[d]^2$ starting at $(1,1)$ such that for all $k\in[d]$, $b_{kk}=\sqrt{\log(ed^2/k^2)+t}$. For instance, one may assume that $b$ is a symetric matrice and set $b_{pq}=\sqrt{\log(ed^2/(pq))+t}$ when $q\geq p$. We also denote by  $(A_{(p,q)}^*:p,q\in[d])$ the non-increasing sequence (for the same lexicographical order over $[d]^2$ used before) of the rearrangement of the absolute values of the entries of $A$, for instance $A_{(d,d)}^* = \min(|A_{pq}|:p,q\in[d])$ and $A_{(1,1)}^* = \max(|A_{pq}|:p,q\in[d])$. We denote by $B_{SLOPE}$ the unit ball of the SLOPE norm.

Let $\rho>0$. We denote by $\norm{\cdot}_\rho$ the following interpolation pseudo-norm onto $\bR^{d\times d}$ defined by 
\begin{equation}\label{eq:norm_interpol_mat_slope}
\norm{A}_\rho= \sup\left(\inr{A,Z}:Z\in \rho B_{SLOPE}\cap B_2\right).
\end{equation}

\begin{Theorem}\label{theo:main_processus_sto_SLOPE}There exists an absolute constant $c_0$ such that the following holds. Let $k\in[d]$ and $\gamma\geq1$.
Grant Assumption~\ref{assum:weak_moment_slope} for some $w$ and $t\geq \max\big(2 \log(\lceil \log(k^2)\rceil), \gamma \log(ed^2/k^2)\big)$ and assume that $N\geq \log(ed^2)+t$. With probability at least $1-2\exp(-t/2)$,
\begin{equation*}
   \norm{\hat \Sigma_N -\Sigma}_\rho\leq \frac{c_0 w^2}{\sqrt{N}}\min(\rho, d).
   \end{equation*}   
\end{Theorem}


\section{Sparse PCA}\label{sec:sparsePCA}

Principal Components analysis (PCA) is one of the most fundamental dimension reduction algorithm as well as one of the most used data visualization tool. It can be efficiently performed via some truncated SVD algorithms on the $N\times d$ data matrix ($N$ being the number of data and $d$ the dimension of the data, that is the number of features) which requires only $\cO(k^2\min(d,N))$ operations to get the first $k$ top eigenvectors \cite{MR2806637,MR3024913}. 

However, principal components are linear mixture of features that may be of very different nature  and as so are for most of the time meaningless. This problem becomes more salient for high-dimensional data (i.e. when $d>N$) where the diversity of features (text, socio-professional categories, geographic location, familiar situation, consumption habits, etc.)  may be very large. Moreover, in the high-dimensional setting, PCA no longer provides meaningful estimates of the principal components of the actual covariance matrix $\Sigma$ as exhibited by the phase transition from \cite{MR2165575}. 

One way to alleviate both interpretation and inconsistency in the high-dimensional setting is to look for principal components which are linear mixture of a small number of features -- that is ''sparse'' principal component. This problem is known as sparse PCA and was introduced in \cite{johnstone2009sparse,MR2751448}. It can be stated as the following optimization problem:
\begin{equation}\label{eq:sPCA}
\hat v_1\in\argmax_{\norm{v}_2 = 1, \norm{v}_0\leq k}\norm{\hat \Sigma_N v}_2
\end{equation}where the $X_i$'s are $i.i.d$ centered vectors in $\bR^d$ with covariance $\bE[X_iX_i^\top]=\Sigma$, $\hat \Sigma_N=(1/N)\sum_{i=1}^N (X_i-\bar X_N)(X_i-\bar X_N)^\top$ is the empirical covariance matrix, $\norm{v}_0$ is the size of the support of $v$ and $k$ is some fixed sparsity level. 

From an algorithmic point of view there are two major issues in the optimization problem \eqref{eq:sPCA}: 1) the objective function that we want to maximize is convex; and it is notoriously difficult to maximize a convex function even on a convex set 2) because of the sparsity constraint '$\norm{v}_0\leq k$', the constraint set is not convex. If the sparsity constraint was not there, then \eqref{eq:sPCA} would be the classical PCA problem for finding a first principal component, that is a top eigenvector of $\hat\Sigma_N$. In that case, even though it is a maximization problem of a convex function on a convex set, this problem can be solved efficiently for instance via the power method and is in fact one of the few situation where maximizing a convex function can be performed efficiently.

The extra sparsity constraint in \eqref{eq:sPCA} somehow emphasis this original issue that the objective function to maximize is convex. One way to overcome this issue is to adapt the power method to this extra constraint, see \cite{MR2600619}. Another way is via SDP relaxation \cite{MR2353806}. We will use this latter approach so we present it in the next subsection in more details.

\subsection{SDP relaxation in sparse PCA}
Let $X\in\bR^d$ be a centered random vector with distribution $P$. Let $X_1, \ldots, X_N\in \bR^d$ be N independant copies of $X$. Define $A:= (1/N)\sum_{i=1}^NX_iX_i^\top$, the empirical covariance matrix of the $X_i$'s. Let $\Sigma := \bE[A] = \bE_{X\sim P}[XX^T]$ be their covariance matrix. We are looking for a first principal component with a support of small cardinality, that is for a vector $v^*\in\bR^d$ with unit-length and cardinality less than a certain integer $k\leq d$, and such that the variance of the $X_i$'s when projected onto $v^*$ is maximal. This can be written as follows:
\begin{equation}\label{pb:card_constraint}
v^* \in \argmax_{v\in \cE}\E [\inr{X, v}^2] \mbox{ where } \cE := \left\{v\in\bR^d:\norm{v}_2 = 1, \norm{v}_0\leq k\right\}.
\end{equation}
This problem is known to be NP-hard in general \cite{magdonismail2015nphardness}, so we are looking to relax it. One way to do this is to replace the cardinality function by the $\ell_1^d$-norm. Another way is via the lifting procedure, which is described for example in \cite{lemarechal1999DSRelaxation} and is based on the principle that quadratic  objective functions and constraint sets of a vector $v$ can be written as linear objective functions and constraint sets of the symmetric rank one matrix $vv^\top$ .

In our case, we first note that $\E [\inr{X, v}^2] = \inr{\bE[A], vv^\top} = \inr{\Sigma, vv^\top}$. Then, if $Z=vv^T$ with $v \in S_2^{d-1}$ and $\norm{v}_0 \leq k$, we have $\mathrm{Tr}(Z) = \norm{v}_2^2= 1$ and $\norm{Z}_0 \leq k^2$. Finding a solution of \eqref{pb:card_constraint} is then equivalent \cite{MR2353806,lemarechal1999DSRelaxation} to finding a top singular vector of $Z^\star$, where $Z^\star$ is solution of the optimization problem
\begin{equation*}
Z^\star \in  \argmax_{Z \in \cC} \inr{\bE[A], Z}  \mbox{ where }  \cC_0 := \left\{Z \in \bR^{d\times d} : Z = vv^T, v \in \bR^{d},  \mathrm{Tr}(Z)=1, \norm{Z}_0 \leq k^2\right\}.
\end{equation*}
In the latter problem, the objective function has now become a linear one thanks to the lifting approach, however the constraint set is not convex. We are now working on that issue to get a full SDP relaxation of (\ref{pb:card_constraint}). First, we may replace the condition \doublequote{$Z = v v^T$} by the equivalent condition \doublequote{$Z\succeq 0$ and $\mathrm{rank}(Z) = 1$} in $\cC_0$. However, $\cC_0 := \left\{Z \in \bR^{d\times d} : Z \succeq 0, \mathrm{Tr}(Z)=1, \norm{Z}_0 \leq k^2, \mathrm{rank}(Z)=1\right\}$ is not convex, because of two non-convex constraints: \textit{the cardinality constraint \doublequote{ $\norm{Z}_0 \leq k^2$}} and \textit{the rank constraint \doublequote{ $\mathrm{rank}(Z)=1$}} that we are just dropping out of $\cC_0$.  By doing so, we end up with the following convex optimization problem:
\begin{equation}\label{pb:card_constrained_convex}
Z^* \in \argmax_{Z\in \cC} \inr{\bE[A], Z} \mbox{ where } \cC := \{Z \in \bR^{d\times d}: Z\succeq 0, \mathrm{Tr}(Z)=1\}.
\end{equation}
We then see $Z^*$ as an oracle for the linear loss function $Z \to \ell_Z(X) = - \inr{XX^\top, Z}$ and its associated risk function $Z\to \bE \ell_Z(X)$ over the model $\cC$, that is $Z^* \in \argmin_{Z\in \cC} P\ell_Z$. This enables us to leverage the methodological tools introduced in Section \ref{sec:general_bounds} to derive estimators for $Z^*$ and provide statistical guarantees onto them.

This configuration allows us to refer to the work of Samworth et al. \cite{MR3546438}. The authors study the sparse PCA problem where the distribution of the data $X_1, \ldots, X_N$ belongs to a class $\cP$ of distributions that all have a sub-exponential tail; it includes, among others, sub-Gaussian distributions (see equation~(4) in \cite{MR3546438} for a definition).  In particular, they propose the following $\ell_1$-regularized ERM estimator
\begin{align}\label{def:estimator_sparsePCA_Berthet}
\hat Z \in \argmin_{Z\in\cC}\left( \inr{\frac{-1}{N}\sum_{i=1}^N X_iX_i^\top, Z} + \lambda\norm{Z}_1\right) \mbox{ where } \cC:=\left\{Z:Z\succeq 0, \mathrm{Tr}(Z) = 1\right\} 
\end{align}
and provide an algorithm for solving it in polynomial time.  We report below their main results for this estimator.  
\begin{Theorem}\label{Theo:Samworth_sparsePCA}[Theorem 5 in \cite{MR3546438}]
Let  $X_1, \ldots, X_N\in \bR^d$ be $i.i.d$ random vectors with distribution in $\cP$ and a covariance matrix satisfying the spiked covariance model: $\bE[X_iX_i^\top] = I_d+\theta\beta^*(\beta^*)^\top$, where $\beta^*$ is a $k$-sparse vector with unit euclidean norm. Let $\lambda = 4\sqrt{\log(d)/N}$, $\epsilon =\log(d)/(4N)$ and consider $\hat{v}_{\lambda, \epsilon} \in \argmax_{\norm{v}_2=1}v^\top\hat{Z}^{\epsilon}v$, where $\hat{Z}^{\epsilon}$ is an $\epsilon$-maximizer of $Z\rightarrow  \inr{\frac{1}{N}\sum_{i=1}^NX_iX_i^\top, Z} - \lambda\norm{Z}_1$ over the model $\cC$ defined in (\ref{def:estimator_sparsePCA_Berthet}). Finally, let $\hat{v}^0_{\lambda, \epsilon}$ be the $k$-sparse vector derived from $\hat{v}_{\lambda, \epsilon}$ by setting all but its largest $k$ coordinates in absolute value to $0$. If $4\log(d)\leq N\leq k^2d^2\theta^{-2}$ and $0<\theta\leq k$, then it holds true that:
\begin{align*}
\bE\left[\sqrt{2}\norm{\hat{v}^0_{\lambda, \epsilon}(\hat{v}^0)_{\lambda, \epsilon}^\top - \beta^*(\beta^*)^\top}_2\right] \leq (32\sqrt{2}+3)\sqrt{\frac{k^2\log(d)}{N\theta^2}}.
\end{align*}
\end{Theorem}

We are now using our methodology to propose several estimators and provide our insights on the sparse PCA problem. In particular, we will extend Theorem~\ref{Theo:Samworth_sparsePCA} to the heavy-tailed framework, provide in-deviation results and improve the rate to the optimal one $k^2\log(ed/k)/N$ (thanks to localization). On top of that, we will construct new estimators based on the MOM principle to handle  robustness issues in sparse PCA.

\subsection{Exactness and curvature in the spiked covariance model.}
We present here two results that will be of crucial importance in the analysis of our estimators (the proofs are postponed to Section \ref{sec:sparsePCA}). The first one concerns the exactness in the spiked covariance model. That is, the oracle $Z^*$ as defined by equation (\ref{pb:card_constrained_convex}), obtained after a lifting and a convex relaxation of the initial problem, turns out to be a matrix of rank one whose unit-norm leading eigenvector is $\pm\beta^*$. 

\begin{Lemma}\label{coro:exactness1_spiked} 
In the spiked covariance model $\Sigma=\theta (\beta^*)(\beta^*)^\top + I_d$ with $\beta^*\in S_2^{d-1}$ and $\beta^*$ is $k$-sparse, we have $Z^*=(\beta^*)(\beta^*)^\top$, for $Z^*$ defined in (\ref{pb:card_constrained_convex}).
\end{Lemma}

The second one concerns the curvature of the excess risk function around the oracle $Z^*$. Following our methodology, we need to understand the behavior of the excess risk around $Z^*$ in order to find a good $G$ function that will be used to be define localized subsets of our model. Then, later, based on the results from Section~\ref{sec:processus_sto} we will compute the  Rademacher complexities of these localized subsets  and then the local complexity fixed points as introduced in Section \ref{sec:general_bounds}. The fixed point is then used to establish statistical bounds on our estimators.  Finding the 'right' curvature function of the excess risk is therefore important in our approach. The following result provides a curvature of the excess risk \simplequote{globally}, that is on the entire set $\cC$ and not just around $Z^*$ (see the proof in Section \ref{coro:curvature_excess_risk}).

\begin{Lemma}\label{coro:curvature_excess_risk}
In the spiked covariance model $\Sigma=\theta (\beta^*)(\beta^*)^\top + I_d$ with $\beta^*\in S_2^{d-1}$ and $\beta^*$ is $k$-sparse, the following holds. For all $Z\in\cC$, we have $P\cL_Z = \inr{\Sigma, Z^*-Z}\geq(\theta/2)\norm{Z^*-Z}_2^2$.
\end{Lemma}

As a consequence, using our terminology, the problem has an excess risk curvature function given by $G:Z\to \norm{Z}_2^2$ - where $\norm{\cdot}_2$ is the canonical Hilbertian norm in $\bR^{d\times d}$. We will therefore use the $\ell_2$-norm to the square to define our localized models for the study of all estimators introduced below.

\subsection{$\ell_1$-Regularized ERM estimator}

Since the parameter we want to estimate has a sparse structure, the choice of estimators regularized by an appropriate norm will enable us to take advantage of this structural property. We start with a regularized ERM estimator, as presented in Section \ref{sec:RERM}, where the $\ell_1$-norm is used as regularization norm:
\begin{equation}\label{pb:RERM_card_penalized_convex}
\hat{Z}^{\mathrm{RERM}}_\lambda \in \argmin_{Z\in \cC} \left( P_N\ell_Z + \lambda\norm{Z}_1\right),\quad \mbox{where}\quad  \cC:=\left\{Z\in \bR^{d\times d}: Z \succeq 0, \mathrm{Tr}(Z)=1\right\}
\end{equation}and $\ell_Z(X) = -\inr{XX^\top, Z}$ and $P_N \ell_Z = (1/N)\sum_{i=1}^N \ell_Z(X_i)$. This puts us in condition to use the results of Section~\ref{sec:RERM} to provide statistical guarantees on $\hat{Z}^{\mathrm{RERM}}_\lambda$.

Lemma \ref{coro:curvature_excess_risk} shows that, for any value of $\rho >0$ and $\delta\in(0, 1)$, Assumption \ref{ass:A-spars} is satisfied with $A=2/\theta$ and $G:Z\in\bR^{d\times d}\to\norm{Z}_2^2$. In order to proceed with our methodology, the next step is then to identify a value of $\rho^*$ which satisfies the $2/\theta$-sparsity equation from Definition \ref{def:RERM_A-sparsity}. This is the purpose of the following Lemma (the proof is given in section \ref{proof:lemma_sparsity_sparsePCA_RMOM_l1}).

\begin{Lemma}\label{Lemma:rho-spars}
Let $A>0$, $\delta \in (0, 1)$, and define $r^*(.):= r^*_{\mathrm{RERM, G}}(A, ., \delta)$. If $\rho \geq 10k\sqrt{r^*(\rho)}$, then $\rho$ satisfies the $A$-sparsity equation from Definition \ref{def:RERM_A-sparsity}. 
\end{Lemma}

The last step is to compute the local complexity fixed point of Definition \ref{def:RERM_Complexity_G_function}, which is what we are working on below. 


\begin{Lemma}\label{Lemma:r2_bound}
Grant Assumption \ref{assum:weak_moment} with $t= \log(ed/10k)$. Suppose that $\beta^*$ is $k$-sparse, with $k\leq ed/200$. Let $A=2/\theta$ and assume that $N\geq 3\log\left(ed/10k\right)$. Then there exists an absolute constant $b>0$ such that, defining:
\begin{align}\label{eq:rstar_RERM_L1}
 \rho^* := 200 bAk^2\sqrt{\frac{1}{N} \log\left(\frac{ed}{k}\right)}  \mbox{ and }   r^*(\rho) := bA\sqrt{\frac{\rho^2}{N}\log\left(\frac{b^2A^2(ed)^4}{N\rho^2}\right)},
\end{align} 
one has $r^*_{\mathrm{RERM, G}}\left(A, \rho^*, 10k/ed\right) \leq r^*(\rho^*)$ and $\rho^*$ satisfies the $A$-sparsity equation from Definition \ref{def:RERM_A-sparsity}.
\end{Lemma}

We are now ready to state our main result concerning the $\ell_1$-regularized ERM estimator for the sparse PCA problem.


\begin{Theorem}\label{theo:main_RERM_SPCA}
Grant Assumption \ref{assum:weak_moment} with $t = \log(ed/10k)$. Suppose that $\beta^*$ is $k$-sparse, with $k\leq ed/200$. Assume that $N\geq 3\log\left(ed/10k\right)$ and that $\lambda$ satisfies the following inequalities:
\begin{align}\label{eq:main_RERM_L1_cond_lambda}
\frac{20}{21}b\sqrt{\frac{1}{N}\log\left(\frac{ed}{200^{1/2}\log(200)^{1/4}k}\right)}  \leq \lambda \leq \frac{2}{\sqrt{3}}b \sqrt{\frac{1}{N}\log\left(\frac{ed}{200^{2/3}k}\right)}
\end{align}
where $b$ is the absolute constant introduced in Lemma \ref{Lemma:r2_bound} above. Let $C = 40b$. Then, with probability at least $1-10k/ed$, it holds true that:
\begin{equation*}
 \norm{\hat{Z}^{\mathrm{RERM}}_{\lambda}-Z^*}_1  \leq  10C k^2\sqrt{\frac{1}{N\theta^2}\log\left(\frac{ed}{k}\right)}, \quad  
 \norm{\hat{Z}^{\mathrm{RERM}}_{\lambda}-Z^*}_2  \leq  C \sqrt{\frac{k^2}{N\theta^2}\log\left(\frac{ed}{k}\right)}
 \end{equation*}
 and
 \begin{equation*}
P\cL_{\hat{Z}^{\mathrm{RERM}}_{\lambda}}  \leq \frac{C^2}{2}\frac{k^2}{N\theta}\log\left(\frac{ed}{k}\right).
   \end{equation*}  
\end{Theorem}
Note that if one is willing to get a better deviation parameter, one can assume $N$ larger than $\Upsilon \log(ed/10k)$, for $\Upsilon$ large enough. \\ 

Up to this point, we have introduced an estimator for $Z^*$ and provided a convergence rate with high probability. However, our primary focus is not on $Z^*$ itself, but rather on its unit-norm leading eigenvectors $\pm\beta^*$. The purpose of the upcoming result is to leverage the preceding one in order to establish properties related to $\beta^*$.

\begin{Corollary}\label{coro:main_RERM_L1}
Let $\hat{\beta}$ $\in\bR^d$ be a leading unit length eigenvector of $\hat{Z}^{\mathrm{RERM}}_{\lambda}$. Under the conditions of Theorem \ref{theo:main_RERM_SPCA}, there exists an absolute constant $D>0$ such that with probability at least  $1-10k/ed$:
\begin{align*}
\norm{\hat{\beta}\hat{\beta}^\top - \beta^*(\beta^*)^\top}_2\leq D\sqrt{\frac{k^2}{N\theta^2}\log\left(\frac{ed}{k}\right)}.
\end{align*}
\end{Corollary}

We therefore obtain a convergence rate of magnitude $\left(k^2\log\left(ed/k\right)/(N\theta^2)\right)^{1/2}$, when our dataset is made up of $i.i.d$ random variables whose distribution satisfies Assumption \ref{assum:weak_moment}, which includes the case of $i.i.d$ sub-Gaussian variables but it goes much beyond up to variables with only $\log d$ moments. The result of \cite{MR3546438} is available for a class of distributions, including sub-Gaussian distributions, whose covariance matrix fits within the spiked covariance model. They obtain a convergence rate of magnitude $\left(k^2\log(d)/(N\theta^2)\right)^{1/2}$, although our result holds with polynomial deviation while theirs is in expectation. We also note that our result does not suffer from any restrictive condition concerning $\theta$. 
We therefore slightly improve the results from \cite{MR3546438}; this improvement is of the same order as the one obtained for the LASSO in \cite{MR3852663} and is due to a careful localization argument. This shows that our analysis is precise enough to catch the subtle difference between the $\log d$ rate from \cite{MR3546438} and the $\log(ed/k)$ obtained in Theorem~\ref{theo:main_RERM_SPCA}.  Our result also extend the scope of Theorem~\ref{Theo:Samworth_sparsePCA} to heavy-tailed data since we only require the existence of $\log d$ moments. However, to get this improvement for the Lasso type estimator (\ref{pb:RERM_card_penalized_convex}), one needs to choose $\lambda$ depending on $k$ in \eqref{eq:main_RERM_L1_cond_lambda}, which is unknown in practice. To solve this issue, we could use a Lepskii's adaptation scheme as in \cite{MR3852663}. However, we will not follow this path but rather consider another regularization norm: the $SLOPE$ norm, that allows to get the same results as in  Theorem~\ref{theo:main_RERM_SPCA} but a choice of $\lambda$ independent of $k$. This will also give us the opportunity to run our methodology one more time for a different regularization  norm.

\subsection{$SLOPE$ regularized ERM estimator}\label{sec:RERM_SLOPE_sparsePCA}

In this section, we study a regularized ERM estimator of $Z^*$ with the $SLOPE$ norm (introduced in Section \ref{subsec:SLOPE_interpolation}, and whose definition is restated below) as the regularization norm. We consider a lexicographical order over $[d]^2$ such that for any $k\in[d]$, the $k^2$ largest elements in $[d]^2$ belong to $[k]^2$. 
We fix $t>0$ (which will be choosen appropriately later) and we define, for $p\leq q$, $b_{pq}(t)=:\sqrt{\log(ed^2/pq)+t}$, and $b_{pq}(t) = b_{qp}(t)$ for $p>q$. For $Z \in \bR^{d\times d}$, we define $Z^\sharp$ the matrix obtained from $Z$ by reordering its element in absolute value in non-increasing order, and we finally define its $SLOPE$ norm by:
\begin{align*}
 \norm{Z}_{SLOPE} := \sum_{p,q=1}^db_{pq}Z^\sharp_{pq}.
\end{align*} 
 Our estimator is then:
\begin{equation}\label{pb:RERM_card_penalized_convex}
\hat{Z}^{\mathrm{RERM}}_{SLOPE} \in \argmin_{Z\in \cC}  \,\left(P_N\ell_Z + \lambda\norm{Z}_{SLOPE}\right) \mbox{ for } \cC := \{Z \in \bR^{d\times d}: Z\succeq 0, \mathrm{Tr}(Z)=1\}
\end{equation}
and a regularization parameter $\lambda > 0$ to be chosen later. This puts us in condition to use the results of Section \ref{sec:RERM} to provide statistical guarantees on $\hat{Z}^{\mathrm{RERM}}_{SLOPE}$.  

As before, the essence of Lemma \ref{coro:curvature_excess_risk} in this context is that, for any value of $\rho >0$ and $\delta\in]0, 1[$, Assumption \ref{ass:A-spars} is satisfied with $A=2/\theta$ and $G:Z\in\bR^{d\times d}\to\norm{Z}_2^2$. In order to proceed with our methodology, our next step is then to identify a value of $\rho^*$ which satisfies the $2/\theta$-sparsity equation. This is the purpose of the following Lemma. 

\begin{Lemma}\label{Lemma:SPCA_SLOPE_A-sparsity}
Assume that $\beta^*$ is $k$-sparse, for some $k\in [d]$. Let $A>0$, $\delta \in (0,1)$ and $t>0$. Define $\Gamma_k(t):=3\sum_{\ell = 1}^k b_{\ell\ell}(t)$. If $\rho \geq 10\Gamma_k(t)\sqrt{r^*_{\mathrm{RERM, G}}(A, \rho, \delta)}$, then $\rho$ satisfies the $A$-sparsity equation from Definition \ref{def:RERM_A-sparsity}. 
\end{Lemma}

Following the path traced by our methodology, all that remains is to calculate the complexity fixed-point parameter

\begin{align*}
r^*_{\mathrm{RERM, G}}(A, \rho, \delta) = \inf\left(r > 0: \bP\left(\sup_{Z\in \cC : \norm{Z-Z^*}_{SLOPE}\leq \rho, \norm{Z-Z^*}_2\leq \sqrt{r}}|(P-P_N)\cL_Z| \leq \frac{r}{3A}\right)\geq 1-\delta\right).
\end{align*}

The next Lemma gives us an upper bound for $r^*_{\mathrm{RERM, G}}(A, \rho, \delta) $, when $\rho$ satisfies the sparsity equation of Definition \ref{def:RERM_A-sparsity}.

\begin{Lemma}\label{Lemma:SPCA_SLOPE_fixed-point}
Grant Assumption \ref{assum:weak_moment_slope} for $t= 2\log(ed^2/k^2)$. Suppose that $\beta^*$ is $k$-sparse, with $k\leq d/(e^2\log(d))$. Let $A>0$, and assume that $N\geq 3\log(ed^2)$. Then, there exists an absolute constant $b>0$ such that, defining:
\begin{align*}
\rho^*:=10\Gamma_k^*\frac{bA}{\sqrt{N}}\min\left(10\Gamma_k^*;d\right) ~~ \mbox{ and } ~~  r^* := \frac{b^2A^2}{N}\min\left(10\Gamma_k^*; d\right)^2
\end{align*}
one has $r_{\mathrm{RERM, G}}^*(A, \rho^*, 2k^2/(ed^2))\leq r^*$ and $\rho^*$ satisfies the A-sparsity equation rom Definition \ref{def:RERM_A-sparsity}, where $\Gamma_k^* = \Gamma_k(2\log(ed^2/k^2))$ is the quantity introduced in Lemma \ref{Lemma:SPCA_SLOPE_A-sparsity}.
\end{Lemma}

We are now ready to state our main result concerning the $SLOPE$ regularized ERM estimator for the sparse PCA problem.

\begin{Theorem}\label{theo:main_SPCA_RSLOPE}
Grant Assumption \ref{assum:weak_moment_slope} for $t = 2\log(ed^2/k^2)$. Suppose that $\beta^*$ is $k$-sparse, with \\ $k\leq \min\left(d/(e^2\log(d)), (e/140\sqrt{2})^2d\right)$. Assume that $N\geq 3\log(ed^2)$ and that $\lambda$ satisfies the following inequalities:
\begin{align}\label{eq:lambda_RMOM_SLOPE}
\frac{10b}{21\sqrt{N}} < \lambda < \frac{2b}{3\sqrt{N}},
\end{align}
where $b$ is the constant previously defined in Lemma \ref{Lemma:SPCA_SLOPE_fixed-point}. Then there exist an absolute constants $C_1>0$ such that one has with probability at least $1 - 2k^2/(ed^2)$:
\begin{align*}
\norm{\hat{Z}^{\mathrm{RERM}}_{SLOPE}-Z^*}_{SLOPE} & \leq C_1\frac{k^2}{\sqrt{N\theta^2}}\log\left(\frac{ed^2}{k^2}\right),\\ 
\norm{\hat{Z}^{\mathrm{RERM}}_{SLOPE}-Z^*}_{2} & \leq C_1\sqrt{\frac{k^2}{N\theta^2}\log\left(\frac{ed^2}{k^2}\right)}
\end{align*} and  
\begin{equation*}
\inr{\Sigma, Z^*-\hat{Z}^{\mathrm{RERM}}_{SLOPE}}   \leq C_1\frac{k^2}{N\theta}\log\left(\frac{ed^2}{k^2}\right).
\end{equation*}
\end{Theorem}

We can now use this result to obtain properties about our object of interest, which is not directly $Z^*$, but its unit-length leading eigenvectors $\pm\beta^*$. 

\begin{Corollary}\label{coro:main_SPCA_RSLOPE}
Let $\hat{\beta}$ $\in\bR^d$ be a leading unit-eigen vector of $\hat{Z}^{\mathrm{RSLOPE}}_{\lambda}$. Under the conditions of Theorem \ref{theo:main_SPCA_RSLOPE}, there exists an absolute constant $C>0$ such that with probability at least  $1-2k^2/ed^2$:
\begin{align*}
\norm{\hat{\beta}\hat{\beta}^\top - \beta^*(\beta^*)^\top}_2\leq C\sqrt{\frac{k^2}{N\theta^2}\log\left(\frac{ed^2}{k^2}\right)}.
\end{align*}
\end{Corollary}

Here again, we obtain a rate of convergence of magnitude $\sqrt{(1/N\theta^2)\log\left(ed^2/k^2\right)}$, holding with polynomial deviation, with no restriction on the value of $\theta$. We note that this result holds with a value of the regularization parameter $\lambda$ that does not depend on the sparsity level $k$ of $\beta^*$.

\subsection{$\ell_1$ regularized minmax MOM estimator.}
Here, we consider the case where data may be corrupted with outliers. We place ourselves in the framework of the adversarial contamination, which is described in Assumption \ref{ass:adv_contamination}: the dataset $\{X_1, \ldots, X_N\}$ used by the statistician may have been  corrupted by an adversary. As a consequence,  on top of the structural learning problem, we now have to face a robustness to data contamination problem. To deal with these issues all together, we use a regularized minmax MOM estimator. 

We therefore consider an equi-partition of $\left\{1, \ldots, N\right\}$ into $B_1\sqcup\cdots \sqcup B_K=[N]$, where $|B_k|= N/K$ for all $k\in [K]$. We consider a $\ell_1$-regularized minmax MOM estimator
\begin{align*}
\hat{Z}^{RMOM}_{K,\lambda} \in \argmin_{Z\in \cC} \sup_{Z^\prime \in \cC} \left(\mathrm{MOM}_K (\ell_Z - \ell_{Z^\prime}) + \lambda(\norm{Z}_1 - \norm{Z^\prime}_1)\right) 
\end{align*}for  $\cC:=\left\{Z\in \bR^{d\times d}:0\preceq Z \preceq I_d, \mathrm{Tr}(Z)=1\right\}$ and a regularization parameter $\lambda$ to be chosen later. 

In what follows, we provide some statistical guarantees on $\hat{Z}^{RMOM}_{K,\lambda}$ based on Theorem~\ref{theo:main_RMOM_G-loc} which is our general result for regularized minmax MOM estimators for a general $G$ function used for localization. Here, following Lemma~\ref{coro:curvature_excess_risk}, we will use $G:Z\to \norm{Z}_2^2$ (and $A=2/\theta$) for such a localization function. Following our methodology, once the curvature of the excess risk is chosen, we have to find an upper bound on the local complexity fixed point $r^*_{\mathrm{RMOM}, G}(\gamma, \rho)$ from Definition~\ref{def:G-loc-fixed-point-RMOM}. But before that we find a sufficient condition on a radius $\rho$ so that it satisfies  the sparsity equation from Definition~\ref{def:sparsity_RMOM_er}.

\begin{Lemma}\label{Lemma:rho-spars-RMOM-SPCA}
Consider $\gamma >0$. If $\rho >0$ is such that $\rho\geq 10k\sqrt{2/\theta}r^*_{\mathrm{RMOM}, G}(\gamma, \rho)$, then $\rho$ satisfies the sparsity equation from Definition~\ref{def:sparsity_RMOM_er}.
\end{Lemma}

Now that we know how to grasp a value of $\rho$ that satisfies the sparsity equation, the subsequent task is to compute the fixed-point parameter $r^*_{\mathrm{RMOM}, G}(\gamma, \rho)$ as introduced in Definition~\ref{def:G-loc-fixed-point-RMOM}, after which, thanks to Theorem \ref{theo:main_RERM_er}, we will be able to provide some statistical bounds on $\hat{Z}^{RMOM}_{K, \lambda}$.


\begin{Lemma}\label{Lemma:r2-MOM}
Grant assumption \ref{assum:weak_moment} for $t=1$. Suppose that $\beta^*$ is $k$-sparse, for some $k\in[d]$. Assume that $N\geq 2\log(ed/k)+1$ and that $\theta \leq k$. Define $G:Z\in\bR^{d\times d}\to (\theta/2)\norm{Z}_2^2$. Consider $\gamma>0$. There exist absolute constants $B$ and $D>0$ such that, defining:
\begin{equation*}
\rho^*(\gamma) := \max\left(\sqrt{480}B\frac{k^2}{\gamma}\sqrt{\frac{1}{N\theta^2}\log\left(\frac{ed}{k}\right)}; 10Dk \sqrt{\frac{2K}{N\theta^2}}\right)
\end{equation*}and  
\begin{equation*}
r^*(\gamma, \rho) := \max\left(\sqrt{\frac{B\rho}{\gamma}}\left(\frac{6}{N}\log\left(\frac{2B(ed)^2}{\gamma\theta\rho}\sqrt{\frac{6}{N}}\right)\right)^{1/4} ; D\sqrt{\frac{K}{N\theta}}\right)
\end{equation*}
one has  $r^*_{\mathrm{RMOM, G}}(\gamma, \rho^*(\gamma)) \leq r^*(\gamma, \rho^*(\gamma))$ and $\rho^*(\gamma)$ satisfies the sparsity equation from Definition \ref{def:sparsity_RMOM_er}. The values of $B$ and $D$ are explicited in Section \ref{proof:Lemma_r2-MOM}. 
\end{Lemma}

We are now ready to state our main result about the $\ell_1$-Regularized MOM estimator for the sparse PCA problem. 

\begin{Theorem}\label{theo:main-RMOM-SPCA}
Grant assumption \ref{assum:weak_moment} for $t=1$. Suppose that $\beta^*$ is $k$-sparse, for some $k\in[d]$. Assume that $N\geq 2\log(ed/k)+1$ and let $K$ be a divisor of $N$ such that $K\geq 100|\cO|$. Let $\gamma = 1/32000$ and $\lambda := 11r^*(\gamma, 2\rho^*(\gamma))/(40\rho^*(\gamma))$, where $r^*(., .)$ and $\rho^*(.)$ are defined in Lemma \ref{Lemma:r2-MOM} above. Then, there exists positive constants $C_1, C_2$ and $C_3$ such that, with probability at least $1-\exp(-72K/625)$, it holds true that: 

\begin{align*}
&\norm{\hat{Z}^{\mathrm{RMOM}}_{K, \lambda}-Z^*}_1 \leq \frac{C_1k}{\sqrt{N\theta^2}} \max\left(k\sqrt{\log\left(\frac{ed}{k}\right)} ; \sqrt{K}\right), \\
&\norm{\hat{Z}^{\mathrm{RMOM}}_{K, \lambda} - Z^*}_2 \leq \frac{C_2}{\sqrt{N\theta^2}} \max\left(k\sqrt{\log\left(\frac{ed}{k}\right)} ; \sqrt{K}\right)
\end{align*} and 
\begin{equation*}
 P\cL_{\hat{Z}^{\mathrm{RMOM}}_{K, \lambda}} \leq \frac{C_3}{N\theta} \max\left(k^2 \log\left(\frac{ed}{k}\right); K\right).
\end{equation*}
\end{Theorem}

Since our primary focus is not on $Z^*$ itself, but its unit-norm leading eigenvector $\beta^*$, we are now in the process of providing a result on $\beta^*$.

\begin{Corollary}\label{coro:main_SPCA_RMOM}
Let $\hat{\beta}$ $\in\bR^d$ be a leading unit length eigenvector of $\hat{Z}^{\mathrm{RMOM}}_{K, \lambda}$. Under the conditions of Theorem \ref{theo:main-RMOM-SPCA}, there exists a universal constant $D>0$ such that with probability at least  $1-\exp(-72K/625)$:
\begin{align*}
\norm{\hat{\beta}\hat{\beta}^\top - \beta^*(\beta^*)^\top}_2\leq \frac{D}{\sqrt{N\theta^2}} \max\left(k\sqrt{\log\left(\frac{ed}{k}\right)} ; \sqrt{K}\right).
\end{align*}
\end{Corollary}

In the case where $K\leq k^2\log\left(ed/k\right)$, we get a rate of convergence of magnitude $\sqrt{k^2/(N\theta^2)\log\left(ed/k\right)}$, with no restrictions on the value of $\theta$. This happens with an exponentially large probability depending on the number of groups $K$ even though we only have $\log d$ moments and a dataset that may have been corrupted by an adversary. A similar analysis of a SLOPE regularization of the minmax MOM estimator will lead to a sparsity parameter free choice of $\lambda$.

\section{Proofs}
\label{sec:proofs}
All the proofs from the previous sections -- general excess risk and estimation bounds as well as applications -- are gathered in this section. 
\subsection{Proofs of section \ref{sec:general_bounds}}
We define the regularized excess risk $\cL_Z^{\lambda} := \cL_Z + \lambda (\norm{Z} - \norm{Z^*}) $, and the regularized loss $\ell_Z^{\lambda} := \ell_Z+\lambda \norm{.}$ for all $Z\in\cC$. 

\subsubsection{Proof of Theorem \ref{theo:main-penalized}}

Let $\delta\in (0, 1)$. Let  $A>0$ and $\rho^*>0$ be such that Assumption \ref{ass:A-spars} holds, and assume that $\rho^*>0$ satisfies the $A$-sparsity equation from Definition \ref{def:RERM_A-sparsity}. Let $\gamma:= 1/(3A)$. In the rest of the proof, we write $r^*(.)$ for $r^*_{\mathrm{RERM, G}}(A, ., \delta)$.
Let us define $\cB:=\left\{Z\in \cC:\norm{Z-Z^*}\leq \rho^* \mbox{ and } G(Z-Z^*)\leq r^*(\rho^*)\right\}$. Consider the following event:
\begin{align*}{}
\Omega := \left\{\forall Z\in \cB, \quad |(P-P_N)\cL_Z|\leq \gamma r^*(\rho^*)\right\}.
\end{align*}
By definition of $r^*(.)$, $\Omega$ holds with probability at least $1-\delta$. Let us now prove the statistical bounds announced in Theorem \ref{theo:main-penalized} on the event $\Omega$.  

Suppose that $\hat{Z}\in \cB$. This means that $\norm{\hat{Z}-Z^*}\leq \rho^*$ and $G(\hat Z-Z^*)\leq r^*(\rho^*)$. Moreover, on $\Omega$ it also means that $|(P-P_N)\cL_{\hat{Z}}|\leq \gamma r^*(\rho^*)$, and then:
\begin{align*}
P\cL_{\hat{Z}} & = (P-P_N)\cL_{\hat{Z}} + P_N \cL_{\hat{Z}} \leq \gamma r^*(\rho^*) + P_N \cL_{\hat{Z}} = \gamma r^*(\rho^*) + P_N(\cL^{\lambda}_{\hat{Z}} - \lambda (\norm{\hat{Z}}-\norm{Z^*})) \\ & = \gamma r^*(\rho^*) + P_N\cL^{\lambda}_{\hat{Z}} + \lambda (\norm{Z^*}-\norm{\hat{Z}}) \overset{(i)}{\leq} \gamma r^*(\rho^*) + \lambda \norm{\hat{Z}-Z^*} \leq \gamma r^*(\rho^*) + \lambda \rho^* \overset{(ii)}{\leq} 3\gamma r^*(\rho^*) = \frac{r^*(\rho^*)}{A}
\end{align*} 
where $(i)$ holds since $P_N\cL^{\lambda}_{\hat{Z}}\leq 0$ by definition of $\hat Z$ and $(ii)$ holds because of the choice of $\lambda$ given in (\ref{hyp:lambda}).

Then, if we can show that $\hat{Z}\in \cB$, we will have the desired bounds on $\Omega$. Since we know that $P_N\cL^{\lambda}_{\hat{Z}} \leq 0$, it is sufficient to prove that for any $Z\in \cC \backslash \cB$, $P_N\cL^{\lambda}_Z >0$.

Let $Z\in \cC \backslash \cB$. Because $\cC$ is star-shaped in $Z^*$ and by the regularity properties assumed for $G$, we have the existence of $Z_0\in \partial \cB$, the border of $\cB$, and $\alpha>1$ such that $Z-Z^* = \alpha (Z_0-Z^*)$. The border of $\cB$, that we denoted by $\partial\cB$ is the set of all $Z\in\cC$ such that either $\norm{Z-Z^*}= \rho^*$ and $G(Z-Z^*)\leq r^*(\rho^*)$ or $\norm{Z-Z^*}\leq \rho^*$ and $G(Z-Z^*) = r^*(\rho^*)$.   By linearity of the loss function, we have
$
P_N\cL_Z = \alpha P_N\cL_{Z_0}.
$
Moreover, we have by the triangular inequality that
\begin{align*}
 \norm{Z}-\norm{Z^*} =\norm{\alpha Z_0 - (\alpha - 1)Z^*}  -\norm{Z^*} \geq \alpha\norm{Z_0} - (\alpha - 1)\norm{Z^*} - \norm{Z^*} \geq \alpha (\norm{Z_0} - \norm{Z^*}) 
\end{align*} 
and so
\begin{align}\label{ineq:alphaLZ0}
P_N\cL^{\lambda}_{Z} = P_N\cL_{Z} + \lambda (\norm{Z}-\norm{Z^*}) \geq \alpha P_N\cL_{Z_0} + \lambda\alpha (\norm{Z_0} - \norm{Z^*}) = \alpha P_N\cL^{\lambda}_{Z_0}.
\end{align}
We showed that for any $Z\in \cC \backslash \cB$, there exist $Z_0 \in \partial\cB$ and $\alpha>1$ such that $P_N\cL^{\lambda}_{Z} > \alpha P_N\cL^{\lambda}_{Z_0}$. Hence, we only have to show that $Z\to P_N\cL^{\lambda}_Z$ is positive on the border of $\cB$ to show that it is positive over $\cC \backslash \cB$.

Let $Z_0\in \partial\cB$. Two cases arise: either $\norm{Z_0-Z^*} = \rho^*$ and $G(Z-Z^*)\leq r^*(\rho^*)$, or $\norm{Z_0-Z^*} \leq \rho^*$ and $G(Z-Z^*) = r^*(\rho^*)$.

\underline{First case}: We assume that $\norm{Z_0-Z^*} = \rho^*$ and $G(Z-Z^*)\leq r^*(\rho^*)$, that is $Z_0 \in H_{\rho^*, A}$. Let $V\in H$ be such that $\norm{Z^*-V}\leq \rho^*/20$ and $\Phi \in \partial \norm{.}(V)$. We have:
\begin{align*}
\norm{Z_0} - \norm{Z^*} & \geq \norm{Z_0}- \norm{V} - \norm{Z^*-V}  \\
& \geq \inr{\Phi, Z_0-V} - \norm{Z^*-V} ~~ (\mbox{ since } \Phi \in \partial \norm{.}(V)) \\
& = \inr{\Phi, Z_0-Z^*} - \inr{\Phi, V-Z^*} - \norm{Z^*-V} \\
& \geq \inr{\Phi, Z_0-Z^*} - 2 \norm{Z^*-V} ~~ (\mbox{ since } \inr{\Phi, U}\leq \norm{U} \mbox{ for any } U \in H) \\
& \geq \inr{\Phi, Z_0-Z^*} - \frac{\rho^*}{10}
\end{align*}
This is true for any $\Phi\in \underset{V \in Z^*+\frac{\rho^*}{20}}{\cup}\partial \norm{.}(V) = \Gamma_{Z^*}(\rho^*)$. Then taking the $\sup$ over $\Gamma_{Z^*}(\rho^*)$ gives:
\begin{align*}
\norm{Z_0} - \norm{Z^*}  \geq \sup_{\Phi \in \Gamma_{Z^*}(\rho^*)} \inr{\Phi, Z_0-Z^*} - \frac{\rho^*}{10}
\end{align*}
and then taking the infimum over $H_{\rho^*, A}$ gives:
\begin{align*}
\norm{Z_0} - \norm{Z^*}  & \geq \inf_{Z_0 \in H_{\rho^*, A}} \norm{Z_0} - \norm{Z^*} \geq \inf_{Z_0 \in H_{\rho^*, A}} \sup_{\Phi \in \Gamma_{Z^*}(\rho^*)} \inr{\Phi, Z_0-Z^*}- \frac{\rho^*}{10}   = \Delta(\rho^*, A) - \frac{\rho^*}{10} \geq \frac{7}{10} \rho^* 
\end{align*}
where the last inequality holds since $\rho^*$ is supposed to satisfy the $A$-sparsity equation. Then, we have:

\begin{align*}
P_N\cL^{\lambda}_{Z_0} &= P_N\cL_{Z_0} + \lambda (\norm{Z_0} - \norm{Z^*}) \geq P_N\cL_{Z_0} + \frac{7}{10} \lambda\rho^* = P\cL_{Z_0} - (P-P_N) \cL_{Z_0} +  \frac{7}{10} \lambda\rho^*
\end{align*}
But on $\Omega$, we have $(P-P_N) \cL_{Z_0} \leq \gamma r^*(\rho^*)$ since $Z_0\in\cB$, and we know by definition of $Z^*$ that $P\cL_{Z_0}\geq 0$. Then we conclude that:
\begin{align*}
P_N\cL^{\lambda}_{Z_0} \geq \frac{7}{10} \lambda\rho^* - \gamma r^*(\rho^*) > 0
\end{align*}
where the last inequality is due to the choice of $\lambda$ given in (\ref{hyp:lambda}). \\

\underline{Second case}: Now we assume that $\norm{Z_0-Z^*} \leq \rho^*$ and $G(Z-Z^*)= r^*(\rho^*)$. We have:
\begin{align*}
P_N\cL^{\lambda}_{Z_0} &= P_N\cL_{Z_0} - \lambda(\norm{Z^*}-\norm{Z_0}) \geq P\cL_{Z_0} - (P-P_N)\cL_{Z_0} - \lambda \norm{Z^*-Z_0} \geq P\cL_{Z_0} - (P-P_N)\cL_{Z_0} - \lambda \rho^*
\end{align*}
But we know from Assuption \ref{ass:A-spars} that $P\cL_{Z_0} \geq A^{-1}G(Z_0-Z^*)$, and on $\Omega$ we have $(P-P_N)\cL_{Z_0}\leq \gamma r^*(\rho^*)$. Then we get:
\begin{align*}
P_N\cL^{\lambda}_{Z_0}  &\geq A^{-1}G(Z_0-Z^*) -\gamma r^*(\rho^*) - \lambda \rho^* =  A^{-1}r^*(\rho^*) - \gamma r^*(\rho^*)- \lambda \rho^* >0
\end{align*}
where the last inequality comes from the choice of $\lambda$ given in (\ref{hyp:lambda}).\\

Then, we proved that $P_N\cL^{\lambda}_{Z_0} > 0$ for any $Z_0 \in \partial(\cB)$ and as we said before, this implies that $P_N\cL^{\lambda}_Z$ is positive over $\cC\backslash\cB$. Since $P_N\cL^{\lambda}_{\hat Z} < 0$ we conclude that on $\Omega$, $\hat{Z}$ necessarily belongs to $\cB$, which proves the bounds announced in Theorem \ref{theo:main-penalized}.

\endproof

\subsubsection{Proof of Theorem \ref{theo:main_MOM_loc_excess_risk}}
The proof of this theorem is broken down into two steps. First, we identify an event $\Omega$ on which the estimator $\hat{Z}^{\mathrm{MOM}}_{\mathrm{K}}$ has the desired properties. Then, we show that this event holds with high probability. For the sake of simplicity, in the rest of the proof we write $r^*$ for $r^*_{\mathrm{MOM}, \mathrm{ER}}(\gamma)$ and $\hat{Z}$ for $\hat{Z}^{\mathrm{MOM}}_{K}$. Let $\gamma = 1/6400$, and consider the set $\cC_{\gamma}:=\left\{Z\in\cC : P\cL_Z\leq (r^*)^2\right\}$. Define the event $\Omega_K$ as follows:
\begin{align*}
\Omega_K:=\left\{\forall Z\in\cC_{\gamma}, \exists J\subset [K]: |J| > K/2\mbox{ and } \forall k\in J, |(P_{B_k}-P)\cL_Z|\leq (r^*)^2/4\right\}.
\end{align*}
 We start with showing that on $\Omega_K$, the estimator $\hat Z$ satisfies the excess risk bound announced in Theorem \ref{theo:main_MOM_loc_excess_risk}.

\begin{Lemma}\label{Lemma:main_MOM_loc_excess_risk_control_OmegaK}
On the event $\Omega_K$, $P\cL_{\hat{Z}} \leq r^*$.
\end{Lemma}

\begin{proof}
Let $Z\in \cC \backslash \cC_{\gamma}$. Let $\alpha := (r^*)^{-2}P\cL_Z>1$, and let $Z_0:=Z^*+\alpha^{-1}(Z-Z^*)$. By the star-shaped property of $\cC$, $Z_0\in \cC$, and by linearity of $\ell$, $P\cL_{Z_0} = \alpha^{-1}P\cL_Z=(r^*)^{2}$, so that $Z_0\in \cC_{\gamma}$. Then, on $\Omega_K$, there exists strictly more than $K/2$ blocks $B_k$ on which $|(P_{B_k}-P)\cL_{Z_0}|\leq (r^*)^{2}/4$, that is $P_{B_k}\cL_{Z_0}\geq P\cL_{Z_0} - (r^*)^{2}/4 = (3/4)(r^*)^{2}$ and so $P_{B_k}\cL_Z=\alpha P_{B_k}\cL_{Z_0}\geq \alpha(3/4)(r^*)^2$ because $\alpha>1$. This holds on strictly more than half of the blocks $B_k$, therefore ${\rm Med}(-P_{B_k}\cL_Z:k\in[K])\geq -(3/4)(r^*)^2$ and this holds for all $Z\in \cC \backslash \cC_{\gamma}$, hence, we have
\begin{align}\label{eq:Lemma_main_MOM_loc_excess_risk_control_OmegaK_01}
\sup_{Z\in \cC\backslash \cC_{\gamma}} \mathrm{MOM}_K(\ell_{Z^*} - \ell_Z) \leq -(3/4)(r^*)^{2}.
\end{align}
Moreover, on $\Omega_K$, for $Z\in \cC_{\gamma}$, there exists strictly more than $K/2$ blocks $B_k$ on which $-P_{B_k}\cL_{Z}\leq (r^*)^{2}/4 - P\cL_{Z}\leq (r^*)^{2}/4$, since $P\cL_{Z}\geq 0$ by definition of $Z^*$. Therefore, we have
\begin{align}\label{eq:Lemma_main_MOM_loc_excess_risk_control_OmegaK_02}
\sup_{Z\in \cC_{\gamma}} \mathrm{MOM}_K(\ell_{Z^*} - \ell_Z) \leq (r^*)^{2}/4.
\end{align}
But by definition of $\hat{Z}$, we have:
\begin{align*}
\mathrm{MOM}_K(\ell_{\hat{Z}} - \ell_{Z^*}) & \leq \sup_{Z\in \cC } \mathrm{MOM}_K(\ell_{Z^*} - \ell_Z) = \mathrm{max}\left(\sup_{Z\in \cC_{\gamma} } \mathrm{MOM}_K(\ell_{Z^*} - \ell_Z), \sup_{Z\in \cC \backslash\cC_{\gamma}} \mathrm{MOM}_K(\ell_{Z^*} - \ell_Z)\right)  \leq \frac{(r^*)^{2}}{4}
\end{align*}
that is, $\mathrm{MOM}_K(\ell_{Z^*} - \ell_{\hat{Z}}) \geq -(1/4)(r^*)^{2} > -(3/4)(r^*)^{2}$. From (\ref{eq:Lemma_main_MOM_loc_excess_risk_control_OmegaK_01}) we conclude that, necessarily, $\hat{Z}\in \cC_{\gamma}$, that is, $P\cL_{\hat{Z}} \leq (r^*)^{2}$.
\end{proof}

At this point, we proved that on the event $\Omega_K$, the estimator $\hat{Z}$ satisfies the statistical bounds announced in Theorem \ref{theo:main_MOM_loc_excess_risk}. Now it remains to prove that $\Omega_K$ holds with high probability.

\begin{Lemma}\label{Lemma:main_MOM_loc_er_OmegaK_proba}
Assume that $|\cO|\leq K/100$. Then $\Omega_K$ holds with probability at least $1-\exp(-72K/625)$.
\end{Lemma}
\begin{proof}
Let $\phi:t\in \bR \rightarrow \mathbb{1}_{\left\{t\geq 1\right\}}+2(t-(1/2))\mathbb{1}_{\left\{1/2\leq t\leq 1\right\}}$, so that for any $t\in \bR$, $\mathbb{1}_{\left\{t\geq 1\right\}}\leq \phi(t) \leq \mathbb{1}_{\left\{t\geq 1/2\right\}}$. For $k\in [K]$, let $W_k:=\left\{X_i:i\in B_k\right\}$ and $F_Z(W_k)=(P_{B_k}-P)\cL_Z$. We also define the counterparts of these quantities constructed with the non-corrupted vectors: $\widetilde{W}_k:=\left\{\widetilde{X}_i:i\in B_k\right\}$ and $F_Z(\widetilde{W}_k)=(\widetilde{P_{B_k}}-{P})\cL_Z$, where $\widetilde{P_{B_k}}\cL_Z := (K/N)\sum_{i\in B_k}\cL_Z(\widetilde{X}_i)$. Let $\psi(Z)=\sum_{k\in [K]}\mathbb{1}_{\left\{|F_Z(W_k)|\leq (r^*)^2/4\right\}}$. We show now that, with high probability, if $Z\in \cC_{\gamma}$, then $\psi(Z) > K/2$. In the contaminated framework, it is sufficient to prove that,  with high probability, for all $Z\in \cC_{\gamma}$,
\begin{align}\label{eq:proof_Main_MOM_ER_loc_001}
\sum_{k\in [K]}\mathbb{1}_{\left\{|F_Z(\widetilde{W}_k)| > \frac{(r^*)^2}{4}\right\}} \leq \frac{49K}{100}.
\end{align}
Indeed, consider $Z\in \cC_{\gamma}$ such that (\ref{eq:proof_Main_MOM_ER_loc_001}) holds. Then, there exist at least $(1-49/100)K = (51/100)K$ blocks $B_k$ on which $|F_Z(\widetilde{W}_k)| \leq (r^*)^2/4$. On the other hand, we know that $|\cO|\leq K/100$, so that among the $(51/100)K$ previous blocks, at most $K/100$ contain corrupted data. The other $(50/100)K = K/2$ contain only non-corrupted data, so we have $F_Z(\widetilde{W}_k) = F_Z(W_k)$ on these blocks. We conclude that $\sum_{k\in [K]}\mathbb{1}_{\left\{|F_Z(W_k)| \leq (r^*)^2/4\right\}} > K/2$, that is $\psi(Z) > K/2$, if (\ref{eq:proof_Main_MOM_ER_loc_001}) holds.

Let $Z\in \cC_{\gamma}$. We have:
\begin{align}\label{eq:psi_bound_Xtilde_ERloc_01}
&\sum_{k\in [K]}\mathbb{1}_{\left\{|F_Z(\widetilde{W}_k)| > \frac{(r^*)^2}{4}\right\}} = \sum_{k\in [K]}\left[\mathbb{1}_{\left\{|F_Z(\widetilde{W}_k)| > \frac{(r^*)^2}{4}\right\}} - \bP\left(|F_Z(\widetilde{W}_k)| > \frac{(r^*)^2}{8}\right) + \bP\left(|F_Z(\widetilde{W}_k)| > \frac{(r^*)^2}{8}\right) \right] \nonumber\\
 &= \sum_{k\in [K]}\left(\mathbb{1}_{\left\{|F_Z(\widetilde{W}_k)| > \frac{(r^*)^2}{4}\right\}} - \bE\left[\mathbb{1}_{\left\{|F_Z(\widetilde{W}_k)| > \frac{(r^*)^2}{8}\right\}}\right] \right) + \sum_{k\in [K]} \bP\left(|F_Z(\widetilde{W}_k)| > \frac{(r^*)^2}{8}\right) \nonumber\\ 
&\leq \sum_{k\in [K]}\left(\Phi\left( \frac{4|F_Z(\widetilde{W}_k)|}{(r^*)^2}\right) -  \bE\left[\Phi\left( \frac{4|F_Z(\widetilde{W}_k)|}{(r^*)^2}\right)\right]  \right) + \sum_{k\in [K]} \bP\left(|F_Z(\widetilde{W}_k)| > \frac{(r^*)^2}{8}\right) \nonumber\\ 
&\leq \sup_{Z\in \cC_{\gamma}} \left( \sum_{k\in [K]}\Phi\left( \frac{4|F_Z(\widetilde{W}_k)|}{(r^*)^2}\right) -  \bE\left[\Phi\left( \frac{4|F_Z(\widetilde{W}_k)|}{(r^*)^2}\right)\right]  \right) + \sum_{k\in [K]} \bP\left(|F_Z(\widetilde{W}_k)| > \frac{(r^*)^2}{8}\right).
\end{align}
We start with bounding the last sum in the previous inequality. For each $k\in[K]$, it follows from Markov's inequality and the definition of $r^*$ that
\begin{align*}
\bP\left(|F_Z(\widetilde{W}_k)| > \frac{(r^*)^2}{8}\right)  &\leq \frac{64}{(r^*)^4} \bE\left[F_Z(\widetilde{W}_k)^2\right]
 = \frac{64}{(r^*)^4} \left(\frac{K}{N}\right)\mathrm{Var}(\cL_Z(\widetilde{X})) =  \frac{64}{(r^*)^4}\left(V_K(r^*)\right)^2 \leq \frac{1}{200}.
\end{align*}
Plugging that into (\ref{eq:psi_bound_Xtilde_ERloc_01}), we get:
\begin{align}\label{eq:psi_bound_Xtilde_ERloc_02}
\sum_{k\in [K]}\mathbb{1}_{\left\{|F_Z(\widetilde{W}_k)| > \frac{(r^*)^2}{4}\right\}} \leq \frac{K}{200}+ \sup_{Z\in \cC_{\gamma}} \left( \sum_{k\in [K]}\Phi\left( \frac{4|F_Z(\widetilde{W}_k)|}{(r^*)^2}\right) -  \bE\left[\Phi\left( \frac{4|F_Z(\widetilde{W}_k)|}{(r^*)^2}\right)\right]  \right).
\end{align}
We now we have to bound this last term. Using Mc Diarmind inequality (Theorem 6.2 in \cite{BoucheronLugosiMassart} for $t=12/25$), we get that with probability at least $1-\exp(-72K/625)$, for all $Z\in \cC_{\gamma}$,
\begin{align}\label{eq:proof_Lemma_MOM_er_OmegaK_proba_03}
\sum_{k\in [K]}\phi\left(\frac{4|F_Z(\widetilde{W}_k)|}{(r^*)^{2}}\right) -\bE\phi\left(\frac{4|F_Z(\widetilde{W}_k)|}{(r^*)^{2}}\right)\leq \frac{12}{25}K + \bE\left[\sup_{Z\in \cC_{\gamma}}\sum_{k\in [K]}\phi\left(\frac{4|F_Z(\widetilde{W}_k)|}{(r^*)^{2}}\right) -\bE\phi\left(\frac{4|F_Z(\widetilde{W}_k)|}{(r^*)^{2}}\right)\right].
\end{align}
Let now $\epsilon_1, \ldots, \epsilon_{K}$ be Rademacher variables independent from the $\widetilde{X}_i$'s. By the symmetrization Lemma, we have:
\begin{align}\label{eq:proof_Lemma_MOM_er_OmegaK_proba_04}
\bE\left[\sup_{Z\in \cC_{\gamma}}\sum_{k\in [K]}\phi\left(\frac{4|F_Z(\widetilde{W}_k)|}{(r^*)^{2}}\right) -\bE\left[\phi\left(\frac{4|F_Z(\widetilde{W}_k)|}{(r^*)^{2}}\right)\right]\right]   \leq 2\bE\left[\sup_{Z\in \cC_{\gamma}}\sum_{k\in [K]}\epsilon_k\phi\left(\frac{4|F_Z(\widetilde{W}_k)|}{(r^*)^{2}}\right)\right].
\end{align}
As $\phi$ is $2$-Lipschitz with $\phi(0)=0$, we can use the contraction Lemma (see \cite{LedTal01}, Theorem 4.12) to get that:
\begin{align}\label{eq:proof_Lemma_MOM_er_OmegaK_proba_05}
\bE\left[\sup_{Z\in \cC_{\gamma}}\sum_{k\in [K]}\epsilon_k\phi\left(\frac{4|F_Z(\widetilde{W}_k)|}{(r^*)^{2}}\right)\right] &\leq 8\bE\left[\sup_{Z\in \cC_{\gamma}}\sum_{k\in [K]}\epsilon_k\frac{F_Z(\widetilde{W}_k)}{(r^*)^{2}}\right] = \frac{8}{(r^*)^{2}}\bE\left[\sup_{Z\in \cC_{\gamma}}\sum_{k\in [K]}\epsilon_k(\widetilde{P_{B_k}}-{P})\cL_Z\right].
\end{align}

Now, let $(\sigma_i)_{i=1,\ldots, N}$ be a family of Rademacher variables independent from the $\widetilde{X}_i$'s and the $\epsilon_i$'s. For any $k\in[K]$ and any $i\in[N]$, the variables $\epsilon_k\sigma_i\cL_Z(X_i)$ and $\sigma_i\cL_Z(X_i)$ have the same distribution, so that we get, using the symmetrization Lemma:
\begin{align*}
\bE\left[\sup_{Z\in \cC_{\gamma}}\sum_{k\in [K]}\epsilon_k(\widetilde{P_{B_k}}-P)\cL_Z\right] & \leq 2 \bE\left[\sup_{Z\in \cC_{\gamma}}\frac{K}{N}\sum_{i=1}^N\sigma_i\cL_Z(\widetilde{X}_i)\right] = 2KE(r^*) \leq 2K\gamma (r^*)^2.
\end{align*}
Combining this with (\ref{eq:proof_Lemma_MOM_er_OmegaK_proba_03}), (\ref{eq:proof_Lemma_MOM_er_OmegaK_proba_04}) and (\ref{eq:proof_Lemma_MOM_er_OmegaK_proba_05}), we finally get that, with probability at least $1-\exp(-72K/625)$
\begin{align}\label{eq:proof_Lemma_MOM_er_OmegaK_proba_06}
\sup_{Z\in \cC_{\gamma}}&\sum_{k\in [K]}\phi\left(\frac{4|F_Z(\widetilde{W}_k)|}{(r^*)^{2}}\right) -\bE\left[\phi\left(\frac{4|F_Z(\widetilde{W}_k)|}{(r^*)^{2}}\right)\right] \leq \left(\frac{12}{25}+32\gamma\right)K.
\end{align}
Plugging that into (\ref{eq:psi_bound_Xtilde_ERloc_02}), we conclude that with probability at least $1-\exp(-72K/625)$, one has
\begin{align*}
\sum_{k\in [K]}\mathbb{1}_{\left\{|F_Z(\widetilde{W}_k)| > \frac{(r^*)^2}{4}\right\}} \leq \left(\frac{1}{200}+\frac{12}{25}+32\gamma\right) K \leq \frac{49}{100}K
\end{align*}
from our choice of parameters. This allows to affirm that $\Omega_K$ holds with probability at least $1-\exp(-72K/625)$, which concludes the proof.

\end{proof}

\subsubsection{Proof of Theorem \ref{theo:main_MOM_unloc}.}

The proof is divided into two parts. First, we identify an event $\Omega_K$ on which the estimator has the desired statistical properties. Second, we prove that this event holds with high probability. For the sake of simplicity, we write $\hat{Z}$ for $\hat{Z}^{\mathrm{MOM}}_K$ and $r^*$ for $r^*_{\mathrm{MOM}, L_2}(\gamma)$ with $ \gamma= 1/3200$. Let $0<A<1$ be such that Assumption \ref{ass:curvature_MOM_L2} holds.   We define $\nu =A^2/\gamma$, $\tau = (2A)^{-1}$, $C_{K,A}=\max\left((r^*)^2, \nu K/N\right)$ and $\cB_{K,A}:=\left\{Z\in \cC: \norm{Z-Z^*}_{L_2}\leq \sqrt{C_{K,A}}\right\}$ - where the $L_2$-norm is defined as $Z\to\norm{Z}_{L_2} = \bE[\inr{\widetilde{X}, Z}^2]^{1/2}$. We consider the following event:
\begin{align*}
\Omega_K:=\left\{\forall Z \in \cB_{K,A}, \exists J\subset \{1, \ldots, K\}: |J|>\frac{K}{2}\mbox{ and }\forall k \in J,\left|(P_{B_k}-P)\cL_Z\right|\leq \tau C_{K,A}\right\}.
\end{align*}We show in the next three lemmas that, on $\Omega_K$, $\hat Z$ satisfies the statistical bounds announced in Theorem~.\ref{theo:main_MOM_unloc}. Then the fourth lemma will prove that $\Omega_K$ holds with large probability, the one announced in Theorem~.\ref{theo:main_MOM_unloc}.

\begin{Lemma}\label{Lemma:eta_unl}
If there exists $\eta>0$ such that:
\begin{align}\label{ineq:lm_eta_1_unl}
\sup_{Z\in \cC \backslash \cB_{K,A}} \mathrm{MOM}_K(\ell_{Z^*} - \ell_Z) < -\eta \quad \mbox{ and }\quad & \sup_{Z\in \cB_{K,A}} \mathrm{MOM}_K(\ell_{Z^*} - \ell_Z)  \leq  \eta 
\end{align}
then $\norm{\hat{Z}-Z^*}_{L_2}^2\leq C_{K,A}$.
\end{Lemma}

\begin{proof}
Assume that (\ref{ineq:lm_eta_1_unl}) holds. Then:
\begin{align}\label{ineq:lm_eta_2_unl}
\inf_{Z\in \cC \backslash \cB_{K,A}} \mathrm{MOM}_K(\ell_{Z} - \ell_{Z^*}) >\eta.
\end{align}
Moreover, if we define $Z\to T_K(Z)=\sup_{Z^{\prime}\in \cC } \mathrm{MOM}_K(\ell_Z-\ell_{Z^{\prime}})$, then:
\begin{align}\label{ineq:lm_eta_3_unl}
T_K(Z^*)=\max\left(\sup_{Z\in \cC \backslash \cB_{K,A}} \mathrm{MOM}_K(\ell_{Z^*} - \ell_Z),\sup_{Z\in \cB_{K,A}} \mathrm{MOM}_K(\ell_{Z^*} - \ell_Z)\right) \leq \eta.
\end{align}
By definition of $\hat{Z}$, we  have $T_K(\hat{Z})=\sup_{Z\in \cC } \mathrm{MOM}_K(\ell_{\hat{Z}}-\ell_Z)\leq T_K(Z^*)\leq \eta$. But by (\ref{ineq:lm_eta_2_unl}), any $Z\in \cC \backslash \cB_{K,A}$ satisfies:
\begin{align*}
T_K(Z) \geq \mathrm{MOM}_K(\ell_Z-\ell_{Z^*})\geq \inf_{Z\in \cC \backslash \cB_{K,A}}\mathrm{MOM}_K(\ell_Z-\ell_{Z^*}) \geq \eta
\end{align*}
which allows us to conclude that, necessarily, $\hat{Z}\in \cB_{K,A}$, i.e. $\norm{Z^*-\hat Z}_{L_2}^2\leq C_{K,A}$.
\end{proof}

\begin{Lemma}\label{Lemma:eta_eff_unl}
Assume that $K\geq 100|\cO|$. Then on $\Omega_K$, (\ref{ineq:lm_eta_1_unl}) holds with $\eta = \tau C_{K,A}$. 
\end{Lemma}
\begin{proof}
Let $Z\in \cC$ be such that $\norm{Z-Z^*}_{L_2}>\sqrt{C_{K,A}}$. By the star-shaped property of $\cC$, there exists $Z_0\in \cC$ and $\alpha>1$ such that $\norm{Z_0-Z^*}_{L_2}=\sqrt{C_{K,A}}$ and $Z-Z^*=\alpha(Z_0-Z^*)$. Now, for each block $B_k$ we have by the linearity of the loss function:
\begin{align}\label{ineq:lm_eta_eff_1_unl}
P_{B_k}\cL_Z = \alpha P_{B_k}\cL_{Z_0}.
\end{align}
As $Z_0\in \cB_{K,A}$, on $\Omega_K$ there exist strictly more than $K/2$ blocks on which  $|(P_{B_k}-P)\cL_{Z_0}| \leq \tau C_{K,A}$. Moreover, since $\norm{Z_0-Z^*}_{L_2}=\sqrt{C_{K,A}}$, we get from Assumption \ref{ass:curvature_MOM_L2} that $P\cL_{Z_0}\geq A^{-1}\norm{Z_0-Z^*}_{L_2}^2=A^{-1}C_{K,A}$. Then, on these blocks, $P_{B_k} (\ell_{Z_0}-\ell_{Z^*}) \geq P\cL_{Z_0} - \tau C_{K,A} \geq (A^{-1}-\tau)C_{K,A}$, which implies that $P_{B_k} (\ell_{Z^*}-\ell_{Z_0}) \leq  -(A^{-1}-\tau)C_{K,A} \leq -\tau C_{K,A}$, since we have $\tau = (2A)^{-1}$. From (\ref{ineq:lm_eta_eff_1_unl}) we conclude that, on $\Omega_K$, there exist srictly more than $K/2$ blocks $B_k$ on which $P_{B_k}(\ell_{Z^*}-\ell_{Z})\leq -\alpha \tau C_{K,A}\leq -\tau C_{K,A}$, since $\alpha \geq 1$. This is true for all $Z\in\cC\backslash\cB_{K,A}$; in other words, we have
\begin{align*}
\sup_{Z\in \cC \backslash \cB_{K,A}} \mathrm{MOM}_K(\ell_{Z^*} - \ell_Z) \leq -\tau C_{K,A}
\end{align*} 
Moreover, on $\Omega_K$, for any $Z\in \cB_{K,A}$, there exist stricly more than $K/2$ blocks $B_k$ such that $|(P_{B_k}-P)\cL_{Z}| \leq \tau C_{K,A}$, so that $P_{B_k}(\ell_Z - \ell_{Z^*})\geq -\tau C_{K,A} + P(\ell_Z - \ell_{Z^*}) \geq -\tau C_{K,A}$, since $P(\ell_Z - \ell_{Z^*})\geq 0$ by definition of $Z^*$. Then, we have $P_{B_k}(\ell_{Z^*} - \ell_{Z}) \leq \tau C_{K,A}$ on stricly more than $K/2$ blocks, which implies that $\mathrm{MOM}_K(\ell_{Z^*} - \ell_Z) \leq \tau C_{K,A}$. This being true for any $Z\in \cB_{K,A}$, we conclude that (\ref{ineq:lm_eta_1_unl}) holds with $\eta = \tau C_{K,A}$.
\end{proof}

\begin{Lemma}\label{Lemma:stat_bound_2_MOM_L2}
Grant Assumption \ref{ass:curvature_MOM_L2} and assume that $K\geq 100|\cO|$. On $\Omega_K$, $P\cL_{\hat{Z}} \leq 2 \tau C_{K,A}$.
\end{Lemma}

\begin{proof}
Assume that $\Omega_K$ holds. From Lemmas \ref{Lemma:eta_unl} and \ref{Lemma:eta_eff_unl} , $\norm{\hat{Z}-Z^*}_{L_2}^2\leq C_{K,A}$, that is $\hat Z\in \cB_{K,A}$. Therefore, on strictly more than $K/2$ blocks $B_k$, we have $\left|(P_{B_k}-P)\cL_{\hat{Z}}\right|\leq \tau C_{K,A}$, and then on these blocks:
\begin{align}\label{eq:Lemma_stat_bound_01}
P\cL_{\hat{Z}} \leq P_{B_k}\cL_{\hat{Z}}  +\tau C_{K,A}.
\end{align} 
In addition, by definition of $\hat{Z}$ and (\ref{ineq:lm_eta_3_unl}) (for $\eta = \tau C_{K,A}$):
\begin{align*}
\mathrm{MOM}_K(\ell_{\hat{Z}}-\ell_{Z^*}) \leq \sup_{Z\in\cC}\mathrm{MOM}_K(\ell_{Z^*}-\ell_{Z}) \leq \tau C_{K,A}
\end{align*}
which implies the existence of $K/2$ blocks (at least) on which:
\begin{align}\label{eq:Lemma_stat_bound_02}
 P_{B_k}\cL_{\hat{Z}}\leq \tau C_{K,A}
 \end{align}
 As a consequence, there exist at least one block $B_k$ on which (\ref{eq:Lemma_stat_bound_01}) and (\ref{eq:Lemma_stat_bound_02}) holds simultaneously. On this block, we have:
$
P\cL_{\hat{Z}} \leq \tau C_{K,A} + \tau C_{K,A} = 2 \tau C_{K,A},
$
which concludes the proof.
\end{proof}

At this point, we proved that on the event $\Omega_K$, the estimator $\hat{Z}$ has the statistical properties announced in Theorem \ref{theo:main_MOM_unloc}. In the final lemma, we show that $\Omega_K$ holds with high probability.

\begin{Lemma}\label{Lemma:OmegaK_unl}
Assume that $|\cO|\leq K/100$. Then $\Omega_K$ holds with probability at least $1-\exp(-72K/625)$.
\end{Lemma}
\begin{proof}
Let $\phi:t\in \bR \rightarrow \mathbb{1}_{\left\{t\geq 1\right\}}+2(t-1/2)\mathbb{1}_{\left\{1/2\leq t\leq 1\right\}}$, so that for any $t\in \bR$, $\mathbb{1}_{\left\{t\geq 1\right\}}\leq \phi(t) \leq \mathbb{1}_{\left\{t\geq 1/2\right\}}$. For $k\in [K]$, let $W_k:=\left\{X_i:i\in B_k\right\}$ and $F_Z(W_k)=(P_{B_k}-P)\cL_Z$. We also define the counterparts of these quantities constructed with the non-corrupted vectors: $\widetilde{W}_k:=\left\{\widetilde{X}_i:i\in B_k\right\}$ and $F_Z(\widetilde{W}_k)=(\widetilde{P_{B_k}}-P)\cL_Z$, where $\widetilde{P_{B_k}}\cL_Z := (K/N)\sum_{i\in B_k}\cL_Z(\widetilde{X}_i)$. Let $\psi(Z)=\sum_{k\in [K]}\mathbb{1}_{\left\{|F_Z(W_k)|\leq \tau C_{K,A}\right\}}$. We are now showing that, with high probability, if $Z\in \cB_{K,A}$, then $\psi(Z) > K/2$. In the adversarial corruption setup,  it is enough to prove that the following inequality occurs with high probability: for all $Z\in \cB_{K,A}$,
\begin{align}\label{eq:psi_non_corrupted}
\sum_{k\in [K]}\mathbb{1}_{\left\{|F_Z(\widetilde{W}_k)| > \tau C_{K,A}\right\}} \leq \frac{49K}{100}.
\end{align}
Indeed, consider $Z\in \cC$ such that (\ref{eq:psi_non_corrupted}) holds. Then, there exist at least $(1-49/100)K = (51/100)K$ blocks $B_k$ on which $|F_Z(\widetilde{W}_k)| \leq \tau C_{K,A}$. On the other hand, we know that $|\cO|\leq K/100$, so that among the $(51/100)K$ previous blocks, at most $K/100$ contain corrupted data. The other $(50/100)K = K/2$ contain only non-corrupted data, so we have $F_Z(\widetilde{W}_k) = F_Z(W_k)$ on these blocks. We conclude that $\sum_{k\in [K]}\mathbb{1}_{\left\{|F_Z(W_k)| \leq \tau C_{K,A}\right\}} > K/2$, that is $\psi(Z) > K/2$, if (\ref{eq:psi_non_corrupted}) holds.

Then, we only have to show that (\ref{eq:psi_non_corrupted}) holds uniformly over all $Z\in \cB_{K,A}$ with high probability. This is what we do now. Let $Z\in \cB_{K,A}$. We have:
\begin{align}\label{eq:psi_bound_Xtilde_01}
&\sum_{k\in [K]}\mathbb{1}_{\left\{|F_Z(\widetilde{W}_k)| > \tau C_{K,A}\right\}} = \sum_{k\in [K]}\left[\mathbb{1}_{\left\{|F_Z(\widetilde{W}_k)| > \tau C_{K,A}\right\}} - \bP\left(|F_Z(\widetilde{W}_k)| > \frac{\tau C_{K,A}}{2}\right) + \bP\left(|F_Z(\widetilde{W}_k)| > \frac{\tau C_{K,A}}{2}\right) \right] \nonumber\\
& = \sum_{k\in [K]}\left(\mathbb{1}_{\left\{|F_Z(\widetilde{W}_k)| > \tau C_{K,A}\right\}} - \bE\left[\mathbb{1}_{\left\{|F_Z(\widetilde{W}_k)| > \frac{\tau C_{K,A}}{2}\right\}}\right] \right) + \sum_{k\in [K]} \bP\left(|F_Z(\widetilde{W}_k)| > \frac{\tau C_{K,A}}{2}\right) \nonumber\\ 
 &\leq \sum_{k\in [K]}\left(\Phi\left( \frac{|F_Z(\widetilde{W}_k)|}{\tau C_{K,A}}\right) -  \bE\left[\Phi\left( \frac{|F_Z(\widetilde{W}_k)|}{\tau C_{K,A}}\right)\right]  \right) + \sum_{k\in [K]} \bP\left(|F_Z(\widetilde{W}_k)| > \frac{\tau C_{K,A}}{2}\right)  \nonumber\\ 
&\leq \sup_{Z\in \cB_{K,A}} \left( \sum_{k\in [K]}\Phi\left( \frac{|F_Z(\widetilde{W}_k)|}{\tau C_{K,A}}\right) -  \bE\left[\Phi\left( \frac{|F_Z(\widetilde{W}_k)|}{\tau C_{K,A}}\right)\right]  \right) + \sum_{k\in [K]} \bP\left(|F_Z(\widetilde{W}_k)| > \frac{\tau C_{K,A}}{2}\right).
\end{align}
We start with bounding the last sum in the previous inequality. For each $k\in[K]$, it follows from Markov's inequality, the definition of $C_{K,A}$ and the linearity of the loss function that
\begin{align*}
&\bP\left(|F_Z(\widetilde{W}_k)| > \frac{\tau C_{K,A}}{2}\right)  \leq \frac{4}{(\tau C_{K,A})^2} \bE\left[F_Z(\widetilde{W}_k)^2\right]  = \frac{4}{(\tau C_{K,A})^2} \left(\frac{K}{N}\right) {\rm Var}(\cL_Z(\widetilde{X}))\\ 
&\leq \frac{4}{(\tau C_{K,A})^2} \left(\frac{K}{N}\right) \bE[\cL_Z(\widetilde{X})^2]= \frac{4}{(\tau C_{K,A})^2} \frac{K}{N} \norm{Z-Z^*}_{L_2}^2 \leq \frac{4}{(\tau C_{K,A})^2} \frac{K}{N} C_{K,A} \leq \frac{4}{\tau^2\nu} = \frac{1}{200}.
\end{align*}
 Plugging the latter result into (\ref{eq:psi_bound_Xtilde_01}), we get:
\begin{align}\label{eq:psi_bound_Xtilde_02}
\sum_{k\in [K]}\mathbb{1}_{\left\{|F_Z(\widetilde{W}_k)| > \tau C_{K,A}\right\}} \leq \frac{K}{200} + \sup_{Z\in \cB_{K,A}} \left( \sum_{k\in [K]}\Phi\left( \frac{|F_Z(\widetilde{W}_k)|}{\tau C_{K,A}}\right) -  \bE\left[\Phi\left( \frac{|F_Z(\widetilde{W}_k)|}{\tau C_{K,A}}\right)\right]  \right).
\end{align}
We now have to bound this last term. Using Mc Diarmind inequality (Theorem 6.2 in \cite{BoucheronLugosiMassart} for taking $t = 12/25$), we get that with probability at least $1-\exp(-72K/625)$, for all $Z\in \cB_{K,A}$,
\begin{align*}
\sum_{k\in [K]}\phi\left(\frac{|F_Z(\widetilde{W}_k)|}{\tau C_{K,A}}\right) -\bE\phi\left(\frac{|F_Z(\widetilde{W}_k)|}{\tau C_{K,A}}\right)\leq \frac{12}{25}K + \bE\left[\sup_{Z\in \cB_{K,A}}\sum_{k\in [K]}\phi\left(\frac{|F_Z(\widetilde{W}_k)|}{\tau C_{K,A}}\right) -\bE\phi\left(\frac{|F_Z(\widetilde{W}_k)|}{\tau C_{K,A}}\right)\right].
\end{align*}
Let now $\epsilon_1, \ldots, \epsilon_{K}$ be Rademacher variables independent from the $\widetilde{X}_i$'s. By the symmetrization Lemma, we have:
\begin{align*}
\bE\left[\sup_{Z\in \cB_{K,A}}\sum_{k\in [K]}\phi\left(\frac{|F_Z(\widetilde{W}_k)|}{\tau C_{K,A}}\right) -\bE\left[\phi\left(\frac{|F_Z(\widetilde{W}_k)|}{\tau C_{K,A}}\right)\right] \right]  \leq 2\bE\left[\sup_{Z\in \cB_{K,A}}\sum_{k\in [K]}\epsilon_k\phi\left(\frac{|F_Z(\widetilde{W}_k)|}{\tau C_{K,A}}\right)\right]
\end{align*}
As $\phi$ is $2$-Lipschitz with $\phi(0)=0$, we can use the contraction Lemma (see \cite{LedTal01}, Theorem 4.3) to get that:
\begin{align*}
\bE\left[\sup_{Z\in \cB_{K,A}}\sum_{k\in [K]}\epsilon_k\phi\left(\frac{|F_Z(\widetilde{W}_k)|}{\tau C_{K,A}}\right)\right] &\leq 2\bE\left[\sup_{Z\in \cB_{K,A}}\sum_{k\in [K]}\epsilon_k\frac{F_Z(\widetilde{W}_k)}{\tau C_{K,A}}\right] = 2 \bE\left[\sup_{Z\in \cB_{K,A}}\sum_{k\in [K]}\epsilon_k\frac{(\widetilde{P_{B_k}}-P)\cL_Z}{\tau C_{K,A}}\right]
\end{align*}
Now, let $(\sigma_i)_{i=1,\ldots, K}$ be a family of Rademacher variables independant from the $\widetilde{X}_i$'s and the $\epsilon_k$'s. Using the symmetrization Lemma one more time, we get
\begin{align*}
\bE\left[\sup_{Z\in \cB_{K,A}}\sum_{k\in [K]}\epsilon_k\frac{(\widetilde{P_{B_k}}-P)\cL_Z}{C_{K,A}}\right] \leq 2 \bE\left[\sup_{Z\in \cB_{K,A}}\frac{K}{N}\sum_{i=1}^N\sigma_i\frac{\cL_Z(\widetilde{X}_i)}{C_{K,A}}\right].
\end{align*}
To bound this last term, we consider two cases: either $C_{K,A}=(r^*)^2$ or $C_{K,A}=\nu K/N$. In the first case, by definition of $r^*$ we have:
\begin{align*}
\bE\left[\sup_{Z\in \cB_{K,A}}\sum_{i=1}^N\sigma_i\frac{\cL_Z(\widetilde{X}_i)}{C_{K,A}}\right] \leq \frac{1}{C_{K,A}}\gamma (r^*)^2N = \gamma N.
\end{align*}
In the second case, we decompose the supremum into two parts:
\begin{align*}
\sup_{Z\in \cB_{K,A}}\sum_{i=1}^N\sigma_i\cL_Z(\widetilde{X}_i) = \max &\left(\sup_{Z\in \cB_{K,A}:\norm{Z-Z^*}_{L_2}\leq r^*}\sum_{i=1}^N\sigma_i\cL_Z(\widetilde{X}_i), \sup_{Z\in \cB_{K,A}:r^*\leq \norm{Z-Z^*}_{L_2}\leq \sqrt{\frac{\nu K}{N}}}\sum_{i = 1}^N\sigma_i \cL_Z(\widetilde{X}_i)\right).
\end{align*}
Let $Z\in \cB_{K,A}$ be such that $r^*\leq \norm{Z-Z^*}_{L_2} \leq \sqrt{\frac{\nu K}{N}}$. Since $\cC$ is star-shapped in $Z^*$, there exists $Z_0\in \cC$ such that $\norm{Z_0-Z^*}_{L_2}=r^*$ and $Z-Z^*=\kappa(Z_0-Z^*)$ for some $\kappa\geq 1$, so that $\kappa =  \frac{\norm{Z-Z^*}_{L_2}}{\norm{Z_0-Z^*}_{L_2}}\leq \sqrt{\frac{\nu K}{N}}\frac{1}{r^*}$. Moreover, we have by linearity of $\cL$ that $\cL_{Z_0} = \kappa \cL_Z$. Therefore, we obtain
\begin{align*}
\sup_{Z\in \cB_{K,A}:r^*\leq \norm{Z-Z^*}_{L_2} \leq \sqrt{\frac{\nu K}{N}}}\sum_{i=1}^N\sigma_i\cL_Z(\widetilde{X}_i) &\leq \sup_{1\leq \kappa\leq \frac{1}{r^*}\sqrt{\frac{\nu K}{N}}}\sup_{Z_0\in \cB_{K,A}:\norm{Z_0-Z^*}_{L_2}\leq r^* }\sum_{i=1}^N\sigma_i\kappa\cL_{Z_0}(\widetilde{X}_i) \\
& =  \sqrt{\frac{\nu K}{N}}\frac{1}{r^*}\sup_{Z_0\in \cB_{K,A}:\norm{Z_0-Z^*}_{L_2}\leq r^* }\sum_{i=1}^N\sigma_i\cL_{Z_0}(\widetilde{X}_i).
\end{align*}
Since $C_{K,A} = \nu K/N\geq (r^*)^2$, we get, using the definition of $r^*$:
\begin{align*}
\bE\left[\sup_{Z\in \cB_{K,A}}\sum_{i=1}^N\sigma_i\cL_Z(\widetilde{X}_i)\right]& \leq \sqrt{\frac{\nu K}{N}}\frac{1}{r^*} \bE\left[\sup_{Z\in \cB_{K,A}:\norm{Z-Z^*}_{L_2}\leq r^* }\sum_{i=1}^N\sigma_i\cL_{Z}(\widetilde{X}_i)\right] \leq \sqrt{\frac{\nu K}{N}}\frac{1}{r^*} \gamma (r^*)^2N \leq C_{K,A}\gamma N.
\end{align*}
Finally, we get that whatever the value of $C_{K,A}$ is:
\begin{align*}
\bE\left[\sup_{Z\in \cB_{K,A}}\sum_{i=1}^N\sigma_i\frac{\cL_Z(\widetilde{X}_i)}{C_{K,A}}\right] \leq \gamma N.
\end{align*}
Combining all these inequalities, we finally get that, with probability at least $1-\exp(-12K/625)$, for all $Z\in \cB_{K,A}$,
\begin{align*}
\sum_{k\in [K]}\mathbb{1}_{\left\{|F_Z(\widetilde{W}_k)| > \tau C_{K,A}\right\}} & \leq  \frac{12}{25}K + \frac{1}{200}K + \frac{8\gamma}{\tau} K \leq \frac{49}{400}K
\end{align*}
From our choice of parameters. This concludes the proof.
\end{proof}

\subsubsection{Proof of Theorem \ref{coro:main-MOM-loc}}
The proof of this theorem follows the same lines as the one of the last Theorem~\ref{theo:main_MOM_loc_excess_risk} and  \ref{theo:main_MOM_unloc}: we start with identifying an event on which our estimator has the desired properties, and then we prove that this event holds with large probability. 

For the sake of simplicity, we write $\hat{Z}$ for $\hat{Z}^{\mathrm{MOM}}_K$ and $r^*$ for $r^*_{\mathrm{MOM, G}}(\gamma)$ for $\gamma = 1/6400$. Consider $A$ and $G:H\rightarrow \bR$ such that Assumption \ref{ass:curvature-MOM} holds. Define $\cC_{\gamma, G} := \left\{Z\in \cC:G(Z-Z^*)\leq (r^*)^2\right\}$.  We consider the following event:
\begin{align*}
\Omega_K = \left\{\forall Z \in \cC_{\gamma, G},\exists J\subset [N]:|J|> K/2 \mbox{ and } \forall k\in J,\left|(P_{B_k}-P)\cL_Z\right|\leq \frac{(r^*)^2}{4} \right\}.
\end{align*}

We first show that on the event $\Omega_K$,  $\hat Z$ satisfies the statistical bounds announced in Theorem \ref{coro:main-MOM-loc}.

\begin{Lemma}\label{Lemma:coro_main-MOM-loc_01}
If there exists $\eta >0$ such that 
\begin{align}\label{eq:coro_main-MOM-loc_01}
\underset{Z\in\cC\backslash\cC_{\gamma, G}}{\mathrm{sup}} \mathrm{MOM}_K(\ell_{Z^*}-\ell_{Z}) < -\eta \quad \mbox{and} \quad
\underset{Z\in \cC_{\gamma, G}}{\mathrm{sup}} \mathrm{MOM}_K(\ell_{Z^*}-\ell_{Z}) \leq \eta 
\end{align}
then $G(\hat{Z}-Z^*)\leq (r^*)^2$.
\end{Lemma}

\begin{proof}
Assume that (\ref{eq:coro_main-MOM-loc_01}) holds. Then:
\begin{align}\label{eq:coro_main-MOM-loc_02}
\inf_{Z\in \cC \backslash \cC_{\gamma, G}} \mathrm{MOM}_K(\ell_{Z} - \ell_{Z^*}) >\eta.
\end{align}
Moreover, define $Z\to T_K(Z)=\sup_{Z^{\prime}\in \cC } \mathrm{MOM}_K(\ell_Z-\ell_{Z^{\prime}})$, we have
\begin{align}\label{ineq:Lemma_eta_MOM_G}
T_K(Z^*)=\max\left(\sup_{Z\in \cC \backslash \cC_{\gamma, G}} \mathrm{MOM}_K(\ell_{Z^*} - \ell_Z),\sup_{Z\in \cC_{\gamma, G}} \mathrm{MOM}_K(\ell_{Z^*} - \ell_Z)\right) \leq \eta
\end{align}
and, by definition of $\hat Z$, we also have $T_K(\hat{Z})=\sup_{Z\in \cC } \mathrm{MOM}_K(\ell_{\hat{Z}}-\ell_Z)\leq \sup_{Z\in \cC } \mathrm{MOM}_K(\ell_{Z^*}-\ell_Z)= T_K(Z^*)\leq \eta$. However, by (\ref{eq:coro_main-MOM-loc_02}), any $Z\in \cC \backslash \cC_{\gamma, G}$ must satisfy
\begin{align*}
T_K(Z) \geq \mathrm{MOM}_K(\ell_Z-\ell_{Z^*})\geq \inf_{Z\in \cC \backslash \cC_{\gamma, G}}\mathrm{MOM}_K(\ell_Z-\ell_{Z^*}) > \eta.
\end{align*}
Therefore, we necessarily have $\hat{Z}\in \cC\cap\cC_{\gamma, G}$, that is $G(\hat{Z}-Z^*)\leq (r^*)^2$.
\end{proof}

\begin{Lemma}\label{Lemma:coro_main-MOM-loc_02}
Assume that $A<2$. On the event $\Omega_K$,  (\ref{eq:coro_main-MOM-loc_01}) holds with $\eta = (r^*)^2/4$.
\end{Lemma}

\begin{proof}
Let $Z$ be such that $G(Z-Z^*)>(r^*)^2$. By the star-shaped property of $\cC$ and the regularity property of $G$, there exist $Z_0\in \partial \cC_{\gamma, G}$ and $\alpha > 1$ such that $Z = Z^* + \alpha (Z_0 - Z^*)$. Since $G(Z_0-Z^*)=(r^*)^2$, we have by Assumption \ref{ass:curvature-MOM} that $P\cL_{Z_0} \geq A^{-1}G(Z_0-Z^*)$. Moreover, on $\Omega_K$, there are at least $K/2$ blocks $B_k$ on which $\left|(P_{B_k}-P)\cL_{Z_0}\right|\leq (r^*)^2/4 $ and so $P_{B_k}\cL_{Z_0}\geq P\cL_{Z_0} - (r^*)^2/4 \geq A^{-1}G(Z_0-Z^*) - (r^*)^2/4 \geq (r^*)^2/4$ since we assumed that $A<2$. Now, by linearity of the loss function, we have on these blocks
\begin{align*}
P_{B_k}\cL_{Z} = \alpha P_{B_k}\cL_{Z_0} \geq \alpha (r^*)^2/4 > (r^*)^2/4.
\end{align*}
We conclude that $\mathrm{MOM}_K(\ell_{Z^*}-\ell_Z) < -(r^*)^2/4$. This being true for any $Z\in \cC\backslash \cC_{\gamma, G}$ we have:
\begin{align*}
\underset{Z\in\cC\backslash\cC_{\gamma, G}}{\mathrm{sup}} \mathrm{MOM}_K(\ell_{Z^*}-\ell_{Z}) \leq -\frac{(r^*)^2}{4}.
\end{align*}This shows the left-hand side inequality of (\ref{eq:coro_main-MOM-loc_01}) for $\eta = (r^*)^2/4$.

Next, let $Z\in \cC$ be such that $G(Z-Z^*)\leq(r^*)^2$. On $\Omega_K$, there are at least $K/2$ blocks $B_k$ on which $\left|(P_{B_k}-P)\cL_{Z_0}\right|\leq (r^*)^2/4 $, that is $-P_{B_k}\cL_Z\leq (r^*)^2/4 - P\cL_{Z} \leq (r^*)^2/4$ since $P\cL_Z\geq 0$ by definition of $Z^*$. Then, $\mathrm{MOM}_K(\ell_{Z^*}-\ell_Z) \leq (r^*)^2/4$. This holds  for all $Z\in  \cC_{\gamma, G}$, in other words, the right-hand side inequality of  (\ref{eq:coro_main-MOM-loc_01}) holds for $\eta = (r^*)^2/4$ and this concludes the proof.
\end{proof}

\begin{Lemma}\label{Lemma:coro_main-MOM-loc_03}
Assume the conditions of Theorem \ref{coro:main-MOM-loc} are met. Then, on $\Omega_K$, $P\cL_{\hat{Z}} \leq (r^*)^2/2$.
\end{Lemma}

\begin{proof}
From Assumption \ref{ass:curvature-MOM} combined with the fact that $A<2$, we have from Lemmas \ref{Lemma:coro_main-MOM-loc_01} and \ref{Lemma:coro_main-MOM-loc_02} that $G(\hat{Z}-Z^*) \leq (r^*)^2$. Then on $\Omega_K$ there exist strictly more than $K/2$ blocks $B_k$ on which $\left|(P_{B_k}-P)\cL_{\hat{Z}}\right|\leq (r^*)^2/4$, that is:
\begin{align}\label{eq:coro_MOM_loc_G_01}
P\cL_{\hat{Z}} \leq P_{B_k}\cL_{\hat{Z}} + \frac{(r^*)^2}{4}
\end{align}
Moreover, by (\ref{ineq:Lemma_eta_MOM_G}) and by definition of $\hat{Z}$, we have:
\begin{align*}
\mathrm{MOM}_K(\ell_{\hat{Z}}-\ell_{Z^*}) \leq \sup_{Z\in \cC} \mathrm{ MOM}_K(\ell_{\hat{Z}}-\ell_{Z}) \leq \sup_{Z\in \cC} \mathrm{ MOM}_K(\ell_{Z^*}-\ell_{Z}) = T_K(Z^*) \leq \eta = \frac{(r^*)^2}{4}.
\end{align*}
As a consequence, there exist at least $K/2$ blocks $B_k$ on which $P_{B_k}(\ell_{\hat{Z}}-\ell_{Z^*})\leq (r^*)^2/4 $, that is:
\begin{align}\label{eq:coro_MOM_loc_G_02}
P_{B_q}\cL_{\hat{Z}} \leq \frac{(r^*)^2}{4}.
\end{align}
So there must be at least one block $B_{k_0}$ on which (\ref{eq:coro_MOM_loc_G_01}) and (\ref{eq:coro_MOM_loc_G_02}) hold simultaneously. On this block, we have:
\begin{align*}
P\cL_{\hat{Z}} \leq P_{B_{k_0}}\cL_{\hat{Z}} + \frac{(r^*)^2}{4} \leq \frac{(r^*)^2}{4} + \frac{(r^*)^2}{4} = \frac{(r^*)^2}{2}.
\end{align*}
\end{proof}

At this stage of the proof, we have shown that on the event $\Omega_K$, the estimator $\hat{Z}$ has the statistical bounds announced in Theorem \ref{coro:main-MOM-loc}. The final ingredient is to show that, under the conditions of Theorem \ref{coro:main-MOM-loc}, $\Omega_K$ holds with exponentially large probability. This is the purpose of the next result that can be proved using the same proof as the one of Lemma \ref{Lemma:main_MOM_loc_er_OmegaK_proba}.

\begin{Lemma}\label{Lemma:coro_main-MOM-loc_04}
Assume the conditions of Theorem \ref{coro:main-MOM-loc} are met, with $A<2$. Then $\Omega_K$ holds with probability at least $1-\exp(-72K/625)$.
\end{Lemma}

\subsubsection{Proof of Theorem \ref{theo:main_RERM_er}}
The proof is structured in the same way as the previous ones: we identify an event on which $\hat{Z}^{\mathrm{RMOM}}_{K, \lambda}$ has the desired statistical properties, then we show that this event holds with high probability. Let $\gamma = 1/32000$. Consider $\rho^*>0$ such that $\rho^*$ satisfies the sparsity equation of Definition \ref{def:sparsity_RMOM_er}. For the sake of simplicity, all along this proof we write $\hat{Z}$ for $\hat{Z}^{\mathrm{RMOM}}_{K, \lambda}$ and $r^*_b:=r^*_{\mathrm{RMOM, ER}}(\gamma, b\rho^*)$ for both $b\in\{1,2\}$. For $b\in \{1, 2\}$, we define $\cB_b:=\left\{Z\in \cC:P\cL_Z\leq (r^*_b)^2\mbox{ and }\norm{Z-Z^*}\leq b\rho^* \right\}$. Then we define:
\begin{align*}
\Omega_{K} = \left\{\forall b \in \{1, 2\}, \forall Z \in \cB_b, \exists J\subset[K], |J|>K/2, \forall k\in J, \left|(P_{B_k}-P)\cL_Z\right|\leq \frac{(r^*_b)^2}{20}\right\}.
\end{align*}
Finally, we consider $\lambda := (11/(40\rho^*))(r^*_2)^2$. We begin the proof by showing that on $\Omega_K$, $\hat Z$ has the statistical properties announced in Theorem \ref{theo:main_RERM_er}.

\begin{Lemma}\label{Lemma:proof_main_RMOM_ER-loc-00}
If there exists $\eta > 0$ such that
\begin{align}\label{eq:Lemma_eta_RMOM_ER-loc-10}
\underset{Z\in\cC\backslash\cB_{2}}{\mathrm{sup}} \mathrm{ MOM}_K(\ell_{Z^*}-\ell_{Z}) + \lambda\left(\norm{Z^*} - \norm{Z}\right) < -\eta 
\end{align}
and
\begin{align}\label{eq:Lemma_eta_RMOM_ER-loc-20}
\sup_{Z\in\cC} \mathrm{ MOM}_K(\ell_{Z^*}-\ell_{Z}) + \lambda\left(\norm{Z^*} - \norm{Z}\right)  \leq \eta 
\end{align}then $P\cL_{\hat{Z}}\leq (r^*_2)^2$ and $\norm{\hat Z-Z^*}\leq 2\rho^*$.
\end{Lemma}

\begin{proof}
For $Z\in \cC$, define $S(Z)=\sup_{Z^\prime\in \cC}\mathrm{MOM}_K(\ell_Z-\ell_{Z^\prime})+\lambda(\norm{Z}-\norm{Z^\prime})$. For all $Z\in \cC\backslash \cB_2$ we have:
\begin{align*}
S(Z) \geq \mathrm{MOM}_K(\ell_Z-\ell_Z^*)+\lambda(\norm{Z}-\norm{Z^*}) \geq \inf_{Z\in \cC \backslash \cB_2} \mathrm{MOM}_K(\ell_Z-\ell_Z^*)+\lambda(\norm{Z}-\norm{Z^*}) > \eta
\end{align*}
since (\ref{eq:Lemma_eta_RMOM_ER-loc-10}) holds. Moreover, we have by definition of $\hat{Z}$:
\begin{align*}
S(\hat{Z}) \leq S(Z^*) = \sup_{Z\in \cC} \mathrm{MOM}_K(\ell_{Z^*} - \ell_Z) + \lambda (\norm{Z^*} - \norm{Z}) \leq  \eta 
\end{align*}
since (\ref{eq:Lemma_eta_RMOM_ER-loc-20}) holds. This shows that necessarily $\hat{Z}\in \cB_2$.
\end{proof}

We are now looking for $\eta>0$ such that (\ref{eq:Lemma_eta_RMOM_ER-loc-10}) and (\ref{eq:Lemma_eta_RMOM_ER-loc-20}) hold, which the following Lemma allows us to do. 

\begin{Lemma}\label{Lemma:proof_main_RMOM_ER-loc-01}
Under the assumptions of Theorem \ref{theo:main_RERM_er} and on the event event $\Omega_K$,  (\ref{eq:Lemma_eta_RMOM_ER-loc-10}) and (\ref{eq:Lemma_eta_RMOM_ER-loc-20}) hold with $\eta = 19(r_2^*)^2/50$.
\end{Lemma}

\begin{proof} Let $b\in\{1,2\}$.
Let $Z\in\cC\backslash\cB_b$. By the star-shaped property of $\cC$, there exist $Z_0\in \partial\cB_b$ and $\alpha>1$ such that $Z=Z^*+\alpha(Z_0-Z^*)$.  As a consequence, by linearity of the loss function and convexity of the regularization norm, for all $k\in[K]$ we have 
\begin{align}\label{eq:top_proof_lemma_6_12}
&P_{B_k}\cL^{\lambda}_{Z} = P_{B_k}\cL_{Z} + \lambda (\norm{Z} - \norm{Z^*}) = \alpha P_{B_k}\cL_{Z_0} + \lambda (\norm{\alpha Z_0 + (1-\alpha)Z^*} - \norm{Z^*})\nonumber\\ 
&\geq \alpha P_{B_k}\cL_{Z_0} + \lambda \alpha(\norm{Z_0} - \norm{Z^*}) = \alpha P_{B_k}\cL^{\lambda}_{Z_0}.
\end{align}

Now, since $Z_0\in\partial\cB_b$, we have either \textit{a)} $P\cL_{Z_0} = (r^*_b)^2$ and $\norm{Z_0-Z^*} < b\rho^*$ or \textit{b)} $P\cL_{Z_0} < (r^*_b)^2$ and $\norm{Z_0-Z^*} = b\rho^*$.

In the first case \textit{a)}, on $\Omega_K$, there are at least $K/2$ blocks $B_k$ on which $P_{B_k}\cL_{Z_0}\geq P\cL_{Z_0} - (r^*_b)^2/20 = (19/20)(r^*_b)^2$. Therefore,  on these blocs, we have
\begin{align}\label{eq:Lemma_nu_ER_01}
  &P_{B_k}\cL_{Z_0}^{\lambda}  = P_{B_k}\cL_{Z_0} + \lambda (\norm{Z_0}-\norm{Z^*}) \geq \frac{19}{20}(r^*_b)^2 - \lambda \norm{Z_0-Z^*}\nonumber\\ 
  &\geq \frac{19}{20}(r^*_b)^2 - \lambda b\rho^*  = \frac{19}{20}(r^*_b)^2 - \frac{11b}{40} (r^*_2)^2 \geq 
\left\{
\begin{array}{cc}
2(r_2^*)^2/5 & \mbox{ for } b=2\\
(r_2^*)^2/5 & \mbox{ for } b=1.
\end{array}
\right.
\end{align}where we used in the case $b=1$ that $r_1^*\geq r_2^*/\sqrt{2}$ thanks to Proposition~\ref{prop:prop_local_fixed_point_appendix} from the Appendix.

In the second case \textit{b)}, we have $Z_0 \in \bar{H}_{b\rho^*}$ from Definition \ref{def:sparsity_RMOM_er}. Since the sparsity equation holds for $\rho=\rho^*$, it also holds for $\rho=b\rho^*$ (see Proposition~\ref{prop:prop_sparsity_equation} in the Appendix). Let $V\in H$ be such that $\norm{Z^*-V}\leq b\rho^*/20$ and $\Phi \in \partial \norm{.}(V)$. We have:
\begin{align*}
\norm{Z_0} - \norm{Z^*} & \geq \norm{Z_0}- \norm{V} - \norm{Z^*-V}  \\
& \geq \inr{\Phi, Z_0-V} - \norm{Z^*-V} ~~ (\mbox{ since } \Phi \in \partial \norm{.}(V)) \\
& = \inr{\Phi, Z_0-Z^*} - \inr{\Phi, V-Z^*} - \norm{Z^*-V} \\
& \geq \inr{\Phi, Z_0-Z^*} - 2 \norm{Z^*-V} ~~ (\mbox{ since } \inr{\Phi,U}\leq \norm{U} \mbox{ for any } U \in H) \\
& \geq \inr{\Phi, Z_0-Z^*} - \frac{b\rho^*}{10}.
\end{align*}
This is true for any $\Phi\in \underset{V \in Z^*+b\rho^*/20}{\cup}\partial \norm{.}(V) = \Gamma_{Z^*}(b\rho^*)$. Then taking the $\sup$ over $\Gamma_{Z^*}(b\rho^*)$ gives:
\begin{align*}
\norm{Z_0} - \norm{Z^*}  \geq \sup_{\Phi \in \Gamma_{Z^*}(b\rho^*)} \inr{\Phi, Z_0-Z^*} - \frac{b\rho^*}{10}
\end{align*}
and then taking the infimum over $\bar{H}_{b\rho^*}$ gives:
\begin{align}\label{eq:Lemma_nu_ER_02}
\norm{Z_0} - \norm{Z^*}  & \geq \inf_{Z_0 \in \bar{H}_{b\rho^*}} \norm{Z_0} - \norm{Z^*} \geq \inf_{Z_0 \in \bar{H}_{b\rho^*}} \sup_{\Phi \in \Gamma_{Z^*}(2\rho^*)} \inr{\Phi, Z_0-Z^*} - \frac{b\rho^*}{10}  = \Delta(b\rho^*) - \frac{b\rho^*}{10} \geq \frac{7}{10} b\rho^*
\end{align}
where the last inequality holds since $b\rho^*$ satisfies the sparsity equation. Then, $\lambda(\norm{Z_0}-\norm{Z^*})\geq (7/10)\lambda b\rho^* = (77/400)b(r^*_2)^2$. Now, since $Z_0\in \cB_b$, on $\Omega_K$ there exist at least $K/2$ blocks $B_k$ such that $|(P_{B_k}-P)\cL_{Z_0}|\leq (r^*_b)^2/20 $ and so $P_{B_k} \cL_{Z_0}\geq (r^*_b)^2/20$ - because $P\cL_{Z_0}\geq 0$. Therefore, on the very same  blocks,
\begin{align}\label{eq:Lemma_nu_ER_03}
 P_{B_k} \cL^{\lambda}_{Z_0} =  P_{B_k} \cL_{Z_0}+\lambda (\norm{Z_0}-\norm{Z^*}) \geq  -\frac{1}{20}  (r^*_b)^2 + \frac{77}{400}b(r^*_2)^2\geq \left\{
\begin{array}{cc}
134(r_2^*)^2/400 & \mbox{ for } b=2\\
29(r_2^*)^2/400 & \mbox{ for } b=1
\end{array}
\right.
\end{align}where we used that $r_1^*\leq r_2^*$ (see Proposition~\ref{prop:prop_local_fixed_point_appendix} in the Appendix). As a consequence, it follows from \eqref{eq:top_proof_lemma_6_12}, the fact that $\alpha>1$,  \eqref{eq:Lemma_nu_ER_01} and \eqref{eq:Lemma_nu_ER_03} for $b=2$ that for all $Z\in\cC\backslash\cB_2$, on more than $K/2$ blocks $B_k$: $P_{B_k}\cL_Z^\lambda\geq (134/400)(r^*_2)^2$ and so \eqref{eq:Lemma_eta_RMOM_ER-loc-10} holds for $\eta \leq  (134/400)(r^*_2)^2$.

Let us now turn to Equation~\eqref{eq:Lemma_eta_RMOM_ER-loc-20}. Let $Z\in\cB_1$. On $\Omega_K$ there exist at least $K/2$ blocks $B_k$ such that $|(P_{B_k}-P)\cL_{Z}|\leq (r^*_1)^2/20$. On these blocks $B_k$, all $P_{B_k} \cL^{\lambda}_Z$'s are such that
\begin{align}\label{eq:Lemma_nu_ER_04}
 &P_{B_k} \cL^{\lambda}_Z = P_{B_k} \cL_{Z}+\lambda (\norm{Z}-\norm{Z^*}) \geq P\cL_{Z}-\frac{1}{20}(r_1^*)^2 - \lambda\norm{Z-Z^*} \geq -\frac{1}{20} (r^*_1)^2 -\lambda\rho^* \nonumber \\  
&= -\frac{1}{20} (r^*_1)^2 - \frac{11}{40}(r^*_2)^2\geq  - \frac{13}{40}(r^*_2)^2
\end{align}because $r_1^*\leq r_2^*$ (see Proposition~\ref{prop:prop_local_fixed_point_appendix} in the Appendix). Next, it follows from \eqref{eq:top_proof_lemma_6_12}, the fact that $\alpha>1$,  (\ref{eq:Lemma_nu_ER_01}) for $b=1$, (\ref{eq:Lemma_nu_ER_03}) for $b=1$ and (\ref{eq:Lemma_nu_ER_04})  that
\begin{align}\label{eq:Lemma_nu_ER_06}
\sup_{Z\in\cC}\mathrm{MOM}_K (\ell_{Z^*}-\ell_Z) + \lambda (\norm{Z^*}-\norm{Z}) & \leq  \max\left(\frac{-1}{5}, \frac{-29}{400}, \frac{13}{40}\right) r_2^2 = \frac{13}{40} (r_2^*)^2 
\end{align}
and so (\ref{eq:Lemma_eta_RMOM_ER-loc-20}) holds for $\eta \geq  13(r_2^*)^2/40$. As a consequence, \eqref{eq:Lemma_eta_RMOM_ER-loc-10} and (\ref{eq:Lemma_eta_RMOM_ER-loc-20}) both hold for $\eta =  132(r_2^*)^2/400$.
\end{proof}

At this stage, we have shown that on the event $\Omega_K$, the estimator $\hat{Z}$ has the statistical properties announced in Theorem \ref{theo:main_RERM_er}. In what follows we prove that in the framework of Theorem \ref{theo:main_RERM_er}, $\Omega_K$ holds with exponentially large probability.

\begin{Lemma}\label{Lemma:OmegaK_proba_RMOM_ER}
Assume that $K\geq 100|\cO|$, and let $\rho^*>0$ be such that it satisfies the sparsity equation from Definition \ref{def:sparsity_RMOM_er}. Then, $\Omega_K$ holds with probability at least $1-2\exp(-72K/625)$.
\end{Lemma}

\begin{proof}
Let $\phi:t\in \bR \rightarrow \mathbb{1}_{\left\{t\geq 1\right\}}+2(t-1/2)\mathbb{1}_{\left\{1/2\leq t\leq 1\right\}}$, so that for any $t\in \bR$, $\mathbb{1}_{\left\{t\geq 1\right\}}\leq \phi(t) \leq \mathbb{1}_{\left\{t\geq 1/2\right\}}$. For $k\in [K]$, let $W_k:=\left\{X_i:i\in B_k\right\}$ and $F_Z(W_k)=(P_{B_k}-P)\cL_Z$. We also define the counterparts of these quantities constructed with the non-corrupted vectors: $\widetilde{W}_k:=\left\{\widetilde{X}_i:i\in B_k\right\}$ and $F_Z(\widetilde{W}_k)=(\widetilde{P_{B_k}}-\widetilde{P})\cL_Z$, where $\widetilde{P_{B_k}}\cL_Z := \frac{K}{N}\sum_{i\in B_k}\cL_Z(\widetilde{X}_i)$ and $\widetilde{P}\cL_Z := \bE[\cL_{Z}(\widetilde{X}_i)]$. For both $b\in\{1, 2\}$,  let $Z\to \psi_b(Z)=\sum_{k\in [K]}\mathbb{1}_{\left\{|F_Z(W_k)|\leq (r^*_b)^2/20\right\}}$. Let $b\in\{1, 2\}$. We want to show that, with high probability, if $Z\in \cB_{b}$, then $\psi_b(Z) > K/2$ which follows  if one can proves that 
\begin{align}\label{eq:proof_main_RMOM_L2_loc_15}
\sum_{k\in [K]}\mathbb{1}_{\left\{|F_Z(\widetilde{W}_k)| > \frac{(r_b^*)^2}{20}\right\}} \leq \frac{49K}{100}.
\end{align}
Indeed, consider $Z\in \cB_{b}$ such that (\ref{eq:proof_main_RMOM_L2_loc_15}) holds. Then, there exist at least $(1-49/100)K = 51K/100$ blocks $B_k$ on which $|F_Z(\widetilde{W}_k)| \leq (r_b^*)^2/20$. On the other hand, we know that $|\cO|\leq K/100$, so that among the $51K/100$ previous blocks, at most $K/100$ contains corrupted data. The other $50K/100 = K/2$ contain only non-corrupted data, so we have $F_Z(\widetilde{W}_k) = F_Z(W_k)$ on these block and so $\psi_b(Z) > K/2$. \\

Let  $Z\in \cB_{b}$. We have:
\begin{align}\label{eq:proof_main_RMOM_L2_loc_16}
&\sum_{k\in \left[K\right]}\mathbb{1}_{\left\{|F_Z(\widetilde{W}_k)| > \frac{(r_b^*)^2}{20}\right\}} = \sum_{k\in \left[K\right]}\left[\mathbb{1}_{\left\{|F_Z(\widetilde{W}_k)| > \frac{(r_b^*)^2}{20}\right\}} - \bP\left(|F_Z(\widetilde{W}_k)| > \frac{(r_b^*)^2}{40}\right) + \bP\left(|F_Z(\widetilde{W}_k)| > \frac{(r_b^*)^2}{40}\right) \right] \nonumber\\
 & = \sum_{k\in \left[K\right]}\left(\mathbb{1}_{\left\{|F_Z(\widetilde{W}_k)| > \frac{(r_b^*)^2}{20}\right\}} - \bE\left[\mathbb{1}_{\left\{|F_Z(\widetilde{W}_k)| > \frac{(r_b^*)^2}{40}\right\}}\right] \right) + \sum_{k\in \left[K\right]} \bP\left(|F_Z(\widetilde{W}_k)| > \frac{(r_b^*)^2}{40}\right) \nonumber\\ 
 &\leq \sum_{k\in \left[K\right]}\left(\phi\left( \frac{20|F_Z(\widetilde{W}_k)|}{(r_b^*)^2}\right) -  \bE\left[\phi\left( \frac{20|F_Z(\widetilde{W}_k)|}{(r_b^*)^2}\right)\right]  \right) + \sum_{k\in \left[K\right]} \bP\left(|F_Z(\widetilde{W}_k)| > \frac{(r_b^*)^2}{40}\right) \nonumber\\ 
 &\leq \sup_{Z\in \cB_{b}} \left( \sum_{k\in \left[K\right]}\phi\left( \frac{20|F_Z(\widetilde{W}_k)|}{(r_b^*)^2}\right) -  \bE\left[\phi\left( \frac{20|F_Z(\widetilde{W}_k)|}{(r_b^*)^2}\right)\right]  \right) + \sum_{k\in \left[K\right]} \bP\left(|F_Z(\widetilde{W}_k)| > \frac{(r_b^*)^2}{40}\right)
\end{align}
We start with bounding the last sum in the previous inequality. For each $k\in\left[K\right]$, Markov's inequality and the definition of $r^*_b$ yield to
\begin{align*}
\bP\left(|F_Z(\widetilde{W}_k)| > \frac{(r_b^*)^2}{40}\right) & \leq \frac{1600}{(r_b^*)^4} \bE\left[F_Z(\widetilde{W}_k)^2\right] = \frac{1600}{(r_b^*)^4} \left(\frac{K}{N}\right)\mathrm{Var}(\cL_Z(\widetilde{X}))\leq   \frac{1600}{(r_b^*)^4} \left(V_K(r_b^*)\right)^2 \leq \frac{1}{200}
\end{align*}

Plugging this last result into (\ref{eq:proof_main_RMOM_L2_loc_16}), we get:
\begin{align}\label{eq:proof_main_RMOM_L2_loc_17}
\sum_{k\in \left[K\right]}\mathbb{1}_{\left\{|F_Z(\widetilde{W}_k)| > \frac{(r_b^*)^2}{20}\right\}} \leq \frac{K}{200}+ \sup_{Z\in \cB_{b}} \left( \sum_{k\in \left[K\right]}\phi\left( \frac{20|F_Z(\widetilde{W}_k)|}{(r_b^*)^2}\right) -  \bE\left[\phi\left( \frac{20|F_Z(\widetilde{W}_k)|}{(r_b^*)^2}\right)\right]  \right).
\end{align}

We now we have to bound this last term. Using Mc Diarmind inequality (Theorem 6.2 in \cite{BoucheronLugosiMassart} with $t=12/25$), we get that with probability at least $1-\exp(-72K/625)$, for all $Z\in \cB_{b}$, 
\begin{align}\label{eq:proof_main_RMOM_L2_loc_18}
\sum_{k\in \left[K\right]}\phi\left( \frac{20|F_Z(\widetilde{W}_k)|}{(r_b^*)^2}\right) -\bE\phi\left( \frac{20|F_Z(\widetilde{W}_k)|}{(r_b^*)^2}\right)\leq \frac{12K}{25} + \bE\left[\sup_{Z\in \cB_{b}}\sum_{k\in \left[K\right]}\phi\left( \frac{20|F_Z(\widetilde{W}_k)|}{(r_b^*)^2}\right) -\bE\phi\left( \frac{20|F_Z(\widetilde{W}_k)|}{(r_b^*)^2}\right)\right]. 
\end{align}
Let now $\epsilon_1, \ldots, \epsilon_{K}$ be Rademacher variables independant from the $\widetilde{X}_i$'s. By the symmetrization Lemma, we have:
\begin{align}\label{eq:proof_main_RMOM_L2_loc_19}
\bE\left[\sup_{Z\in \cB_{b}}\sum_{k\in \left[K\right]}\phi\left( \frac{20|F_Z(\widetilde{W}_k)|}{(r_b^*)^2}\right) -\bE\left[\phi\left( \frac{20|F_Z(\widetilde{W}_k)|}{(r_b^*)^2}\right)\right]\right]   \leq 2\bE\left[\sup_{Z\in \cB_{b}}\sum_{k\in \left[K\right]}\epsilon_k\phi\left( \frac{20|F_Z(\widetilde{W}_k)|}{(r_b^*)^2}\right)\right].
\end{align}
As $\phi$ is Lipschitz with $\phi(0)=0$, we can use the contraction Lemma (see \cite{LedTal01}, chapter 4) to get that:
\begin{align}\label{eq:proof_main_RMOM_L2_loc_20}
\bE\left[\sup_{Z\in \cB_{b}}\sum_{k\in \left[K\right]}\epsilon_k\phi\left( \frac{20|F_Z(\widetilde{W}_k)|}{(r_b^*)^2}\right)\right] &\leq 2\bE\left[\sup_{Z\in \cB_{b}}\sum_{k\in \left[K\right]}\epsilon_k\frac{20F_Z(\widetilde{W}_k)}{(r_b^*)^{2}}\right]= \frac{40}{(r_b^*)^{2}}\bE\left[\sup_{Z\in \cB_{b}}\sum_{k\in \left[K\right]}\epsilon_k(\widetilde{P_{B_k}}-\widetilde{P})\cL_Z\right]
\end{align}

Now, let $(\sigma_i)_{i=1,\ldots, N}$ be a family of Rademacher variables independant from the $\widetilde{X}_i$'s and the $\epsilon_i$'s. Using the symmetrization Lemma again, we get:
\begin{align*}
\bE\left[\sup_{Z\in \cB_{b}}\sum_{k\in \left[K\right]}\epsilon_k(\widetilde{P_{B_k}}-\widetilde{P})\cL_Z\right] & \leq 2 \bE\left[\sup_{Z\in \cB_{b}}\frac{K}{N}\sum_{i=1}^N\sigma_i\cL_Z(\widetilde{X}_i)\right] \leq  2KE(r_b^*, b\rho^*) \leq 2K\gamma (r_b^*)^2.
\end{align*}

Combining this with (\ref{eq:proof_main_RMOM_L2_loc_18}), (\ref{eq:proof_main_RMOM_L2_loc_19}) and (\ref{eq:proof_main_RMOM_L2_loc_20}), we finally get that, with probability at least $1-\exp(-72K/625)$:
\begin{align}\label{eq:proof_main_RMOM_L2_loc_21}
\sup_{Z\in \cC_{\gamma}}&\sum_{k\in \left[K\right]}\phi\left( \frac{20|F_Z(\widetilde{W}_k)|}{(r_b^*)^2}\right) -\bE\left[\phi\left( \frac{20|F_Z(\widetilde{W}_k)|}{(r_b^*)^2}\right)\right] \leq \left(\frac{12}{25}+160\gamma\right)K
\end{align}
Plugging that into (\ref{eq:proof_main_RMOM_L2_loc_17}), we conclude that, with probability at least $1-\exp(-72K/625)$, for all $Z\in\cB_b$,
\begin{align*}
\sum_{k\in \left[K\right]}\mathbb{1}_{\left\{|F_Z(\widetilde{W}_k)| > \frac{(r_b^*)^2}{20}\right\}} \leq \left(\frac{1}{200}+\frac{12}{25}+160\gamma\right) K \leq \frac{49}{100}K
\end{align*}
for our choice of parameters. Now, in order for $\Omega_K$ to hold, this inequality must be verified for both $b=1$ and $2$. Then, we finally conclude that $\Omega_{K}$ holds with probability $1-2\exp(-72K/625)$, which concludes the proof.

\end{proof}

\subsubsection{Proof of Theorem \ref{theo:main-RMOM-unl}}

Let $K>0$ be a divisor of $N$ such that $K\geq 100|\cO|$. Let $\gamma = 1/32000$. Let $A\in(0, 1]$ and $\rho^* >0$ be such that Assumption \ref{ass:curvature_RMOM_L2} holds and satisfying the sparsity equation from Definition~\ref{def:A-spar-MOM}. Define $\nu = 320000A^2$. 

For the sake of simplicity, we write all along this proof $\hat{Z}$ for $\hat{Z}^{\mathrm{RMOM}}_{K, \lambda}$.  For $b\in \{1, 2\}$, we define  $r_b^* = r^*_{\mathrm{RMOM}, L_2}(\gamma, b \rho^*)$,
\begin{align*}
C_{K, b} := \max\left(\nu\frac{K}{N}, (r^*_b)^2 \right) = C_K(\gamma, b\rho^*, A),
\end{align*}and   the localized models
 $\cB_{K, b}:=\left\{Z\in \cC:\norm{Z-Z^*}\leq b\rho^* \mbox{ and } \norm{Z-Z^*}_{L_2}\leq \sqrt{C_{K,b}}\right\}$ - we recall that the $L_2$-norm associated with the good data $\widetilde{X}$ is defined as $\norm{Z}_{L_2} = \bE[\inr{\widetilde{X}, Z}^2]^{1/2}$. With these notation, we have $\lambda:=(11/(40\rho^*)) C_{K, 2}$.  Finally, we define the event onto which $\hat Z$ will have the desired properties:
\begin{align*}
\Omega_K = \left\{\forall b \in \{1, 2\}, \forall Z \in \cB_{K, b}, ~~ \sum_{k=1}^K\mathbb{1}_{\left\{\left|(P_{B_k}-P)\cL_Z\right|\leq \frac{C_{K, b}}{20} \right\}} > \frac{K}{2}\right\}.
\end{align*}

First, we show that on $\Omega_K$, $\hat{Z}$ has the statistical properties announced in Theorem \ref{theo:main-RMOM-unl}. Then, we show that $\Omega_K$ holds with high probability.

\begin{Lemma}\label{Lemma:in_Bb}
If there exists $\eta >0$ such that
\begin{align}\label{eta1}
 \sup_{Z\in \cC \backslash \cB_{K, 2}} \mathrm{MOM}_K(\ell_{Z^*} - \ell_Z) + \lambda (\norm{Z^*} - \norm{Z}) < -\eta
\end{align}
and
\begin{align}\label{eta2}
 \sup_{Z\in \cC} \mathrm{MOM}_K(\ell_{Z^*} - \ell_Z) + \lambda (\norm{Z^*} - \norm{Z}) \leq  \eta 
\end{align}
then $\norm{Z-Z^*}\leq 2\rho^*$ and $\norm{Z-Z^*}_{L_2}\leq \sqrt{C_{K,2}}$.
\end{Lemma}

\begin{proof}
Assume that such an $\eta$ exists. For $Z\in \cC$, define $S(Z)=\sup_{Z^\prime\in \cC}\mathrm{MOM}_K(\ell_Z-\ell_{Z^\prime})+\lambda(\norm{Z}-\norm{Z^\prime})$. For $Z\in \cC\backslash \cB_{K, 2}$ we have:
\begin{align*}
S(Z) \geq \mathrm{MOM}_K(\ell_Z-\ell_{Z^*})+\lambda(\norm{Z}-\norm{Z^*}) \geq \inf_{Z\in \cC \backslash \cB_{K, 2}} \mathrm{MOM}_K(\ell_Z-\ell_{Z^*})+\lambda(\norm{Z}-\norm{Z^*}) > \eta
\end{align*}
since (\ref{eta1}) holds. Moreover, we have by definition of $\hat{Z}$:
\begin{align*}
S(\hat{Z}) \leq S(Z^*) = \sup_{Z\in \cC} \mathrm{MOM}_K(\ell_{Z^*} - \ell_Z) + \lambda (\norm{Z^*} - \norm{Z}) \leq  \eta 
\end{align*}
since (\ref{eta2}) holds. This shows that necessarily $\hat{Z}\in \cB_{K, 2}$.
\end{proof}

We are now looking for $\eta >0$ such that (\ref{eta1}) and (\ref{eta2}) hold. In the following result we identify such a $\eta$ on the event $\Omega_K$. 

\begin{Lemma}\label{lemma:proof_main_RMOM_L2_eta}
Under the conditions of Theorem \ref{theo:main-RMOM-unl} and on the event $\Omega_K$, (\ref{eta1}) and (\ref{eta2}) hold with $\eta = 33C_{K, 2}/100$.
\end{Lemma}

\begin{proof}
Consider $b\in \{1, 2\}$ and  $Z\in \cC \backslash \cB_{K, b}$. From the star-shaped property of $\cC$, we have the existence of $Z_0\in\partial\cB_{K,b}$ and $\alpha > 1$ such that $Z = Z^*+\alpha(Z_0 - Z^*)$. As a consequence, by linearity of the loss function and convexity of the regularization norm, for all $k\in[K]$ we have 
\begin{align}\label{eq:proof_main_RMOM_L2_eta_01}
&P_{B_k}\cL^{\lambda}_{Z} = P_{B_k}\cL_{Z} + \lambda (\norm{Z} - \norm{Z^*}) = \alpha P_{B_k}\cL_{Z_0} + \lambda (\norm{\alpha Z_0 + (1-\alpha)Z^*} - \norm{Z^*})\nonumber\\ 
&\geq \alpha P_{B_k}\cL_{Z_0} + \lambda \alpha(\norm{Z_0} - \norm{Z^*}) = \alpha P_{B_k}\cL^{\lambda}_{Z_0}.
\end{align}

Now, since $Z_0\in\partial\cB_{K,b}$, we have either \textit{a)} $\norm{Z_0-Z^*}_{L_2}=\sqrt{C_{K,b}}$ and $\norm{Z_0-Z^*} < b\rho^*$ or \textit{b)} $\norm{Z_0-Z^*}_{L_2} < \sqrt{C_{K,b}}$ and $\norm{Z_0-Z^*} = b\rho^*$.

In the first case \textit{a)}, on $\Omega_K$, there are at least $K/2$ blocks $B_k$ on which $P_{B_k}\cL_{Z_0}\geq P\cL_{Z_0} - C_{K,b}/(20)$. But from Assumption \ref{ass:curvature_RMOM_L2}, we have in this case that $AP\cL_{Z_0}\geq \norm{Z_0-Z^*}_{L_2}^2 = C_{K, b}$, so that, on the same blocks of data, $P_{B_k}\cL_{Z_0}\geq (1/A)C_{K,b} - (1/20)C_{K,b} \geq (19/20)C_{K,b}$, since we assumed that $0<A\leq 1$. Therefore,  on these blocs, we have
\begin{align*}
  &P_{B_k}\cL_{Z_0}^{\lambda}  = P_{B_k}\cL_{Z_0} + \lambda (\norm{Z_0}-\norm{Z^*}) \geq \frac{19}{20}C_{K,b} - \lambda \norm{Z_0-Z^*}\nonumber\\ 
  &\geq \frac{19}{20}C_{K,b} - \lambda b\rho^*  = \frac{19}{20}C_{K,b} - \frac{11b}{40} C_{K,2}. 
\end{align*}
But thanks to Proposition~\ref{prop:prop_local_fixed_point_appendix} from the Appendix, we have that $r_1^*\geq r_2^*/\sqrt{2}$, from which we deduce that $C_{K,1}\geq C_{K,2}/2$. As a consequence, on the previous blocks, we have
\begin{align}\label{eq:proof_main_RMOM_L2_eta_02}
P_{B_k}\cL_{Z_0}^{\lambda} \geq \left\{
\begin{array}{cc}
7C_{K, 2}/40 & \mbox{ for } b=1\\
16C_{K, 2}/40  & \mbox{ for } b=2.
\end{array}
\right.
\end{align}

In the second case \textit{b)}, we have $Z_0 \in \widetilde{H}_{b\rho^*, A}$ from Definition \ref{def:A-spar-MOM}. Since the sparsity equation is satisfied by $\rho^*$, it is also satisfied by $b\rho^*$ as well (see Proposition~\ref{prop:prop_sparsity_equation} in the Appendix). Let $V\in H$ be such that $\norm{Z^*-V}\leq b\rho^*/20$ and $\Phi \in \partial \norm{.}(V)$. We have:
\begin{align*}
\norm{Z_0} - \norm{Z^*} & \geq \norm{Z_0}- \norm{V} - \norm{Z^*-V}  \\
& \geq \inr{\Phi, Z_0-V} - \norm{Z^*-V} ~~ (\mbox{ since } \Phi \in \partial \norm{.}(V)) \\
& = \inr{\Phi, Z_0-Z^*} - \inr{\Phi, V-Z^*} - \norm{Z^*-V} \\
& \geq \inr{\Phi, Z_0-Z^*} - 2 \norm{Z^*-V} ~~ (\mbox{ since } \inr{\Phi,U}\leq \norm{U} \mbox{ for any } U \in H) \\
& \geq \inr{\Phi, Z_0-Z^*} - \frac{b\rho^*}{10}.
\end{align*}
This is true for any $\Phi\in \underset{V \in Z^*+b\rho^*/20}{\cup}\partial \norm{.}(V) = \Gamma_{Z^*}(b\rho^*)$. Then taking the $\sup$ over $\Gamma_{Z^*}(b\rho^*)$ gives:
\begin{align*}
\norm{Z_0} - \norm{Z^*}  \geq \sup_{\Phi \in \Gamma_{Z^*}(b\rho^*)} \inr{\Phi, Z_0-Z^*} - \frac{b\rho^*}{10}
\end{align*}
and then taking the infimum over $\widetilde{H}_{b\rho^*, A}$ gives:
\begin{align}\label{eq:proof_main_RMOM_L2_eta_03}
\norm{Z_0} - \norm{Z^*}  & \geq \inf_{Z_0 \in\widetilde{H}_{b\rho^*, A}} \norm{Z_0} - \norm{Z^*} \geq \inf_{Z_0 \in \widetilde{H}_{b\rho^*, A}} \sup_{\Phi \in \Gamma_{Z^*}(2\rho^*)} \inr{\Phi, Z_0-Z^*} - \frac{b\rho^*}{10}  = \Delta(b\rho^*) - \frac{b\rho^*}{10} \geq \frac{7}{10} b\rho^*
\end{align}
where the last inequality holds since $b\rho^*$ satisfies the sparsity equation. Then, $\lambda(\norm{Z_0}-\norm{Z^*})\geq (7/10)\lambda b\rho^* = (77/400)bC_{K,2}$. Now, since $Z_0\in \cB_{K,b}$, on $\Omega_K$ there exist at least $K/2$ blocks $B_k$ such that $|(P_{B_k}-P)\cL_{Z_0}|\leq C_{K,b}/(20) $ and so $P_{B_k} \cL_{Z_0}\geq -C_{K,b}/(20)$ (because $P\cL_{Z_0}\geq 0$). Therefore, on the very same  blocks,
\begin{align}\label{eq:proof_main_RMOM_L2_eta_04}
 P_{B_k} \cL^{\lambda}_{Z_0} =  P_{B_k} \cL_{Z_0}+\lambda (\norm{Z_0}-\norm{Z^*}) \geq  -\frac{1}{20}  C_{K,b} + \frac{77b}{400}C_{K,2}\geq \left\{
\begin{array}{cc}
57(r_2^*)^2/400 & \mbox{ for } b=1\\
134(r_2^*)^2/400 & \mbox{ for } b=2
\end{array}
\right.
\end{align}where we used that $C_{K, 1}\leq C_{K, 2}$ because $r^*_1\leq r^*_2$ (see Proposition~\ref{prop:prop_local_fixed_point_appendix} in the Appendix). As a consequence, it follows from \eqref{eq:proof_main_RMOM_L2_eta_01}, the fact that $\alpha>1$,  \eqref{eq:proof_main_RMOM_L2_eta_02} and \eqref{eq:proof_main_RMOM_L2_eta_04} for $b=2$ that, for all $Z\in\cC\backslash\cB_{K,2}$, on more than $K/2$ blocks $B_k$: $P_{B_k}\cL_Z^\lambda\geq (134/400)C_{K, 2}$ and so \eqref{eta1} holds for $\eta <  (134/400)C_{K, 2}$. \\ 

Let us now turn to Equation~\eqref{eta2}. Let $Z\in\cB_{K,1}$. On $\Omega_K$ there exist at least $K/2$ blocks $B_k$ such that $|(P_{B_k}-P)\cL_{Z}|\leq C_{K, 1}/20$. On these blocks $B_k$, all $P_{B_k} \cL^{\lambda}_Z$'s are such that
\begin{align}\label{eq:proof_main_RMOM_L2_eta_05}
 &P_{B_k} \cL^{\lambda}_Z = P_{B_k} \cL_{Z}+\lambda (\norm{Z}-\norm{Z^*}) \geq P\cL_{Z}-\frac{1}{20}C_{K,1} - \lambda\norm{Z-Z^*} \geq -\frac{1}{20} C_{K,1} -\lambda\rho^* \nonumber \\  
&= -\frac{1}{20} C_{K,1} - \frac{11}{40}C_{K,2}\geq  - \frac{13}{40}C_{K,2}
\end{align}where we used the fact that, thanks to Proposition~\ref{prop:prop_local_fixed_point_appendix} in the Appendix, $C_{K,1}\leq C_{K,2}$. Next, it follows from \eqref{eq:proof_main_RMOM_L2_eta_01}, the fact that $\alpha>1$,  (\ref{eq:proof_main_RMOM_L2_eta_02}) and (\ref{eq:proof_main_RMOM_L2_eta_01}) for $b=1$ and (\ref{eq:proof_main_RMOM_L2_eta_05})  that
\begin{align}\label{eq:proof_main_RMOM_L2_eta_06}
\sup_{Z\in\cC}\mathrm{MOM}_K (\ell_{Z^*}-\ell_Z) + \lambda (\norm{Z^*}-\norm{Z}) & \leq  \max\left(\frac{-7}{40}, \frac{-57}{400}, \frac{13}{40}\right) C_{K,2} = \frac{13}{40} C_{K,2}
\end{align}
and so (\ref{eta2}) holds for $\eta \geq  13C_{K,2}/40$. As a consequence, \eqref{eta1} and (\ref{eta2}) both hold for $\eta =  132C_{K,2}/400$.

\end{proof}

From Lemmas \ref{Lemma:in_Bb} and \ref{lemma:proof_main_RMOM_L2_eta}, we conclude that on the event $\Omega_K$, $\hat Z \in \cB_{K, 2}$. We use this information to upper bound the excess risk of $\hat Z$ in the following result.

\begin{Lemma}\label{lemma:proof_main_RMOM_L2_excess_risk}
Under the conditions of Theorem \ref{theo:main-RMOM-unl} and on the event $\Omega_K$, we have $P\cL_{\hat Z} \leq (27/100)C_{K,2}$.
\end{Lemma}

\begin{proof}
From Lemmas \ref{Lemma:in_Bb} and \ref{lemma:proof_main_RMOM_L2_eta}, we have that $\hat Z \in \cB_{K, 2}$. On $\Omega_K$, this implies the existence of stricly more than $K/2$ blocks $B_k$ on which
\begin{align}\label{eq:proof_main_RMOM_L2_eta_061}
P\cL_{\hat{Z}}\leq P_{B_k}\cL_{\hat{Z}}+ C_{K, 2}/20.
\end{align}
Now, by definition of $\hat{Z}$, ($\ref{eta2}$) and Lemma~\ref{lemma:proof_main_RMOM_L2_eta} we get
\begin{align*}
\mathrm{MOM}_K(\ell_{\hat{Z}}- \ell_{Z^*}) + \lambda(\norm{\hat{Z}}-\norm{Z^*}) \leq \sup_{Z\in \cC} \mathrm{MOM}_K(\ell_{Z^*}- \ell_{Z}) + \lambda(\norm{Z^*}-\norm{Z}) \leq 33C_{K, 2}/100.
\end{align*} 
This means that there exist at least $K/2$ blocks $B_k$ on which $P_{B_k}\cL_{\hat Z} + \lambda(\norm{\hat{Z}}-\norm{Z^*}) \leq 33C_{K, 2}/100$. Since $\lambda(\norm{Z^*}-\norm{\hat{Z}})\leq \lambda\norm{Z^*-\hat{Z}}\leq 2\lambda\rho^* = 11C_{K, 2}/20$, we have on these blocks
\begin{align}\label{eq:proof_main_RMOM_L2_eta_062}
P_{B_k}\cL_{\hat{Z}}\leq 33C_{K, 2}/100 + 11C_{K, 2}/20 = 22 C_{K, 2}/100.
\end{align}
Therefore, there exist at least a block $B_{k_0}$ on which (\ref{eq:proof_main_RMOM_L2_eta_061}) and (\ref{eq:proof_main_RMOM_L2_eta_062}) hold simultaneously. On this block, we can write
\begin{align*}
P\cL_{\hat{Z}}\leq P_{B_{k_0}}\cL_{\hat{Z}}+ C_{K, 2}/20 \leq 22 C_{K, 2}/100 + C_{K, 2}/20 = 27C_{K, 2}/100.
\end{align*}
\end{proof}

At this stage, we have shown that on the event $\Omega_K$, the regularized minmax MOM-estimator $\hat{Z}$ has the statistical properties announced in Theorem \ref{theo:main-RMOM-unl}. In what follows, we prove that, in the framework of Theorem \ref{theo:main-RMOM-unl}, $\Omega_K$ holds with exponentially large probability.

\begin{Proposition}\label{prop:OmegaK}
Consider $\rho^*$ that satisfies the sparsity equation from Definition \ref{def:A-spar-MOM}. Assume that $K\geq 100|\cO|$. Then, $\Omega_K$ holds with probability at least $1-2\exp(-72K/625)$.
\end{Proposition}

\begin{proof}
Let $b\in \{1, 2\}$. Let $\phi:t\in \bR \rightarrow \mathbb{1}_{\left\{t\geq 1\right\}}+2(t-1/2)\mathbb{1}_{\left\{1/2\leq t\leq 1\right\}}$, so that for any $t\in \bR$, $\mathbb{1}_{\left\{t\geq 1\right\}}\leq \phi(t) \leq \mathbb{1}_{\left\{t\geq 1/2\right\}}$. For $k\in \left[K\right]$, let $W_k:=\left\{X_i:i\in B_k\right\}$ and $F_Z(W_k)=(P_{B_k}-P)\cL_Z$. We also define the counterparts of these quantities constructed with the non-corrupted vectors: $\widetilde{W}_k:=\left\{\widetilde{X}_i:i\in B_k\right\}$ and $F_Z(\widetilde{W}_k)=(\widetilde{P_{B_k}}-P)\cL_Z$, where $\widetilde{P_{B_k}}\cL_Z := (K/N)\sum_{i\in B_k}\cL_Z(\widetilde{X}_i)$. Let $\psi_b(Z)=\sum_{k\in \left[K\right]}\mathbb{1}_{\left\{|F_Z(W_k)|\leq  C_{K,b}/(20)\right\}}$. We would like to show that, if $Z\in \cB_{K, b}$, then $\psi_b(Z) > K/2$ with high probability. As we showed in the proof of Lemma \ref{Lemma:main_MOM_loc_er_OmegaK_proba}, in our framework this is true if we show that with high probability, for all $Z\in\cB_{K,b}$,
\begin{align}\label{eq:psi_non_corrupted_RMOM_L2}
\sum_{k\in \left[K\right]}\mathbb{1}_{\left\{|F_Z(\widetilde{W}_k)| > \frac{C_{K,b}}{20}\right\}} \leq \frac{49K}{100}
\end{align}
and this is what we do now. Let $Z\in \cB_{K, b}$. We have:
\begin{align}\label{eq:proof_main_RMOM_L2_eta_07}
&\sum_{k\in \left[K\right]}\mathbb{1}_{\left\{|F_Z(\widetilde{W}_k)| > \frac{C_{K,b}}{20}\right\}} = \sum_{k\in \left[K\right]}\left[\mathbb{1}_{\left\{|F_Z(\widetilde{W}_k)| > \frac{C_{K,b}}{20}\right\}} - \bP\left(|F_Z(\widetilde{W}_k)| > \frac{C_{K,b}}{40}\right) + \bP\left(|F_Z(\widetilde{W}_k)| > \frac{C_{K,b}}{40}\right) \right] \nonumber\\
 & = \sum_{k\in \left[K\right]}\left(\mathbb{1}_{\left\{|F_Z(\widetilde{W}_k)| > \frac{C_{K,b}}{20}\right\}} - \bE\left[\mathbb{1}_{\left\{|F_Z(\widetilde{W}_k)| > \frac{C_{K,b}}{40}\right\}}\right] \right) + \sum_{k\in \left[K\right]} \bP\left(|F_Z(\widetilde{W}_k)| > \frac{C_{K,b}}{40}\right) \nonumber\\ 
 &\leq \sum_{k\in \left[K\right]}\left(\phi\left( \frac{20|F_Z(\widetilde{W}_k)|}{C_{K,b}}\right) -  \bE\left[\phi\left( \frac{20|F_Z(\widetilde{W}_k)|}{C_{K,b}}\right)\right]  \right) + \sum_{k\in \left[K\right]} \bP\left(|F_Z(\widetilde{W}_k)| > \frac{C_{K,b}}{40}\right) \nonumber\\ 
 &\leq \sup_{Z\in \cB_{K, b}} \left( \sum_{k\in \left[K\right]}\phi\left( \frac{20|F_Z(\widetilde{W}_k)|}{C_{K,b}}\right) -  \bE\left[\phi\left( \frac{20|F_Z(\widetilde{W}_k)|}{C_{K,b}}\right)\right]  \right) + \sum_{k\in \left[K\right]} \bP\left(|F_Z(\widetilde{W}_k)| > \frac{C_{K,b}}{40}\right)
\end{align}
We start with bounding the last sum in the previous inequality. For each $k\in\left[K\right]$, Markov's inequality and the definition of $C_{K,b}$ yield to
\begin{align*}
\bP\left(|F_Z(\widetilde{W}_k)| > \frac{ C_{K,b}}{40}\right) & \leq \left(\frac{40}{C_{K,b}}\right)^2 \bE\left[|F_Z(\widetilde{W}_k)|^2\right] = \left(\frac{40}{C_{K,b}}\right)^2   \left(\frac{K}{N}\right)^2 \bE\left[\left(\sum_{i\in B_k}\cL_Z(\widetilde{X}_i) - \bE\left[\cL_Z(\widetilde{X}_i)\right]\right)^2\right] \\ 
& \leq \left(\frac{40}{C_{K,b}}\right)^2  \frac{K}{N} \norm{Z-Z^*}_{L_2}^2 \leq \left(\frac{40}{C_{K,b}}\right)^2   \frac{K}{N} C_{K,b} \leq 40^2 \nu^{-1} = \frac{1}{200}.
\end{align*}

Plugging this last result into (\ref{eq:proof_main_RMOM_L2_eta_07}), we get:
\begin{align}\label{eq:proof_main_RMOM_L2_eta_08}
&\sum_{k\in \left[K\right]}\mathbb{1}_{\left\{|F_Z(\widetilde{W}_k)| > \frac{C_{K,b}}{20}\right\}} \leq \frac{K}{200}+ \sup_{Z\in \cB_{K, b}} \left( \sum_{k\in \left[K\right]}\phi\left( \frac{20|F_Z(\widetilde{W}_k)|}{C_{K,b}}\right) -  \bE\left[\phi\left( \frac{20|F_Z(\widetilde{W}_k)|}{C_{K,b}}\right)\right]  \right).
\end{align}

We now have to bound this last term. Using Mc Diarmind inequality (Theorem 6.2 in \cite{BoucheronLugosiMassart} with $t=12/25$), we get that with probability at least $1-\exp(-72K/625)$, for all $Z\in \cB_{K, b}$, 
\begin{align}\label{eq:proof_main_RMOM_L2_eta_09}
\sum_{k\in \left[K\right]}\phi\left( \frac{20|F_Z(\widetilde{W}_k)|}{C_{K,b}}\right) -\bE\phi\left( \frac{20|F_Z(\widetilde{W}_k)|}{C_{K,b}}\right)\leq \frac{12K}{25} + \bE\left[\sup_{Z\in \cB_{K, b}}\sum_{k\in \left[K\right]}\phi\left( \frac{20|F_Z(\widetilde{W}_k)|}{C_{K,b}}\right) -\bE\phi\left( \frac{20|F_Z(\widetilde{W}_k)|}{C_{K,b}}\right)\right]. 
\end{align}
Let now $\epsilon_1, \ldots, \epsilon_{K}$ be Rademacher variables independant from the $\widetilde{X}_i$'s. By the symmetrization Lemma, we have:
\begin{align}\label{eq:proof_main_RMOM_L2_eta_10}
\bE\left[\sup_{Z\in \cB_{K, b}}\sum_{k\in \left[K\right]}\phi\left( \frac{20|F_Z(\widetilde{W}_k)|}{C_{K,b}}\right) -\bE\left[\phi\left( \frac{20|F_Z(\widetilde{W}_k)|}{C_{K,b}}\right)\right]\right]   \leq 2\bE\left[\sup_{Z\in \cB_{K, b}}\sum_{k\in \left[K\right]}\epsilon_k\phi\left( \frac{20|F_Z(\widetilde{W}_k)|}{C_{K,b}}\right)\right].
\end{align}
As $\phi$ is Lipschitz with $\phi(0)=0$, we can use the contraction Lemma (see \cite{LedTal01}, chapter 4) to get that:
\begin{align}\label{eeq:proof_main_RMOM_L2_eta_11}
\bE\left[\sup_{Z\in \cB_{K, b}}\sum_{k\in \left[K\right]}\epsilon_k\phi\left( \frac{20|F_Z(\widetilde{W}_k)|}{C_{K,b}}\right)\right] &\leq 2\bE\left[\sup_{Z\in \cB_{K, b}}\sum_{k\in \left[K\right]}\epsilon_k\frac{20F_Z(\widetilde{W}_k)}{C_{K, b}}\right]= \frac{40}{C_{K, b}}\bE\left[\sup_{Z\in \cB_{K, b}}\sum_{k\in \left[K\right]}\epsilon_k(\widetilde{P_{B_k}}-\widetilde{P})\cL_Z\right]
\end{align}

Now, let $(\sigma_i)_{i=1,\ldots, N}$ be a family of Rademacher variables independant from the $\widetilde{X}_i$'s and the $\epsilon_i$'s. Using the symmetrization Lemma again, we get:
\begin{align*}
\bE\left[\sup_{Z\in \cB_{K, b}}\sum_{k\in \left[K\right]}\epsilon_k(\widetilde{P_{B_k}}-\widetilde{P})\cL_Z\right] & \leq 2 \bE\left[\sup_{Z\in \cB_{K, b}}\frac{K}{N}\sum_{i=1}^N\sigma_i\cL_Z(\widetilde{X}_i)\right] \leq  2K(r_b^*)^2 \leq 2K\gamma C_{K,b}.
\end{align*}

Combining this with (\ref{eq:proof_main_RMOM_L2_eta_09}), (\ref{eq:proof_main_RMOM_L2_eta_10}) and (\ref{eeq:proof_main_RMOM_L2_eta_11}), we finally get that, with probability at least $1-\exp(-72K/625)$:
\begin{align}\label{eq:proof_main_RMOM_L2_eta_12}
\sup_{Z\in \cB_{K, b}}&\sum_{k\in \left[K\right]}\phi\left( \frac{20|F_Z(\widetilde{W}_k)|}{C_{K,b}}\right) -\bE\left[\phi\left( \frac{20|F_Z(\widetilde{W}_k)|}{C_{K,b}}\right)\right] \leq \left(\frac{12}{25}+160\gamma \right)K
\end{align}
Plugging that into \eqref{eq:proof_main_RMOM_L2_eta_08}, we conclude that, with probability at least $1-\exp(-72K/625)$, for all $Z\in\cB_{K,b}$,
\begin{align*}
\sum_{k\in \left[K\right]}\mathbb{1}_{\left\{|F_Z(\widetilde{W}_k)| > \frac{C_{K,b}}{20}\right\}} \leq \left(\frac{1}{200}+\frac{12}{25}+160\gamma \right) K \leq \frac{49}{100}K
\end{align*}
for our choice of parameters. Now, in order for $\Omega_K$ to hold, this inequality must be verified for both $b=1$ and $2$. Then, we finally conclude that $\Omega_{K}$ holds with probability $1-2\exp(-72K/625)$, which concludes the proof.

\end{proof}

\subsubsection{Proof of Theorem \ref{theo:main_RMOM_G-loc}}

The proof is structured in the same way as the previous ones: we identify an event on which $\hat{Z}^{\mathrm{RMOM}}_{K, \lambda}$ has the desired statistical properties, then we show that this event holds with high probability. We place ourselves under the conditions of Theorem \ref{theo:main_RMOM_G-loc}, i.e., we assume the existence of $A\in (0,1]$ such that Assumption \ref{ass:curvature_RMOM_G-loc} holds, $\gamma = 1/32000$ and $\rho^*$  which satisfies the sparsity equation from Definition~\ref{def:sparsity_RMOM_G-loc}.  For $b\in \{1, 2\}$ we define $r^*_b = r^*_{\mathrm{RMOM, G}}(\gamma, 2\rho^*)$ and $\cB_b:=\left\{Z\in \cC:G(Z-Z^*)\leq (r^*_b)^2\mbox{ and }\norm{Z-Z^*}\leq b\rho^* \right\}$. With these notation, $\lambda=(11/(40\rho^*))r^*_2$. We consider the event
\begin{align*}
\Omega_{K,G} = \left\{\forall b \in \{1, 2\}, \forall Z \in \cB_b, ~~ \sum_{k=1}^K\mathbb{1}\left(\left|(P_{B_k}-P)\cL_Z\right|\leq \frac{1}{20}(r^*_b)^2\right) > \frac{K}{2}\right\}
\end{align*}
For the sake of simplicity, in the rest of the proof we write $\hat{Z} = \hat{Z}^{\mathrm{RMOM}}_{K, \lambda}$.

\begin{Lemma}\label{Lemma:proof_main_RMOM_G-loc-00}
If there exists $\eta > 0$ such that
\begin{align}\label{eq:Lemma_eta_RMOM_G-loc_01}
&\underset{Z\in\cC\backslash\cB_{2}}{\mathrm{sup}} \mathrm{MOM}_K(\ell_{Z^*}-\ell_{Z}) + \lambda (\norm{Z^*} - \norm{Z}) < -\eta 
\end{align}
and
\begin{align}\label{eq:Lemma_eta_RMOM_G-loc_02}
 \underset{Z\in\cC}{\mathrm{sup}} ~ \mathrm{MOM}_K(\ell_{Z^*}-\ell_{Z}) + \lambda (\norm{Z^*} - \norm{Z}) \leq \eta 
\end{align}
then $\norm{Z-Z^*}\leq 2\rho^*$ and $G(Z-Z^*)\leq r^*_{\mathrm{RMOM, G}}(\gamma, 2\rho^*)^2$.
\end{Lemma}

\begin{proof}
Let $\eta$ be such that \eqref{eq:Lemma_eta_RMOM_G-loc_01} and \eqref{eq:Lemma_eta_RMOM_G-loc_02} hold. For all $Z\in \cC$, define $S(Z)=\sup_{Z^\prime\in \cC}\mathrm{MOM}_K(\ell_Z-\ell_{Z^\prime})+\lambda(\norm{Z}-\norm{Z^\prime})$. It follows from \eqref{eq:Lemma_eta_RMOM_G-loc_01} that for all $Z\in \cC\backslash \cB_{2}$, 
\begin{align*}
S(Z) \geq \mathrm{MOM}_K(\ell_Z-\ell_{Z^*})+\lambda(\norm{Z}-\norm{Z^*})  > \eta
\end{align*}
Moreover, it follows from the definition of $\hat{Z}$ and (\ref{eq:Lemma_eta_RMOM_G-loc_02}) that
\begin{align*}
S(\hat{Z}) \leq S(Z^*) = \sup_{Z\in \cC} \mathrm{MOM}_K(\ell_{Z^*} - \ell_Z) + \lambda (\norm{Z^*} - \norm{Z}) \leq  \eta 
\end{align*}
This shows that necessarily $\hat{Z}\in \cB_{2}$.
\end{proof}

\begin{Lemma}\label{Lemma:proof_main_RMOM_G-loc-01}
Under the conditions of Theorem \ref{theo:main_RMOM_G-loc} and on the event $\Omega_{K, G}$, (\ref{eq:Lemma_eta_RMOM_G-loc_01}) and (\ref{eq:Lemma_eta_RMOM_G-loc_02}) hold with $\eta = (33/100)(r^*_2)^2$.
\end{Lemma}

\begin{proof} Let $b\in\{1,2\}$.
Let $Z\in\cC\backslash\cB_b$. By the star-shaped property of $\cC$ and the regularity property of $G$, there exist $Z_0\in \partial\cB_b$ and $\alpha>1$ such that $Z=Z^*+\alpha(Z_0-Z^*)$.  As a consequence, by linearity of the loss function and convexity of the regularization norm, for all $k\in[K]$ we have 
\begin{align}\label{eq:proof_lemma_main_Gloc_eta_01}
&P_{B_k}\cL^{\lambda}_{Z} = P_{B_k}\cL_{Z} + \lambda (\norm{Z} - \norm{Z^*}) = \alpha P_{B_k}\cL_{Z_0} + \lambda (\norm{\alpha Z_0 + (1-\alpha)Z^*} - \norm{Z^*})\nonumber\\ 
&\geq \alpha P_{B_k}\cL_{Z_0} + \lambda \alpha(\norm{Z_0} - \norm{Z^*}) = \alpha P_{B_k}\cL^{\lambda}_{Z_0}.
\end{align}

Now, since $Z_0\in\partial\cB_b$, we have either \textit{a)} $G(Z_0-Z^*) = (r^*_b)^2$ and $\norm{Z_0-Z^*} < b\rho^*$ or \textit{b)} $G(Z_0-Z^*) < (r^*_b)^2$ and $\norm{Z_0-Z^*} = b\rho^*$.

In the first case \textit{a)}, on $\Omega_{K, G}$, there are at least $K/2$ blocks $B_k$ on which $P_{B_k}\cL_{Z_0}\geq P\cL_{Z_0} - (r^*_b)^2/(20)$. But we also have from Assumption \ref{ass:curvature_RMOM_G-loc} that $AP\cL_{Z_0}\geq G(Z_0-Z^*) = (r^*_b)^2$, so that $P_{B_k}\cL_{Z_0}\geq (1/A) (r^*_b)^2 - (1/20)(r^*_b)^2 \geq(19/20)(r^*_b)^2$, since we assumed that $A\leq 1$. Therefore, on these blocs, we have
\begin{align}\label{eq:proof_lemma_main_Gloc_eta_02}
  &P_{B_k}\cL_{Z_0}^{\lambda}  = P_{B_k}\cL_{Z_0} + \lambda (\norm{Z_0}-\norm{Z^*}) \geq \frac{19}{20}(r^*_b)^2 - \lambda \norm{Z_0-Z^*}\nonumber\\ 
  &\geq \frac{19}{20}(r^*_b)^2 - \lambda b\rho^*  = \frac{19}{20}(r^*_b)^2 - \frac{11b}{40} (r^*_2)^2 \geq 
\left\{
\begin{array}{cc}
(r_2^*)^2/5 & \mbox{ for } b=1\\
2(r_2^*)^2/5 & \mbox{ for } b=2.
\end{array}
\right.
\end{align}
where we used in the case $b=1$ that $r_1^*\geq r_2^*/\sqrt{2}$ thanks to Proposition~\ref{prop:prop_local_fixed_point_appendix} from the Appendix. 

In the second case \textit{b)}, we have $Z_0 \in \bar{H}_{b\rho^*}$ from Definition \ref{def:sparsity_RMOM_G-loc}. Since the sparsity equation holds for $\rho=\rho^*$, it also holds for $\rho=b\rho^*$ (see Proposition~\ref{prop:prop_sparsity_equation} in the Appendix). Let $V\in H$ be such that $\norm{Z^*-V}\leq b\rho^*/20$ and $\Phi \in \partial \norm{.}(V)$. We have:
\begin{align*}
\norm{Z_0} - \norm{Z^*} & \geq \norm{Z_0}- \norm{V} - \norm{Z^*-V}  \\
& \geq \inr{\Phi, Z_0-V} - \norm{Z^*-V} ~~ (\mbox{ since } \Phi \in \partial \norm{.}(V)) \\
& = \inr{\Phi, Z_0-Z^*} - \inr{\Phi, V-Z^*} - \norm{Z^*-V} \\
& \geq \inr{\Phi, Z_0-Z^*} - 2 \norm{Z^*-V} ~~ (\mbox{ since } \inr{\Phi,U}\leq \norm{U} \mbox{ for any } U \in H) \\
& \geq \inr{\Phi, Z_0-Z^*} - \frac{b\rho^*}{10}.
\end{align*}
This is true for any $\Phi\in \underset{V \in Z^*+b\rho^*/20}{\cup}\partial \norm{.}(V) = \Gamma_{Z^*}(b\rho^*)$. Then taking the $\sup$ over $\Gamma_{Z^*}(b\rho^*)$ gives:
\begin{align*}
\norm{Z_0} - \norm{Z^*}  \geq \sup_{\Phi \in \Gamma_{Z^*}(b\rho^*)} \inr{\Phi, Z_0-Z^*} - \frac{b\rho^*}{10}
\end{align*}
and then taking the infimum over $\bar{H}_{b\rho^*, A}$ gives:
\begin{align}\label{eq:proof_lemma_main_Gloc_eta_03}
\norm{Z_0} - \norm{Z^*}  & \geq \inf_{Z_0 \in \bar{H}_{b\rho^*, A}} \norm{Z_0} - \norm{Z^*} \geq \inf_{Z_0 \in \bar{H}_{b\rho^*, A}} \sup_{\Phi \in \Gamma_{Z^*}(2\rho^*)} \inr{\Phi, Z_0-Z^*} - \frac{b\rho^*}{10}  = \Delta(b\rho^*) - \frac{b\rho^*}{10} \geq \frac{7}{10} b\rho^*
\end{align}
where the last inequality holds since $b\rho^*$ satisfies the sparsity equation. Then, $\lambda(\norm{Z_0}-\norm{Z^*})\geq (7/10)\lambda b\rho^* = (77/400)b(r^*_2)^2$. Now, since $Z_0\in \cB_b$, on $\Omega_{K, G}$ there exist at least $K/2$ blocks $B_k$ such that $|(P_{B_k}-P)\cL_{Z_0}|\leq (r^*_b)^2/(20) $ and so $P_{B_k} \cL_{Z_0}\geq - (r^*_b)^2/(20)$ (because $P\cL_{Z_0}\geq 0$). Therefore, on the very same  blocks,
\begin{align}\label{eq:proof_lemma_main_Gloc_eta_04}
 P_{B_k} \cL^{\lambda}_{Z_0} =  P_{B_k} \cL_{Z_0}+\lambda (\norm{Z_0}-\norm{Z^*}) \geq  -\frac{1}{20}  (r^*_b)^2 + \frac{77}{400}b(r^*_2)^2\geq \left\{
\begin{array}{cc}
57(r_2^*)^2/(400) & \mbox{ for } b=1\\
134(r_2^*)^2/(400) & \mbox{ for } b=2
\end{array}
\right.
\end{align}where we used that $r_1^*\leq r_2^*$ (see Proposition~\ref{prop:prop_local_fixed_point_appendix} in the Appendix). As a consequence, it follows from \eqref{eq:proof_lemma_main_Gloc_eta_01}, the fact that $\alpha>1$,  \eqref{eq:proof_lemma_main_Gloc_eta_02} and \eqref{eq:proof_lemma_main_Gloc_eta_04} for $b=2$ that, for all $Z\in\cC\backslash\cB_2$, on more than $K/2$ blocks $B_k$: $P_{B_k}\cL_Z^\lambda\geq (134/400)(r^*_2)^2$ and so \eqref{eq:Lemma_eta_RMOM_G-loc_01} holds for $\eta <  (134/400)(r^*_2)^2$.

Let us now turn to Equation~\eqref{eq:Lemma_eta_RMOM_G-loc_02}. Let $Z\in\cB_1$. On $\Omega_{K, G}$ there exist at least $K/2$ blocks $B_k$ such that $|(P_{B_k}-P)\cL_{Z}|\leq (r^*_1)^2/(20)$. On these blocks $B_k$, all $P_{B_k} \cL^{\lambda}_Z$'s are such that
\begin{align}\label{eq:proof_lemma_main_Gloc_eta_05}
 &P_{B_k} \cL^{\lambda}_Z = P_{B_k} \cL_{Z}+\lambda (\norm{Z}-\norm{Z^*}) \geq P\cL_{Z}-\frac{1}{20}(r_1^*)^2 - \lambda\norm{Z-Z^*} \geq -\frac{1}{20} (r^*_1)^2 -\lambda\rho^* \nonumber \\  
&= -\frac{1}{20} (r^*_1)^2 - \frac{11}{40}(r^*_2)^2\geq  - \frac{13}{40}(r^*_2)^2
\end{align}because $r_1^*\leq r_2^*$ (see Proposition~\ref{prop:prop_local_fixed_point_appendix} in the Appendix). Next, it follows from \eqref{eq:proof_lemma_main_Gloc_eta_01}, the fact that $\alpha>1$,  (\ref{eq:proof_lemma_main_Gloc_eta_02}) and (\ref{eq:proof_lemma_main_Gloc_eta_04}) for $b=1$ and (\ref{eq:proof_lemma_main_Gloc_eta_05}) that
\begin{align}\label{eq:proof_lemma_main_Gloc_eta_06}
\sup_{Z\in\cC}\mathrm{MOM}_K (\ell_{Z^*}-\ell_Z) + \lambda (\norm{Z^*}-\norm{Z}) & \leq  \max\left(\frac{-1}{5}, \frac{-57}{400}, \frac{13}{40}\right) r_2^2 = \frac{13}{40} (r_2^*)^2 
\end{align}
and so (\ref{eq:Lemma_eta_RMOM_G-loc_02}) holds for $\eta \geq  13(r_2^*)^2/(40)$. As a consequence, \eqref{eq:Lemma_eta_RMOM_G-loc_01} and (\ref{eq:Lemma_eta_RMOM_G-loc_02}) both hold for $\eta =  132(r_2^*)^2/(400)$.
\end{proof}

From Lemmas \ref{Lemma:proof_main_RMOM_G-loc-00} and \ref{Lemma:proof_main_RMOM_G-loc-01}, we conclude that on the event $\Omega_{K, G}$, $\hat Z \in \cB_{2}$, that is $\norm{\hat Z-Z^*}\leq 2\rho^*$ and $G(\hat Z - Z^*)\leq (r_2^*)^2$. The following lemma gives us an upper bound on the excess risk $P\cL_{\hat{Z}}$.

\begin{Lemma}\label{Lemma:proof_main_RMOM_G_loc_001}
Under the conditions of Theorem \ref{theo:main_RMOM_G-loc}, and on the event $\Omega_{K, G}$, we have $P\cL_{\hat{Z}} \leq (93/100)(r^*_2)^2$.
\end{Lemma}

\begin{proof}
From Lemmas \ref{Lemma:proof_main_RMOM_G-loc-00} and \ref{Lemma:proof_main_RMOM_G-loc-01}, we get that on $\Omega_{K, G}$, $\hat Z\in\cB_2$. This implies the existence of stricly more than $K/2$ blocks $B_k$ on which $\left|(P_{B_k}-P)\cL_{\hat{Z}}\right|\leq (r^*_2)^2/(20)$, that is:
\begin{align}\label{eq:proof_main_RMOM_G_loc_002}
P\cL_{\hat{Z}} \leq P_{B_k}\cL_{\hat{Z}} + (r^*_2)^2/(20).
\end{align}
Moreover, by (\ref{eq:Lemma_eta_RMOM_G-loc_02}), the definition of $\hat{Z}$ and \eqref{Lemma:proof_main_RMOM_G-loc-01}, we have:
\begin{align*}
\mathrm{MOM}_K(\ell_{\hat{Z}}-\ell_{Z^*}) + \lambda\left(\norm{\hat Z} - \norm{Z^*}\right) & \leq \underset{Z\in \cC}{\mathrm{sup}}\quad \mathrm{MOM}_K(\ell_{\hat{Z}}-\ell_{Z}) + \lambda\left(\norm{\hat Z} - \norm{Z}\right) \\ 
& \leq \underset{Z\in \cC}{\mathrm{sup}}\quad \mathrm{MOM}_K(\ell_{Z^*}-\ell_{Z}) + \lambda\left(\norm{Z^*} - \norm{Z}\right) \leq \frac{33}{100} (r_2^*)^2 .
\end{align*}
As a consequence, there exist at least $K/2$ blocks $B_k$ on which
\begin{align}\label{eq:proof_main_RMOM_G_loc_003}
P_{B_k}\cL_{\hat{Z}}   \leq \frac{33}{100} (r_2^*)^2 - \lambda\left(\norm{\hat Z} - \norm{Z^*}\right) \leq \frac{33}{100} (r_2^*)^2 + \lambda \norm{\hat Z - Z^*} \leq \frac{33}{100} (r_2^*)^2 + 2\lambda\rho^* = \frac{88}{100} (r_2^*)^2.
\end{align}
So there must be at least a block $B_{k_0}$ on which (\ref{eq:proof_main_RMOM_G_loc_002}) and (\ref{eq:proof_main_RMOM_G_loc_003}) hold simultaneously. On this block, we have
\begin{align*}
P\cL_{\hat{Z}} \leq P_{B_{k_0}}\cL_{\hat{Z}} + \frac{1}{20}(r_2^*)^2\leq \frac{88}{100} (r_2^*)^2 + \frac{1}{20}(r_2^*)^2 = \frac{93}{100}(r_2^*)^2.
\end{align*}
\end{proof}

At this stage, we have shown that on the event $\Omega_{K,G}$, the estimator $\hat{Z}$ has the statistical properties announced in Theorem \ref{theo:main_RMOM_G-loc}. In what follows we prove that under the conditions of Theorem \ref{theo:main_RMOM_G-loc}, $\Omega_{K,G}$ holds with exponentially large probability.

\begin{Lemma}\label{Lemma:main_MOM_G-loc_OmegaK_proba}
Assume that $K\geq 100|\cO|$, and let $\rho^*>0$ be such that it satisfies the sparsity equation from Definition \ref{def:sparsity_RMOM_G-loc}. Then, $\Omega_K$ holds with probability at least $1-2\exp(-72K/625)$.
\end{Lemma}

\begin{proof}
Let $\phi:t\in \bR \rightarrow \mathbb{1}_{\left\{t\geq 1\right\}}+2(t-1/2)\mathbb{1}_{\left\{1/2\leq t\leq 1\right\}}$, so that for any $t\in \bR$, $\mathbb{1}_{\left\{t\geq 1\right\}}\leq \phi(t) \leq \mathbb{1}_{\left\{t\geq 1/2\right\}}$. For $k\in [K]$, let $W_k:=\left\{X_i:i\in B_k\right\}$ and $F_Z(W_k)=(P_{B_k}-P)\cL_Z$. We also define the counterparts of these quantities constructed with the non-corrupted vectors: $\widetilde{W}_k:=\left\{\widetilde{X}_i:i\in B_k\right\}$ and $F_Z(\widetilde{W}_k)=(\widetilde{P_{B_k}}-\widetilde{P})\cL_Z$, where $\widetilde{P_{B_k}}\cL_Z := \frac{K}{N}\sum_{i\in B_k}\cL_Z(\widetilde{X}_i)$ and $\widetilde{P}\cL_Z := \bE[\cL_{Z}(\widetilde{X}_i)]$. For both $b\in\{1, 2\}$,  let $Z\to \psi_b(Z)=\sum_{k\in [K]}\mathbb{1}_{\left\{|F_Z(W_k)|\leq (r^*_b)^2/(20)\right\}}$. Let $b\in\{1, 2\}$. We want to show that, with high probability, if $Z\in \cB_{b}$, then $\psi_b(Z) > K/2$. As we showed in the proof of Lemma \ref{Lemma:main_MOM_loc_er_OmegaK_proba}, in our framework this is equivalent to proving that the following inequality occurs with high probability:
\begin{align}\label{eq:lemma_main_RMOM_Gloc_OmegaKG_01}
\sum_{k\in [K]}\mathbb{1}_{\left\{|F_Z(\widetilde{W}_k)| > \frac{(r_b^*)^2}{20}\right\}} \leq \frac{49K}{100},
\end{align}
and this is what we do now. Let  $Z\in \cB_{b}$. We have:
\begin{align}\label{eq:lemma_main_RMOM_Gloc_OmegaKG_02}
&\sum_{k\in \left[K\right]}\mathbb{1}_{\left\{|F_Z(\widetilde{W}_k)| > \frac{(r_b^*)^2}{20}\right\}} = \sum_{k\in \left[K\right]}\left[\mathbb{1}_{\left\{|F_Z(\widetilde{W}_k)| > \frac{(r_b^*)^2}{20}\right\}} - \bP\left(|F_Z(\widetilde{W}_k)| > \frac{(r_b^*)^2}{40}\right) + \bP\left(|F_Z(\widetilde{W}_k)| > \frac{(r_b^*)^2}{40}\right) \right] \nonumber\\
 & = \sum_{k\in \left[K\right]}\left(\mathbb{1}_{\left\{|F_Z(\widetilde{W}_k)| > \frac{(r_b^*)^2}{20}\right\}} - \bE\left[\mathbb{1}_{\left\{|F_Z(\widetilde{W}_k)| > \frac{(r_b^*)^2}{40}\right\}}\right] \right) + \sum_{k\in \left[K\right]} \bP\left(|F_Z(\widetilde{W}_k)| > \frac{(r_b^*)^2}{40}\right) \nonumber\\ 
 &\leq \sum_{k\in \left[K\right]}\left(\phi\left( \frac{20|F_Z(\widetilde{W}_k)|}{(r_b^*)^2}\right) -  \bE\left[\phi\left( \frac{20|F_Z(\widetilde{W}_k)|}{(r_b^*)^2}\right)\right]  \right) + \sum_{k\in \left[K\right]} \bP\left(|F_Z(\widetilde{W}_k)| > \frac{(r_b^*)^2}{40}\right) \nonumber\\ 
 &\leq \sup_{Z\in \cB_{b}} \left( \sum_{k\in \left[K\right]}\phi\left( \frac{20|F_Z(\widetilde{W}_k)|}{(r_b^*)^2}\right) -  \bE\left[\phi\left( \frac{20|F_Z(\widetilde{W}_k)|}{(r_b^*)^2}\right)\right]  \right) + \sum_{k\in \left[K\right]} \bP\left(|F_Z(\widetilde{W}_k)| > \frac{(r_b^*)^2}{40}\right)
\end{align}
We start with bounding the last sum in the previous inequality. For each $k\in\left[K\right]$, Markov's inequality and the definition of $r^*_b$ yield to
\begin{align*}
\bP\left(|F_Z(\widetilde{W}_k)| > \frac{(r_b^*)^2}{40}\right) & \leq \frac{1600^2}{(r_b^*)^4} \bE\left[F_Z(\widetilde{W}_k)^2\right] = \frac{1600^2}{(r_b^*)^4} \left(\frac{K}{N}\right)\mathrm{Var}(\cL_Z(\widetilde{X}))\leq   \frac{1600^2}{(r_b^*)^4} \left(V_K(r_b^*)\right)^2 \leq \frac{1}{200}
\end{align*}

Plugging this last result into (\ref{eq:lemma_main_RMOM_Gloc_OmegaKG_02}), we get:
\begin{align}\label{eq:lemma_main_RMOM_Gloc_OmegaKG_03}
\sum_{k\in \left[K\right]}\mathbb{1}_{\left\{|F_Z(\widetilde{W}_k)| > \frac{(r_b^*)^2}{20}\right\}} \leq \frac{K}{200}+ \sup_{Z\in \cB_{b}} \left( \sum_{k\in \left[K\right]}\phi\left( \frac{20|F_Z(\widetilde{W}_k)|}{(r_b^*)^2}\right) -  \bE\left[\phi\left( \frac{20|F_Z(\widetilde{W}_k)|}{(r_b^*)^2}\right)\right]  \right).
\end{align}

We now have to bound this last term. Using Mc Diarmind inequality (Theorem 6.2 in \cite{BoucheronLugosiMassart} with $t=12/25$), we get that with probability at least $1-\exp(-72K/625)$, for all $Z\in \cB_{b}$, 
\begin{align}\label{eq:lemma_main_RMOM_Gloc_OmegaKG_04}
\sum_{k\in \left[K\right]}\phi\left( \frac{20|F_Z(\widetilde{W}_k)|}{(r_b^*)^2}\right) -\bE\phi\left( \frac{20|F_Z(\widetilde{W}_k)|}{(r_b^*)^2}\right)\leq \frac{12K}{25} + \bE\left[\sup_{Z\in \cB_{b}}\sum_{k\in \left[K\right]}\phi\left( \frac{20|F_Z(\widetilde{W}_k)|}{(r_b^*)^2}\right) -\bE\phi\left( \frac{20|F_Z(\widetilde{W}_k)|}{(r_b^*)^2}\right)\right]. 
\end{align}
Let now $\epsilon_1, \ldots, \epsilon_{K}$ be Rademacher variables independant from the $\widetilde{X}_i$'s. By the symmetrization Lemma, we have:
\begin{align}\label{eq:lemma_main_RMOM_Gloc_OmegaKG_05}
\bE\left[\sup_{Z\in \cB_{b}}\sum_{k\in \left[K\right]}\phi\left( \frac{20|F_Z(\widetilde{W}_k)|}{(r_b^*)^2}\right) -\bE\left[\phi\left( \frac{20|F_Z(\widetilde{W}_k)|}{(r_b^*)^2}\right)\right]\right]   \leq 2\bE\left[\sup_{Z\in \cB_{b}}\sum_{k\in \left[K\right]}\epsilon_k\phi\left( \frac{20|F_Z(\widetilde{W}_k)|}{(r_b^*)^2}\right)\right].
\end{align}
As $\phi$ is Lipschitz with $\phi(0)=0$, we can use the contraction Lemma (see \cite{LedTal01}, chapter 4) to get that:
\begin{align}\label{eq:lemma_main_RMOM_Gloc_OmegaKG_05}
\bE\left[\sup_{Z\in \cB_{b}}\sum_{k\in \left[K\right]}\epsilon_k\phi\left( \frac{20|F_Z(\widetilde{W}_k)|}{(r_b^*)^2}\right)\right] &\leq 2\bE\left[\sup_{Z\in \cB_{b}}\sum_{k\in \left[K\right]}\epsilon_k\frac{20F_Z(\widetilde{W}_k)}{(r_b^*)^{2}}\right]= \frac{40}{(r_b^*)^{2}}\bE\left[\sup_{Z\in \cB_{b}}\sum_{k\in \left[K\right]}\epsilon_k(\widetilde{P_{B_k}}-\widetilde{P})\cL_Z\right]
\end{align}

Now, let $(\sigma_i)_{i=1,\ldots, N}$ be a family of Rademacher variables independant from the $\widetilde{X}_i$'s and the $\epsilon_i$'s. Using the symmetrization Lemma again, we get:
\begin{align*}
\bE\left[\sup_{Z\in \cB_{b}}\sum_{k\in \left[K\right]}\epsilon_k(\widetilde{P_{B_k}}-\widetilde{P})\cL_Z\right] & \leq 2 \bE\left[\sup_{Z\in \cB_{b}}\frac{K}{N}\sum_{i=1}^N\sigma_i\cL_Z(\widetilde{X}_i)\right] \leq  2KE_G(r_b^*, b\rho^*) \leq 2K\gamma (r_b^*)^2.
\end{align*}

Combining this with (\ref{eq:lemma_main_RMOM_Gloc_OmegaKG_04}), (\ref{eq:lemma_main_RMOM_Gloc_OmegaKG_05}) and (\ref{eq:lemma_main_RMOM_Gloc_OmegaKG_06}), we finally get that, with probability at least $1-\exp(-72K/625)$:
\begin{align}\label{eq:lemma_main_RMOM_Gloc_OmegaKG_06}
\sup_{Z\in \cC_{\gamma}}&\sum_{k\in \left[K\right]}\phi\left( \frac{20|F_Z(\widetilde{W}_k)|}{(r_b^*)^2}\right) -\bE\left[\phi\left( \frac{20|F_Z(\widetilde{W}_k)|}{(r_b^*)^2}\right)\right] \leq \left(\frac{12}{25}+160\gamma \right)K
\end{align}
Plugging that into (\ref{eq:lemma_main_RMOM_Gloc_OmegaKG_03}), we conclude that, with probability at least $1-\exp(-72K/625)$, for all $Z\in\cB_b$,
\begin{align*}
\sum_{k\in \left[K\right]}\mathbb{1}_{\left\{|F_Z(\widetilde{W}_k)| > \frac{(r_b^*)^2}{20}\right\}} \leq \left(\frac{1}{200}+\frac{12}{25}+160\gamma \right) K \leq \frac{49}{100}K
\end{align*}
for our choice of parameters. Now, in order for $\Omega_{K, G}$ to hold, this inequality must be verified for both $b=1$ and $2$. Then, we finally conclude that $\Omega_{K, G}$ holds with probability $1-2\exp(-72K/625)$, which concludes the proof.

\end{proof}

\subsection{Proofs of section \ref{sec:processus_sto}}\label{proofs:processus_sto}
\subsubsection{Proof of Theorem \ref{theo:main_processus_sto}}
The proof of Theorem~\ref{theo:main_processus_sto} relies on several Lemmas. We first recall that $k B_1\cap B_2\subset 2 {\rm conv}(U_{k^2}\cap \cS_2^{d\times d})$ where $\cS_2^{d\times d}$ is the unit sphere of $\ell_2^{d\times d}$ and $U_{k^2}$ is the set of all matrices in $\bR^{d\times d}$ with at most $k^2$ non zero entries (see, for instance, equation~(3.1) in \cite{MR2373017}). Hence, for all $A\in\bR^{d\times d}$, we have
\begin{equation*}
\norm{A}\leq 2\sup_{I\subset [d]\times [d]:|I|= k^2} \left(\sum_{(p,q)\in I}A_{pq}^2\right)^{1/2}.
\end{equation*}We therefore need to find a high probability upper bound on the $\ell_2$ norm of the $k^2$ largest entries of $\hat \Sigma_N -\Sigma$. To that end, we start with the following result.

\begin{Lemma}\label{lem:proc_sto_1}
Let $(z_{pq}:p,q\in[d])$ be real-valued random variables (not necessarily independent) and $\lambda,t\geq1$ be two positive constants. We assume that for $r =  2\log(ed/k)+t$, we have $\norm{z_{pq}}_{L_{r}}\leq \lambda \sqrt{r}$ for all $p,q\in[d]$. Then, with probability at least $1-\exp(-t)$, 
\begin{equation*}
\sup_{I\subset [d]\times [d]:|I| = k^2} \left(\sum_{(p,q)\in I}z_{pq}^2\right)^{1/2}\leq e^{2}\lambda \sqrt{2k^2\left(\log(ed/k)+t\right)}.
\end{equation*}
Moreover:
\begin{align*}
\bE\left[\sup_{I\subset [d]\times [d]:|I| = k^2} \left(\sum_{(p,q)\in I}z_{pq}^2\right)^{1/2}\right]\leq e^{2}\lambda \sqrt{6k^2\log(ed/k)}
\end{align*}
\end{Lemma}

\beginproof
We define for all $p,q\in[d]$,
\begin{equation*}
Z_{pq}=z_{pq} I(|z_{pq}|\leq e \lambda \sqrt{r}) \mbox{ and } Y_{pq}=z_{pq} I(|z_{pq}|> e \lambda \sqrt{r})
\end{equation*}so that we have $|z_{pq}|^{t} = |Z_{pq}|^{t} + |Y_{pq}|^{t}$. As a consequence and by convexity of $x\in\bR^+\to x^{t/2}$, we have for all $I\subset [d]\times[d]$
\begin{equation}\label{eq:prop_sto_Z_first}
\left(\frac{1}{|I|}\sum_{(p,q)\in I}z_{pq}^2\right)^{t/2}\leq \frac{1}{|I|}\sum_{(p,q)\in I}|z_{pq}|^{t} = \frac{1}{|I|}\sum_{(p,q)\in I}|Z_{pq}|^{t} + \frac{1}{|I|}\sum_{(p,q)\in I}|Y_{pq}|^{t}.
\end{equation}Let $I\subset [d]\times[d]$ be such that $|I| = k^2$. We have
\begin{equation}\label{eq:prop_sto_Z}
\frac{1}{|I|}\sum_{(p,q)\in I}|Z_{pq}|^{t} \leq (e\lambda\sqrt{r})^{t}.
\end{equation}For the second term in the right hand side inequality of \eqref{eq:prop_sto_Z_first}, we have
\begin{equation*}
\frac{1}{|I|}\sum_{(p,q)\in I}|Y_{pq}|^{t} \leq \frac{1}{|I|}\sum_{(p,q)\in [d]\times [d]}|Y_{pq}|^{t}
\end{equation*}and for $\theta:=r/t$ and all $p,q\in[d]$, we have
\begin{align*}
\bE[|Y_{pq}|^{t}]& = \bE\left[|z_{pq}|^{t} I(|z_{pq}|>e\lambda \sqrt{r})\right]\leq \bE\left[|z_{pq}|^{t\theta}\right]^{1/\theta}\bP\left[|z_{pq}|>e\lambda \sqrt{r}\right]^{1-1/\theta} \leq (\lambda \sqrt{r})^{t}\left(\frac{\norm{z_{pq}}_{L_{r}}}{e\lambda \sqrt{r}}\right)^{r-t}\\
&\leq (\lambda \sqrt{r})^{t} e^{r-t} = (\lambda \sqrt{r})^{t} \frac{k^2}{e^2d^2}.
\end{align*}It follows that 
\begin{equation*}
\bE \sup_{I\subset [d]\times [d]:|I|= k^2}\frac{1}{|I|}\sum_{(p,q)\in I}|Y_{pq}|^{t} \leq\frac{d^2}{k^2} (\lambda \sqrt{r})^{t} \frac{k^2}{e^{2}d^2} = \frac{(\lambda \sqrt{r})^{t}}{e^2}.
\end{equation*}Hence, using \eqref{eq:prop_sto_Z_first}, \eqref{eq:prop_sto_Z} and the last inequality, we get $\bE \cZ^{t}\leq  (e\lambda\sqrt{r})^{t} + (\lambda \sqrt{r})^{t}/e^2\leq 2(e\lambda \sqrt{r})^{t}$ where
\begin{equation*}
\cZ :=  \sup_{I\subset [d]\times [d]:|I|= k^2} \left(\frac{1}{|I|}\sum_{(p,q)\in I}z_{pq}^2\right)^{1/2}.
\end{equation*}As a consequence, $\norm{\cZ}_{L_{t}}\leq e\lambda \sqrt{2r}$ and so, for $t\geq 2$ we get by Markov's inequality that $\cZ\leq e^2\lambda \sqrt{2r}$ with probability at least $1-\exp(-t)$.

Finally, by taking  $t=1$ above we get:
\begin{align*}
\norm{\cZ}_{L_{1}} \leq e\lambda \sqrt{2r} = e\lambda \sqrt{2(2\log(ed/k)+1)} \leq e\lambda \sqrt{6\log(ed/k)}
\end{align*}
since $k\leq d$. As a consequence, $\bE[\cZ] = \norm{\cZ}_{L_{1}} \leq e\lambda \sqrt{6\log(ed/k)}$, which concludes the proof.
\endproof

The proof of Theorem~\ref{theo:main_processus_sto} will follow from Lemma~\ref{lem:proc_sto_1} if one can apply the latter to the variables $z_{pq}= \hat \Sigma_{pq}-\Sigma_{pq}$. We therefore have to check that $(\hat \Sigma_{pq}-\Sigma_{pq}: p,d\in[d])$ satisfies the assumptions of Lemma~\ref{lem:proc_sto_1}. In other words, it only remains to show that for all $p,q\in[d]$, $\hat \Sigma_{pq}-\Sigma_{pq}$ has $r :=  2\log(ed/k)+t$ sub-gaussian moment under Assumption~\ref{assum:weak_moment}. To that end we use a version (see Lemma~2.8 in \cite{LeM14}) of a result due to Lata{\l}a taken from \cite{latala1997estimation} (see Theorem~2 and Remark~2 in  \cite{latala1997estimation}) which states the following: 

\begin{Lemma} \cite{latala1997estimation}\label{lem:moment-sum-variid}
There exists an absolute constant $c_0$ for which the following holds. Let $z$ be a mean-zero random variable and $z_1,\ldots,z_N$ be $N$ independent copies of $z$. Let $p_0\geq 2$ and assume that there exists $\kappa_1>0$ and $\alpha \geq 1/2$ for which $\|z\|_{L_p}\leq \kappa_1 p^\alpha$ for every $2\leq p\leq p_0$. If $N\geq p_0^{\max\{2\alpha-1,1\}}$  then for every $2\leq p\leq p_0$,
    \begin{equation*}
      \Big\|\frac{1}{\sqrt{N}}\sum_{i=1}^N z_i\Big\|_{L_p}\leq c_1(\alpha) \kappa_1 \sqrt{p},
    \end{equation*}
  where $c_1(\alpha)=c_0\exp((2\alpha-1))$.
\end{Lemma}

We use Lemma~\ref{lem:moment-sum-variid} to prove the following moment growth condition on the $\hat \Sigma_{pq}-\Sigma_{pq}, p,q\in[d]$.

\begin{Lemma}\label{lem:proc_sto_2}There exists an absolute constant $c_0$ such that the following holds.
Grant Assumption~\ref{assum:weak_moment} with parameters $w$ and $t\geq2$. For all $p,q\in [d]$ and all $2\leq r \leq 2\log(ed/k)+t$, if $N\geq 2\log(ed/k)+t$ then  $\norm{\hat \Sigma_{pq}-\Sigma_{pq}}_{L_r}\leq (c_0 w^2/\sqrt{N}) \sqrt{r}$.
\end{Lemma}
\beginproof
Let $p,q\in [d]$. It follows from Assumption~\ref{assum:weak_moment}  and Lemma~\ref{lem:moment-sum-variid} that for all $p,q\in [d]$ and all $2\leq r \leq 2\log(ed/k)+t$,
\begin{equation*}
\norm{\hat \Sigma_{pq}-\Sigma_{pq}}_{L_r}\leq \frac{1}{\sqrt{N}} \norm{\frac{1}{\sqrt{N}}\sum_{i=1}^N X_{ip}X_{iq} - \bE X_{ip}X_{iq}}_{L_r}\leq \frac{c_0 w^2}{\sqrt{N}}\sqrt{r}. 
\end{equation*}
\endproof

\paragraph{Proof of Theorem~\ref{theo:main_processus_sto}} 
We set for all $p,q\in [d]$, $z_{pq} = \hat \Sigma_{pq}-\Sigma_{pq}$. It follows from Lemma~\ref{lem:proc_sto_2} for $\alpha=1$ that for all $2\leq r \leq 2\log(ed/k)+t$, $\norm{z_{pq}}_{L_r}\leq \lambda \sqrt{r}$ where $\lambda = c_0 w^2/\sqrt{N}$. The result now follows from Lemma~\ref{lem:proc_sto_1}.
\endproof

\subsubsection{Proof of Theorem \ref{theo:main_processus_sto_SLOPE}}
The proof of Theorem~\ref{theo:main_processus_sto} relies on several Lemmas. We first use a decomposition similar to the one from \cite{MR3782379}. We have
\begin{align}\label{eq:first_SLOPE}
\notag \norm{\hat \Sigma_N -\Sigma}_\rho &\leq \min\left(\sup_{Z\in B_2} \sum_{p,q=1}^d Z_{pq}(\hat \Sigma_N -\Sigma)_{pq}, \sup_{Z\in\rho B_{SLOPE}}\sum_{p,q=1}^d Z_{pq}(\hat \Sigma_N -\Sigma)_{pq}\right)\\
\notag&=\min\left(\sqrt{\sum_{p,q=1}^d(\hat \Sigma_N -\Sigma)_{pq}^2}, \rho\sup_{Z\in\rho B_{SLOPE}}\sum_{p,q=1}^d Z_{(p,q)}^*\beta_{pq}\frac{(\hat \Sigma_N -\Sigma)_{(p,q)}^*}{\beta_{pq}}\right)\\
&=\min\left(\sqrt{\sum_{p,q=1}^d(\hat \Sigma_N -\Sigma)_{pq}^2}, \rho\max_{p,q\in[d]}\frac{(\hat \Sigma_N -\Sigma)_{(p,q)}^*}{\beta_{pq}}\right).
\end{align}
We already proved a high probability upper bound on the $\ell_2$ norm of the $k^2$ largest entries of $\hat \Sigma_N -\Sigma$ in the previous section under a weaker assumption than the one in Assumption~\ref{assum:weak_moment_slope}. We just have to use it for $k=d$ to handle the left-hand side term of \eqref{eq:first_SLOPE}.  Therefore, with probability at least $1-\exp(-t)$,
\begin{equation*}
\sqrt{\sum_{p,q=1}^d(\hat \Sigma_N -\Sigma)_{pq}^2}\leq c_0 w^2 \sqrt{\frac{d}{N}}.
\end{equation*} It only remains to handle the second term in the right-hand side inequality of \eqref{eq:first_SLOPE}.  To that end, we start with the following result.

\begin{Lemma}\label{lem:proc_sto_1_slope}
Let $\textbf{z}:=(z_{pq}:p,q\in[d])$ be real-valued random variables (not necessarily independent) and $\lambda,t\geq1$ be two positive constants. We denote by $(z_{(p,q)}^*:p,q\in[d])$ the non-increasing sequence (for the same lexicographical order over $[d]^2$ used before) of the rearrangement of the absolute values of the entries of $\textbf{z}$. Let $p_0,q_0\in[d]$. We assume that for $r =  \log[ed^2/(p_0q_0)]+t$, we have $\norm{z_{pq}}_{L_{r}}\leq \lambda \sqrt{r}$ for all $p,q\in[d]$. Then, 
\begin{equation*}
\norm{\frac{z_{(p_0,q_0)}^*}{\beta_{p_0q_0}}}_{L_t}\leq e^2 \lambda.
\end{equation*}
\end{Lemma}

\beginproof To make the presentation of the proof simpler, we index the entries of $d\times d$ matrices by $[d^2]$. We therefore have $d^2$ random variables $(z_j)_j$ (not necessarily independent) and $\beta_j=\sqrt{\log(ed^2/j)+t}$ for all $j\in[d^2]$. Let $j_0\in[d^2]$ and set $r_0 = \log(ed^2/j_0)+t$. We assume that $\norm{z_j}_{L_{r_0}}\leq \lambda \sqrt{r_0}$ for all $p,q\in[d]$. We want to prove that $\norm{z_{j_0}^*/\beta_{j_0}}_{L_t}\leq e^2 \lambda$. We first remark that
\begin{equation}\label{eq:first_slope}
\frac{z_{j_0}^*}{\beta_{j_0}}\leq \max_{I\subset [d^2]:|I|=j_0}\frac{1}{\beta_{j_0} |I|}\sum_{j\in I} |z_j|:=\cZ.
\end{equation}

We define for all $j\in[d^2]$,
\begin{equation*}
Z_{j}=z_{j} I\left(|z_{j}|\leq e \lambda \sqrt{r_0}\right) \mbox{ and } Y_{j}=z_{j} I\left(|z_{j}|> e \lambda \sqrt{r_0}\right).
\end{equation*}It follows from the convexity of $x\in\bR_+\to x^t$ and the definitions above that
\begin{equation}\label{eq:split_Slope}
\bE\cZ^t\leq \max_{I\subset [d^2]:|I|=j_0}\frac{1}{\beta_{j_0}^t |I|}\sum_{j\in I} |z_j|^t\leq \left(\frac{e\lambda \sqrt{r_0}}{\beta_{j_0}}\right)^t+ \frac{1}{j_0}\sum_{j=1}^{d^2}\frac{\bE |Y_j|^t}{\beta_{j_0}^t}.
\end{equation}

Next, for the second term in the right-hand side inequality of \eqref{eq:split_Slope} for $\theta:=r_0/t$ and all $j\in[d^2]$, we have
\begin{align*}
\bE|Y_{j}|^{t}& = \bE\left[|z_{j}|^{t} I(|z_{j}|>e\lambda \sqrt{r_0})\right]\leq \bE\left[|z_{j}|^{t\theta}\right]^{1/\theta}\bP\left[|z_{j}|>e\lambda \sqrt{r_0}\right]^{1-1/\theta}\\ 
&\leq \left(\lambda \sqrt{r_0}\right)^{t}\left(\frac{\norm{z_{j}}_{L_{r_0}}}{e\lambda \sqrt{r_0}}\right)^{r_0-t}\leq (e\lambda)^t r_0^{t/2}e^{-r_0} = (e\lambda)^t \beta_{j_0}^t e^{-t} \frac{j_0}{ed^2}.
\end{align*}We end up in \eqref{eq:split_Slope} with $\bE\cZ^t\leq \left(e\lambda\right)^t+\lambda^t\leq (e^2\lambda)^t$.
\endproof

\begin{Lemma}\label{lem:slope_max}
Let $\textbf{z}:=(z_{pq}:p,q\in[d])$ be real-valued random variables (not necessarily independent) and $\lambda\geq0,t\geq3$ be two constants. We denote by $(z_{(p,q)}^*:p,q\in[d])$ the non-increasing sequence (for the same lexicographical order over $[d]^2$ used before) of the rearrangement of the absolute values of the entries of $\textbf{z}$. Let $r_0 =  \log(ed^2)+t$ and assume that $\norm{z_{pq}}_{L_{r}}\leq \lambda \sqrt{r}$ for all $p,q\in[d]$ and $2\leq r\leq r_0$. Let $k\in[d]$ and $\gamma\geq1$. Then, when $t\geq \max\big(2 \log(\lceil \log(k^2)\rceil), \gamma \log(ed^2/k^2)\big)$, with probability at least $1-\exp(-t/2)$,
\begin{equation*}
\max_{p,q\in[d^2]}\left(\frac{z_{(p,q)}^*}{\beta_{pq}}\right)\leq \sqrt{2}e^3 \lambda.
\end{equation*}
\end{Lemma}
\beginproof
We use the same 'vectorial' notation as the one introduced in the proof of Lemma~\ref{lem:proc_sto_1_slope}. We remark that for all $j\in[d^2]$, we have 
$(1/\beta_{2j})\leq\sqrt{2} /\beta_{j}$ when $t\geq3$ and for all $j\geq k^2$, $1/\beta_{j}\leq\sqrt{2}/\beta_{k^2}$ when $t\geq \gamma \log(ed^2/k^2)$,  hence, 
\begin{equation*}
\max_{j\in[d^2]}\left(\frac{z_j^*}{\beta_j}\right)\leq \sqrt{2}\max\left(\frac{z_{2^j}^*}{\beta_{2^j}}:j=0,1,\ldots,\lceil \log(k^2)\rceil\right).
\end{equation*}Il follows from Lemma~\ref{lem:proc_sto_1_slope} that for all $j=0,1,\ldots,\lceil \log(k^2)\rceil$, we have $\norm{z_{2^j}^*/\beta_{2^j}}_{L_t}\leq e^2 \lambda$ and so by Markov's inequality with probability at least $1-\exp(-t), z_{2^j}^*/\beta_{2^j}\leq e^3\lambda$. The union bound yields that with probability at least $1-\lceil \log(k^2)\rceil\exp(-t)$, $\max\left(z_{2^j}^*/\beta_{2^j}:j=0,1,\ldots,\lceil \log(k^2)\rceil\right)\leq e^3\lambda$.

\endproof

The proof of Theorem~\ref{theo:main_processus_sto} will follow from Lemma~\ref{lem:slope_max} if one can apply the latter to the variables $z_{pq}= \hat \Sigma_{pq}-\Sigma_{pq}$. We therefore have to check that the family of random variables $(\hat \Sigma_{pq}-\Sigma_{pq}: p,d\in[d])$ satisfies the assumptions of Lemma~\ref{lem:slope_max}. In other words, it only remains to show that for all $p,q\in[d]$, $\hat \Sigma_{pq}-\Sigma_{pq}$ has $r :=  \log(ed^2)+t$ sub-gaussian moment under Assumption~\ref{assum:weak_moment_slope}. To that end we use Lemma~\ref{lem:moment-sum-variid} to prove the following moment growth condition on the $\hat \Sigma_{pq}-\Sigma_{pq}, p,q\in[d]$.

\begin{Lemma}\label{lem:proc_sto_2_slope}There exists an absolute constant $c_0$ such that the following holds.
Grant Assumption~\ref{assum:weak_moment_slope} with parameters $w$ and $t\geq3$. For all $p,q\in [d]$ and all $2\leq r \leq 2\log(ed^2)+t$, if $N\geq 2\log(ed^2)+t$ then  $\norm{\hat \Sigma_{pq}-\Sigma_{pq}}_{L_r}\leq (c_0 w^2/\sqrt{N}) \sqrt{r}$.
\end{Lemma}
\beginproof
Let $p,q\in [d]$. It follows from Assumption~\ref{assum:weak_moment_slope}  and Lemma~\ref{lem:moment-sum-variid} that for all $p,q\in [d]$ and all $2\leq r \leq 2\log(ed^2)+t$,
\begin{equation*}
\norm{\hat \Sigma_{pq}-\Sigma_{pq}}_{L_r}\leq \frac{1}{\sqrt{N}} \norm{\frac{1}{\sqrt{N}}\sum_{i=1}^N X_{ip}X_{iq} - \bE X_{ip}X_{iq}}_{L_r}\leq \frac{c_0 w^2}{\sqrt{N}}\sqrt{r}. 
\end{equation*}
\endproof

\paragraph{Proof of Theorem~\ref{theo:main_processus_sto}} 
We set for all $p,q\in [d]$, $z_{pq} = \hat \Sigma_{pq}-\Sigma_{pq}$. It follows from Lemma~\ref{lem:proc_sto_2_slope} for $\alpha=1$ that for all $2\leq r \leq 2\log(ed^2)+t$, $\norm{z_{pq}}_{L_r}\leq \lambda \sqrt{r}$ where $\lambda = c_0 w^2/\sqrt{N}$. The result now follows from Lemma~\ref{lem:slope_max}.
\endproof

\subsection{Proofs of section \ref{sec:sparsePCA}}

\subsubsection{Proof of Lemma \ref{coro:exactness1_spiked}}

Let $Z\in\cC$ and consider its SVD $Z=\sum_i \sigma_i u_i u_i^\top$. We have
\begin{align*}
&\inr{\Sigma, (\beta^*)(\beta^*)^\top -Z} = \theta\inr{(\beta^*)(\beta^*)^\top, (\beta^*)(\beta^*)^\top-Z} + \inr{I_d, (\beta^*)(\beta^*)^\top-Z} \ove{=}{(i)} \theta\inr{(\beta^*)(\beta^*)^\top, (\beta^*)(\beta^*)^\top - \sum_i\sigma_i u_i u_i^\top} \\{}
&  = \theta\left(1-\sum_{i}\sigma_i \inr{u_i, \beta^*}^2\right)= \theta\sum_i\sigma_i (1-\inr{u_i, \beta^*}^2)\ove{\geq}{(ii)}0
\end{align*}
where we used in \textit{(i)} that $\inr{I_d, (\beta^*)(\beta^*)^\top-Z} = \Tr{(\beta^*)(\beta^*)^\top}-\Tr{Z}=0$, and in \textit{(ii)} that $|\inr{u_i,\beta^*}|\leq 1$ (by Cauchy-Schwart). Hence, $\beta^*(\beta^*)^\top$ is a solution to the problem $\max\left(\inr{\Sigma, Z}, Z\in\cC\right)$. Moreover, using the latter computation, it is straightforward to check that it is unique, that is $\sigma_1 = 1$ and $u_1u_1^\top = \beta^*(\beta^*)^\top$, otherwise inequality \textit{(ii)} would be strict.

\subsubsection{Proof of Lemma \ref{coro:curvature_excess_risk}}

Let $Z\in\cC$ and consider its SVD $Z=\sum_i \sigma_i u_i u_i^\top$. In the proof of Lemma~\ref{coro:exactness1_spiked}, we proved that 
\begin{equation*}
 \inr{\Sigma, Z^* -Z} =  \theta\sum_i\sigma_i(1-\inr{u_i, \beta^*}^2).
  \end{equation*}On the other-hand, we have
  \begin{align*}
    &\norm{Z^*-Z}_2^2 = \Tr{(Z^*-Z)(Z^*-Z)^\top} = \Tr{\left(\sum_i\sigma_i((\beta^*)(\beta^*)^\top - u_i u_i^\top)\right)^2}\\
    &=\sum_{i,j}\sigma_i \sigma_j \Tr{(u_i u_i^\top - (\beta^*)(\beta^*)^\top)(u_j u_j^\top - (\beta^*)(\beta^*)^\top)}=\sum_i\sigma_i\left(\sigma_i - 2\inr{u_i, \beta^*}^2 + 1\right)\\
    &=2\sum_i\sigma_i(1-2\inr{u_i, \beta^*}^2)+ \left(\sum_i \sigma_i^2 - \sigma_i\right) = \frac{2}{\theta}\inr{\Sigma, Z^*-Z} + \left(\norm{Z}_2^2 - \norm{Z}_*\right)\leq \frac{2}{\theta}\inr{\Sigma, Z^*-Z}.
    \end{align*}

\subsubsection{Proof of Lemma \ref{Lemma:rho-spars}}\label{proof:lemma_sparsity_sparsePCA_RMOM_l1}

It follows from the $k$-sparsity of $\beta^*$ that $Z^* = \beta^*(\beta^*)^\top$ is $k^2$-sparse. Let us denote $I:=\mathrm{supp}(Z^*)$: we have $|I|\leq k^2$. Consider $\rho>0$. To solve the sparsity equation, we will use the following result on the sub-differential of a norm: if $\norm{.}$ is a norm over $\bR^{d\times d}$, we have for $Z\in \bR^{d\times d}$: 
\begin{align*}
\partial\norm{.}(Z) = \left\{ \begin{array}{cl}
\{\Phi \in S^*: \inr{\Phi, Z} =  \norm{Z}\}& \mbox{if } Z\ne 0\\
B^*& \mbox{ if } Z=0
\end{array}\right.
\end{align*}
where $S^*$ (resp. $B^*$) is the unit-sphere (resp. unit-ball) for the dual norm associated with $\norm{.}$, that is $Z\in \bR^{d\times d}\to \norm{Z}^* = \sup_{\norm{H} = 1}\inr{Z, H}$. Here, we consider the $\ell_1$-norm, whose dual norm is the $\ell_{\infty}$ norm.\\

Since $Z^*\in Z^*+(\rho/20)B$, we have $$\partial\norm{.}_1(Z^*)\subset \Gamma_{Z^*}(\rho) := \underset{V \in Z^*+(\rho/20)B }{\bigcup}\partial\norm{.}_1(V).$$ Then, there exists $\Phi^*\in \Gamma_{Z^*}(\rho) $ which is norming for $Z^*$, that is $\norm{\Phi^*}_{\infty} = 1$ and $\inr{\Phi^*, Z^*}=\norm{Z^*}_1$. Let $Z\in H_{\rho, A}:= Z^*+(\rho S_1\cap\sqrt{r^*(\rho)}B_2)$. For $J\subset [d]^2$, let $P_{J}$ be the coordinate projection on $J$. Since the supports of $P_{I^c}Z$ and $Z^*$ are disjoints, we can choose $\Phi^*$ such that it is also norming for $P_{I^c}Z$. Then, we have:
\begin{align*}
\inr{\Phi^*, Z-Z^*} &= \inr{\Phi^*, P_I(Z-Z^*)} + \inr{\Phi^*, P_{I^c}(Z-Z^*)} \geq -|\inr{\Phi^*, P_I(Z-Z^*)}| + \norm{P_{I^c}Z}_1 \\ 
& \geq -\norm{\Phi^*}_{\infty}\norm{P_I(Z-Z^*)}_1+ \norm{P_{I^c}Z}_1  = -\norm{P_I(Z-Z^*)}_1 +  \norm{P_{I^c}(Z-Z^*)}_1 \\
& = \norm{Z-Z^*}_1 - 2\norm{P_I(Z-Z^*)}_1 = \rho - 2\norm{P_I(Z-Z^*)}_1.
\end{align*}
Now, we have $\norm{P_I(Z-Z^*)}_1\leq k \norm{P_I(Z-Z^*)}_2 \leq k \norm{Z-Z^*}_2 \leq k\sqrt{r^*(\rho)}$. We conclude that $\inr{\Phi^*, Z-Z^*} \geq \rho - 2k\sqrt{r^*(\rho)}$. Then, $\underset{\Phi \in \Gamma_{Z^*}(\rho) }{\mathrm{sup}} \inr{\Phi, Z-Z^*}\geq \inr{\Phi^*, Z-Z^*}\geq \rho - 2k\sqrt{r^*(\rho)}$. Since this is true for any $Z\in H_{\rho, A}$, we conclude that $c\geq \rho - 2k\sqrt{r^*(\rho)}$, where $P_I(Z-Z^*)$ is the quantity introduced in Definition \ref{def:RERM_A-sparsity}. Then, if we choose $\rho$ such that $\rho\geq 10k\sqrt{r^*(\rho)}$, we have $\Delta(\rho, A)\geq (4/5)\rho$, and the $A$-sparsity equation is satisfied by such a $\rho$.

\subsubsection{Proof of Lemma \ref{Lemma:r2_bound}}

From Lemma \ref{coro:curvature_excess_risk}, we get that Assumption \ref{ass:A-spars} holds with $G:Z\in\bR^{d\times d}\rightarrow \norm{Z}_2^2$ and $A=2/\theta$, for any $\rho>0$ and $\delta\in(0, 1)$. Moreover, Assumption \ref{assum:weak_moment} is granted for $t = \log(ed/10k)$ and $w\geq 0$. Let then $c_0>0$ be the constant provided by Theorem \ref{theo:main_processus_sto}, and define $b = 3c_0w^2$. Let us define the following function:
\[
r : \rho > 0 \to bA \sqrt{\frac{\rho^2}{N}\log\left(\frac{b^2A^2(ed)^4}{N\rho^2}\right)}.
\]
We also consider 
\[
\rho^* := 200Ab k^2\sqrt{\frac{1}{N}\log\left(\frac{ed}{k}\right)}, 
\]
as well as $r^* = r(\rho^*)$. We have:
\begin{align}\label{eq:rG_bound_RERM_L1_001}
100k^2 r^* = 100k^2 bA\sqrt{\frac{(\rho^*)^2}{N}\log\left(\frac{(ed)^4}{4.10^4 k^4\log(\frac{ed}{k})}\right)} \leq   200 k^2 bA\sqrt{\frac{(\rho^*)^2}{N}\log\left(\left(\frac{ed}{k}\right)^4\right)}  = (\rho^*)^2,
\end{align}
so that $\rho^* \geq 10k \sqrt{r^*}$. Let us then define $k^* := \rho^*/\sqrt{r^*}$. Since $k^*>k$, any $2\leq r \leq 2\log(ed/k^*) + t$ satisfies $2\leq r\leq 2\log(ed/k) + t$, so that Assumption \ref{assum:weak_moment} holds with $w$, $t$ and $k^*$. We are then in measure to apply Theorem \ref{theo:main_processus_sto} with those parameters. As a consequence, as soon as $N\geq 2\log(ed/k^*) + t$, one has with probability at least $1-\exp(-t)$:
\begin{align}\label{eq:rG_bound_RERM_L1_01}
\norm{\hat \Sigma_N - \Sigma}_{k^*} \leq c_0w^2\sqrt{\frac{(k^*)^2\left(\log\left(\frac{ed}{k^*}\right)+t\right)}{N}}
\end{align}
where $\norm{.}_{k^*}$ is the $\ell_1/\ell_2$ interpolation norm defined in (\ref{eq:norm_interpol_mat}). Now, we have:
\begin{align}\label{eq:rG_bound_RERM_L1}
\sup_{Z\in \cC\cap(Z^*+\rho^* B_1\cap \sqrt{r^*}B_2)}|(P-P_N)\cL_Z| & \leq \sqrt{r^*} \sup_{Z\in \cC\cap(Z^*+k^* B_1\cap B_2)}|(P-P_N)\cL_Z|\nonumber\\ 
& = \sqrt{r^*} \bignorm{\frac{1}{N}\sum_{i=1}^NX_iX_i^{\top}-\bE[X_iX_i^{\top}]}_{k^*} = \sqrt{r^*} \norm{\hat \Sigma_N - \Sigma}_{k^*}. 
\end{align}
Combining it with (\ref{eq:rG_bound_RERM_L1_01}), we get that with probability at least $1-\exp(-t)$:
\begin{align}\label{eq:rG_bound_RERM_L1_02}
\sup_{Z\in \cC\cap(Z^*+\rho^* B_1\cap \sqrt{r^*}B_2)}|(P-P_N)\cL_Z| & \leq c_0w^2\sqrt{\frac{(\rho^*)^2\left(\log\left(\frac{ed}{k^*}\right)+t\right)}{N}} \leq c_0w^2\sqrt{\frac{2(\rho^*)^2}{N}\log\left(\frac{ed}{10k}\right)}
\end{align}
since $k^*\geq 10k$. Now, we have:
\begin{align*}
r^* &= bA \sqrt{\frac{(\rho^*)^2}{N}\log\left(\frac{b^2A^2(ed)^4}{N(\rho^*)^2}\right)}  = bA \sqrt{\frac{(\rho^*)^2}{N}\log\left(\frac{(ed)^4}{200^2k^4\log(ed/10k)}\right)} \\ 
& \geq bA \sqrt{\frac{(\rho^*)^2}{N}\log\left(\frac{(ed)^3}{200^2k^3}\right)} \geq bA \sqrt{\frac{(\rho^*)^2}{N}\log\left(\frac{ed}{10k}\right)}
\end{align*}
where the last inequality holds since we assumed that $k\leq ed/200$. Combining it with (\ref{eq:rG_bound_RERM_L1_02}), we conclude that:
\begin{align*}
\sup_{Z\in \cC\cap(Z^*+\rho^* B_1\cap \sqrt{r^*}B_2)}|(P-P_N)\cL_Z| & \leq \frac{r^*}{3A}
\end{align*}
which allows us to conclude that $r^*_{RERM,G}(A, \rho^*, e^{-t}) \leq r^*$. Moreover, we have from (\ref{eq:rG_bound_RERM_L1_001}) that
\begin{align*}
\rho^* \geq 10k\sqrt{r^*} \geq 10k\sqrt{r^*_{RERM,G}(A, \rho^*, e^{-t}) }
\end{align*}
that is, $\rho^*$ satisfies the $A$-sparsity equation from Definition \ref{def:RERM_A-sparsity}. These results are valid provided that $N\geq 2\log(ed/k^*) + t$, which is ensured by the assumption that $N\geq 3\log(ed/10k)$, given that $k^*\geq 10k$. This concludes the proof.

\subsubsection{Proof of Theorem \ref{theo:main_RERM_SPCA}}

From Lemma \ref{coro:curvature_excess_risk}, we get that Assumption \ref{ass:A-spars} holds with $G:Z\in\bR^{d\times d}\rightarrow \norm{Z}_2^2$ and $A=2/\theta$, for any $\rho>0$ and $\delta\in(0, 1)$. Moreover, since we assumed that $N\geq 3\log(ed/10k)$, Lemma \ref{Lemma:r2_bound} applies and so, for $\rho^*$ and $r^*(\rho^*)$ (defined in (\ref{eq:rstar_RERM_L1})) we have $r^*_{\mathrm{RERM, G}}(A, \rho^*, 10k/ed) \leq r^*(\rho^*)$ and $\rho^*$ satisfies the $A$-sparsity equation from Definition \ref{def:RERM_A-sparsity}. We are then in position to apply Theorem \ref{theo:main-penalized}, provided that $\lambda$ satisfies (\ref{hyp:lambda}). Now, we have:
\begin{align*}
\frac{r^*(\rho^*)}{\rho^*} = bA \sqrt{\frac{1}{N}\log\left(\frac{(bA)^2(ed)^4}{N(\rho^*)^2}\right)} = bA \sqrt{\frac{1}{N}\log\left(\frac{(ed)^4}{200^2k^4\log(\frac{ed}{k})}\right)},
\end{align*}
so that:
\[
bA \sqrt{\frac{3}{N}\log\left(\frac{ed}{200^{2/3}k}\right)} \leq \frac{r^*(\rho^*)}{\rho^*} \leq bA \sqrt{\frac{4}{N}\log\left(\frac{ed}{200^{1/2}\log(200)^{1/4}k}\right)},
\]since we assumed that $k\leq ed/200$. As a consequence, (\ref{hyp:lambda}) is satisfied as soon as:
\[
\frac{20}{21}b\sqrt{\frac{1}{N}\log\left(\frac{ed}{200^{1/2}\log(200)^{1/4}k}\right)}  \leq \lambda \leq \frac{2}{\sqrt{3}}b \sqrt{\frac{1}{N}\log\left(\frac{ed}{200^{2/3}k}\right)}
\]
which is the assumption made in (\ref{eq:main_RERM_L1_cond_lambda}). We only have to check that this authorized interval for $\lambda$ is not empty, which is ensured as soon as $ed/k\geq 200^{48/47}/\log(200)^{25/47}$, which is granted by the assumption that $k\leq ed/200$. \\

We are then in measure to apply Theorem \ref{theo:main-penalized}, which enables us to state that, with probability at least $1-10k/ed$:
\begin{align*}
 \norm{\hat{Z}^{\mathrm{RERM}}_{\lambda}-Z^*}_1 & \leq \rho^* = 200Ab k^2\sqrt{\frac{1}{N}\log\left(\frac{ed}{k}\right)} = 400b k^2\sqrt{\frac{1}{N\theta^2}\log\left(\frac{ed}{k}\right)}, \nonumber \\ 
 \norm{\hat{Z}^{\mathrm{RERM}}_{\lambda}-Z^*}_2 & \leq \sqrt{r^*_{\mathrm{RERM, G}}\left(A, \rho^*, 10k/ed\right)} \leq \frac{\rho^*}{10k} = 40b \sqrt{\frac{k^2}{N\theta^2}\log\left(\frac{ed}{k}\right)} \\  
\end{align*}  and
\begin{equation*}
P\cL_{\hat{Z}^{\mathrm{RERM}}_{\lambda}}  \leq A^{-1}r^*_{\mathrm{RERM, G}}\left(A, \rho^*, 10k/ed\right)\leq \frac{(\rho^*)^2}{100k^2A} = 800b^2\frac{k^2}{N\theta}\log\left(\frac{ed}{k}\right).
\end{equation*}
This concludes the proof.

\subsubsection{Proof of Corollary \ref{coro:main_RERM_L1}}

From Theorem \ref{theo:main_RERM_SPCA}, we get the existence of a universal constant $C > 0$ such that with probability at least $1-20\left(k/ed\right)^{3/4}$, $\norm{\hat{Z}^{\mathrm{RERM}}_{\lambda}-Z^*}_2  \leq C\sqrt{k^2(N\theta^2)\log(ed/k)}$. Now, we can use Davis-Kahan sin-theta theorem (see Corollary~1 in \cite{yu2014useful}) to get the existence of a universal constant $c_0>0$ such that $\mathrm{sin}(\Theta(\hat{\beta},\beta^*)) = (1/\sqrt{2})\norm{\hat{\beta}\hat{\beta}^\top-\beta^*(\beta^*)^\top}_2  \leq (c_0/g) \norm{\hat{Z}^{\mathrm{RERM}}_{\lambda}-Z^*}_2 $ where $g := \lambda_1 - \lambda_2$ ($\lambda_i$ being the $i^{\mathrm{th}}$ largest eigen value of $Z^*$) is the spectral gap of $Z^*$. Here, we know that $Z^*=\beta^*(\beta^*)^\top$ is rank one, with $1$ as order one eigen value and $0$ as order $(d-1)$ eigen value. Then we get $g=1$, which leads us to the desired result, with $D=\sqrt{2}c_0\times C$.

\subsubsection{Proof of Lemma \ref{Lemma:SPCA_SLOPE_A-sparsity}}

Let $A$, $\delta$ and $t>0$. In the rest of the proof, we write $r^*_G(.)$ for $r^*_{\mathrm{RERM,G}}(A, ., \delta)$, $b_{pq}$ for $b_{pq}(t)$ and $\Gamma_k$ for $\Gamma_k(t)$. We consider a lexicographical order on $[d]^2$, $b\in\bR^{d\times d}$ and the norm $\norm{.}_{SLOPE}$ as they are defined in section \ref{sec:RERM_SLOPE_sparsePCA}.

Let $I := \textup{supp}\big((Z^*)^{\sharp}\big)$ be the set of non-zero coefficients of $(Z^*)^{\sharp}$. Since $Z^* = \beta^* (\beta^*)^{\top}$ is $k^2$-sparse, whe have by construction that $|I|\leq k^2$. Let $P_I$ (resp. $P_{I^{c}}$) be the coordinate projection on $I$ (resp. on $I^{c}$). \smallskip

We know that for $Z\neq 0$:
\begin{align*}
\partial\norm{.}_{SLOPE}(Z)=\left\{\Phi \in S_{SLOPE}^*:\inr{\Phi,Z} = \norm{Z}_{SLOPE}\right\},
\end{align*}
where we denoted $S_{SLOPE}^*$ the unit-sphere of the dual norm of the $SLOPE$ norm. Since $Z^*\in Z^*+\frac{\rho}{20}B_{SLOPE}$, we know that $\partial\norm{.}_{SLOPE}(Z^*)\subset \Gamma_{Z^*}(\rho)$. Then:
\begin{align*}
\sup_{\Phi \in \Gamma_{Z^*}(\rho)} \inr{\Phi,Z-Z^* } \geq \sup_{\Phi \in \partial\norm{.}_{SLOPE}(Z^*)} \inr{\Phi,Z-Z^* }.
\end{align*}

Let $\sigma$, $\pi$ be the permutations of $[d]^2$ such that, for any $(p,q)\in [d]^2$, $(Z^*)^{\sharp}_{p,q}= |Z^*_{\sigma(p,q)}|$ and $(Z-Z^*)^{\sharp}_{p,q}= |(Z-Z^*)_{\pi(p,q)}|$. Notice that we have by assumption $\sigma([k]^2)=I$. We then define $\Phi^*$ and $\widetilde{\Phi}^*$ as follows: for all $1\leq p,q\leq d$,
\begin{equation*}
\Phi^*_{p,q}=\begin{cases} \mathrm{sgn}(Z^*_{p,q}) \,b_{\sigma^{-1}(p,q)}\quad \textup{if }(p,q)\leq (k,k); \\
\mathrm{sgn}((Z-Z^*)_{p,q}) \,b_{\pi^{-1}(p,q)} \quad \textup{otherwise}\end{cases}
\end{equation*}and
\begin{equation*}
\widetilde{\Phi}^*_{p,q}= \mathrm{sgn}((Z-Z^*)_{p,q}) \,b_{\pi^{-1}(p,q)}.
\end{equation*}

We easily check that such a $\Phi^*$ belongs to $\partial ||\cdot ||_{SLOPE}(Z^*)$ and $\widetilde{\Phi}^*$ to $\partial ||\cdot ||_{SLOPE}(Z-Z^*)$. Now let $Z\in H_{\rho, A}$. We have:
\begin{equation}\label{eq:ineg:sharp}
\begin{split}
\left\langle\Phi^*,Z-Z^* \right\rangle &= \left\langle\Phi^*,P_{I}(Z-Z^*) \right\rangle + \left\langle\Phi^*,P_{I^{c}}(Z-Z^*) \right\rangle\\
&= \sum\limits_{p,q =1}^{k} \mathrm{sgn}(Z^*_{\sigma(p,q)}) b_{p,q} (Z-Z^{*})_{\sigma(p,q)} + \left\langle\Phi^*,P_{I^{c}}(Z-Z^*) \right\rangle.
\end{split}
\end{equation}
Regarding the first term, we have:
\[
\left|\sum\limits_{p,q=1}^{k}\mathrm{sgn}(Z^*_{\sigma(p,q)}) b_{p,q} (Z-Z^{*})_{\sigma(p,q)}\right| \leq \sum\limits_{p,q=1}^{k} b_{p,q}  \left|(Z-Z^{*})_{\pi(p,q)}\right| = \sum\limits_{p,q=1}^{k} b_{p,q} (Z-Z^*)^{\sharp}_{p,q},
\]
where the first inequality comes from the fact that the operator $(\cdot)^{\sharp}$ orders the absolute values of $(Z-Z^*)$ in non-increasing order (notice that the inequality holds only for the sum, not for each independent term of the sum). Therefore:
\begin{equation}\label{eq:ineg:sharp:1}
\sum\limits_{p,q=1}^{k}\mathrm{sgn}(Z^*_{\sigma(p,q)}) b_{p,q} (Z-Z^{*})_{\sigma_{p,q}}\geq -\sum\limits_{p,q=1}^{k} b_{p,q} (Z-Z^*)^{\sharp}_{p,q}.
\end{equation}
Concerning the second term in \eqref{eq:ineg:sharp}:
\begin{equation}\label{eq:ineg:sharp:2}\begin{split}
\left\langle\Phi^*,P_{I^{c}}(Z-Z^*) \right\rangle &= \left\langle\widetilde{\Phi}^*,P_{I^{c}}(Z-Z^*) \right\rangle= \left\langle\widetilde{\Phi}^*, Z-Z^* \right\rangle-\left\langle\widetilde{\Phi}^*,P_{I}(Z-Z^*) \right\rangle\\
&=||Z-Z^*||_{SLOPE}- \sum\limits_{p,q=1}^{k} b_{\pi^{-1}\circ \sigma(p,q)}(Z-Z^*)^{\sharp}_{\pi^{-1}\circ \sigma(p,q)}\\ 
&\geq ||Z-Z^*||_{SLOPE}- \sum\limits_{p,q=1}^{k} b_{p,q}(Z-Z^*)^{\sharp}_{p,q}.
\end{split}\end{equation}
Putting \eqref{eq:ineg:sharp}, \eqref{eq:ineg:sharp:1} and \eqref{eq:ineg:sharp:2} together, we obtain
\begin{align}\label{eq:ineg:sharp:2bis}
\left\langle\Phi^*,Z-Z^* \right\rangle  \geq ||Z-Z^*||_{SLOPE} - 2 \sum\limits_{p,q=1}^{k} b_{p,q} (Z-Z^*)^{\sharp}_{p,q} = \rho - 2 \sum\limits_{p,q=1}^{k} b_{p,q} (Z-Z^*)^{\sharp}_{p,q}.
\end{align}
Now, since $\norm{Z-Z^*}_2\leq \sqrt{r_G^*}$, we can show that for any $k\in [d]$, $(Z-Z^*)^{\sharp}_{kk}\leq \frac{\sqrt{r_G^*}}{k}$. Indeed, assume the existence of $k_0\in [d]$ such that $(Z-Z^*)^{\sharp}_{k_0k_0} > \frac{\sqrt{r_G^*}}{k_0}$. Then by construction we have that for any $(p, q)\leq (k_0, k_0)$, $(Z-Z^*)^{\sharp}_{pq} \geq (Z-Z^*)^{\sharp}_{k_0k_0}$, so that
\[
\norm{Z-Z^*}_2^2 = \norm{(Z-Z^*)^\sharp}_2^2 \geq \sum_{(p,q)\leq (k_0,k_0)}((Z-Z^*)^\sharp_{pq})^2 > \sum_{(p,q) \leq (k_0,k_0)} \frac{r_G^*}{k_0^2} = r_G^*,
\]
 since the $k_0^2$ largest elements of $(Z-Z^*)^\sharp$ belong to $[k_0]^2$, as a result of which $|\left\{(p,q):\,(p,q)\leq(k_0,k_0)\right\}|\leq k_0^2$. This is inconsistent with the fact that $\norm{Z-Z^*}_2\leq \sqrt{r_G^*}$.\smallskip

As a consequence, we have:
\begin{align*} 
&\sum_{p,q = 1}^k b_{pq}(Z-Z^*)^{\sharp}_{pq}  = \sum_{\ell = 1}^{k-1} \sum_{(\ell, \ell) \leq (p,q) < (\ell+1,\ell+1)}b_{pq}(Z-Z^*)^{\sharp}_{\ell\ell} + b_{kk}(Z-Z^*)^{\sharp}_{kk}\\
  &\leq \sum_{\ell= 1}^{k-1} \left|\left\{(\ell, \ell)\leq (p,q)<(\ell+1, \ell+1)\right\}\right|b_{\ell\ell}(Z-Z^*)^{\sharp}_{\ell\ell} + b_{kk} \frac{\sqrt{r_G^*}}{k}\\
  &\leq \sum_{\ell = 1}^{k-1} (2\ell+1)b_{\ell\ell}\frac{\sqrt{r_G^*}}{\ell} + b_{kk} \frac{\sqrt{r_G^*}}{k} \leq 3 \sqrt{r_G^*}\sum_{\ell = 1}^{k-1} b_{\ell\ell} + b_{kk} \frac{\sqrt{r_G^*}}{k} 
 \leq  3 \sqrt{r_G^*}\sum_{\ell = 1}^{k} b_{\ell\ell} = \sqrt{r_G^*}\Gamma_k.
\end{align*}

Then, under the assumption that $\rho \geq 10\Gamma_k\sqrt{r_G^*(A, \rho, \delta)}$, we get from (\ref{eq:ineg:sharp:2bis}) that $\inr{\Phi^*, Z-Z^*} \geq (4/5)\rho$.
and then: 
\begin{align*}
\sup_{\Phi \in S_{SLOPE}^*} \inr{\Phi, Z-Z^*}\geq \inr{\Phi^*, Z-Z^*} \geq \frac{4}{5}\rho.
\end{align*} Since this is true for any $Z\in H(\rho,A)$, we conclude that:
\begin{align*}
\Delta(\rho,A) = \inf_{Z\in H(\rho,A)}\sup_{\Phi \in S_{SLOPE}^*} \inr{\Phi^*, Z-Z^*} \geq \frac{4}{5}\rho.
\end{align*}
that is, $\rho$ satisfies the $A$-sparsity equation from Definition \ref{def:RERM_A-sparsity}.

\subsubsection{Proof of Lemma \ref{Lemma:SPCA_SLOPE_fixed-point}}

From Lemma \ref{coro:curvature_excess_risk}, we get that Assumption \ref{ass:A-spars} holds with $G:Z\in\bR^{d\times d}\rightarrow \norm{Z}_2^2$ and $A=2/\theta$, for any $\rho>0$ and $\delta\in(0, 1)$. 

For $r$ and $\rho>0$, we define $\cC_{r, \rho}:=\left\{Z\in \cC : \norm{Z-Z^*}_{SLOPE}\leq \rho, \norm{Z-Z^*}_2\leq \sqrt{r}\right\}$. Let $A>0$. For any $\rho$ and $r>0$. We have 
\begin{align}\label{proof:Lemma_SPCA_SLOPE_fixed_point_01}
&\sup_{Z\in \cC_{r, \rho}}|\inr{\Sigma-\hat{\Sigma}_N, Z-Z^*}| \leq \sup_{Z\in (\rho B_{SLOPE}\cap \sqrt{r}B_2)} |\inr{\Sigma-\hat{\Sigma}_N, Z}| \nonumber \\
 &= \sqrt{r}\sup_{Z\in (\frac{\rho}{\sqrt{r}} B_{SLOPE}\cap B_2)} |\inr{\Sigma-\hat{\Sigma}_N, Z}| 
= \sqrt{r} \norm{\Sigma-\hat{\Sigma}_N}_{\frac{\rho}{\sqrt{r}}}
\end{align}
where $\norm{.}_{\rho/\sqrt{r}}$ is the $SLOPE$/$\ell_2$ interpolation norm defined in (\ref{eq:norm_interpol_mat_slope}). Assumption \ref{assum:weak_moment_slope} is granted for $t = 2\log(ed^2/k^2)$. Let us now check that $k\leq d/(e^2\log(d))$: we have that $k^2\log(ek^2) \leq ed^2$, hence,
\begin{align*}
2\log(\lceil\log(k^2)\rceil) & \leq 2\log(\log(k^2)+1) = 2\log(\log(ek^2)) \leq 2 \log\left(\frac{ed^2}{k^2}\right),
\end{align*}
that is, $t\geq \max\big(2 \log(\lceil \log(k^2)\rceil), 2 \log(ed^2/k^2)\big)$. We are then in position to apply Theorem \ref{theo:main_processus_sto_SLOPE} with $\gamma = 2$ and $t = 2\log(ed^2/k^2)$: there exists a universal constant $c_0>0$ such that, provided that $N\geq \log\left(ed^2\right)+t$, one has with probability at least $1-2\exp(-t/2)$:
\[
\norm{\Sigma-\hat{\Sigma}_N}_{\frac{\rho}{\sqrt{r}}} \leq \frac{c_0w^2}{\sqrt{N}}\min\left(\frac{\rho}{\sqrt{r}}, d\right).
\]
Plugging this last result into (\ref{proof:Lemma_SPCA_SLOPE_fixed_point_01}), we get that:
\begin{align}\label{proof:Lemma_SPCA_SLOPE_fixed_point_02}
\sup_{Z\in \cC_{r, \rho}}|\inr{\Sigma-\hat{\Sigma}_N, Z-Z^*}| \leq \sqrt{r}\frac{c_0w^2}{\sqrt{N}}\min\left(\frac{\rho}{\sqrt{r}}, d\right)
\end{align}
with probability at least $1-2\exp(-t/2)$. Next, let us define $b:=3c_0w^2$ and for $\rho>0$, consider 
\begin{align*}
r^*(\rho) := \frac{bA}{\sqrt{N}}\min\left(bA\frac{d^2}{\sqrt{N}};\rho\right).
\end{align*}
One can check that for this choice of  $r^*$, one has $(\sqrt{r^*(\rho)}c_0w^2/\sqrt{N})\min\left(\rho/\sqrt{r^*(\rho)}, d\right) \leq r^*(\rho)/3A$ whatever the value of $\rho$ is. From ($\ref{proof:Lemma_SPCA_SLOPE_fixed_point_02}$) we then deduce that $r^*_{\mathrm{RERM, G}}(A, \rho, 2e^{-t/2})\leq r^*(\rho)$. Let us now consider 
\begin{align*}
\rho^*:=10\Gamma_k^* \frac{bA}{\sqrt{N}}\min\left(10\Gamma_k^*;d\right), 
\end{align*}
where  $\Gamma_k^*:=3\sum_{\ell = 1}^{k} b_{\ell\ell}(t)$. It is straighforward to verify that $\rho^* \geq 10\Gamma_k^* r^*(\rho^*)^{1/2}\geq 10\Gamma_k^* r^*_{\mathrm{RERM, G}}\left(A, \rho^*, 2e^{-t^*/2}\right)^{1/2}$ which, according to Lemma \ref{Lemma:SPCA_SLOPE_A-sparsity}, guarantees that $\rho^*$ satisfies the $A$-sparsity equation from Definition \ref{def:RERM_A-sparsity}. Finally, plugging the expression of $\rho^*$ into the one of $r^*(\rho^*)$, we get that $r^*(\rho^*) = (b^2A^2/N)\min\left(d, 10\Gamma_k^*\right)^2$. Finally, the previous results hold provided that $N\geq \log\left(ed^2\right)+t$, which is granted by the assumption that $N\geq 3\log(ed^2)$. This concludes the proof, noting that $2\exp(-t/2) = 2k^2/(ed^2)$.

\subsubsection{Proof of Theorem \ref{theo:main_SPCA_RSLOPE} }

From Lemmas \ref{coro:exactness1_spiked} and \ref{coro:curvature_excess_risk}, we get that Assumption \ref{ass:A-spars} holds with $G : Z\rightarrow \norm{Z}_2^2$ and $A=2/\theta$. From Lemma \ref{Lemma:SPCA_SLOPE_fixed-point}, we get the existence of a constant $b>0$ such that, provided that $N\geq 3\log(ed^2)$, defining $\rho^*:=10\Gamma_k^*(bA/\sqrt{N})\min\left(10\Gamma_k^*;d\right)$ and $r^* = (b^2A^2/N)\min\left(d, 10\Gamma_k^*\right)^2$, with $\Gamma_k^*=\Gamma_k(2\log(ed^2/k^2))$, one has $r^*_{\mathrm{RERM, G}}(A, \rho^*, 2k^2/ed^2) \leq r^*$ and $\rho^*$ satsifies the $A$-sparsity equation from Definition \ref{def:RERM_A-sparsity}. Let us now upper bound $\Gamma_k^*$:
\begin{align}\label{eq:Gammak_bound_01}
\Gamma_k^* = 3\sum_{\ell =1}^k b_{\ell\ell}\left(2\log\left(\frac{ed^2}{k^2}\right)\right) \leq 3 \left(\sum_{\ell =1}^k \sqrt{\log\left(\frac{ed^2}{\ell^2}\right)} + \sum_{\ell =1}^k \sqrt{2\log\left(\frac{ed^2}{k^2}\right)}\right) \nonumber \\ 
\leq 3 \sum_{\ell =1}^k \sqrt{\log\left(\frac{ed^2}{\ell^2}\right)} + 3k\sqrt{2\log\left(\frac{ed^2}{k^2}\right)}. 
\end{align} 
Concerning the first term in this last inequality, we have:
\begin{align}\label{eq:Gammak_bound_02}
\left(\sum_{\ell =1}^k \sqrt{\log\left(\frac{ed^2}{\ell^2}\right)} \right) ^2 & = 2 \sum_{m<\ell}\sqrt{\log\left(\frac{ed^2}{\ell^2}\right)}\sqrt{\log\left(\frac{ed^2}{m^2}\right)} + \sum_{\ell=1}^k\log\left(\frac{ed^2}{\ell^2}\right) \nonumber \\ 
& \leq 2 \sum_{m<\ell} \log\left(\frac{ed^2}{\ell^2}\right) + \sum_{\ell=1}^k\log\left(\frac{ed^2}{\ell^2}\right) \leq 3k \sum_{\ell=1}^k\log\left(\frac{ed^2}{\ell^2}\right). 
\end{align}
Moreover, we have:
\begin{align}\label{eq:Gammak_bound_03}
\sum_{\ell=1}^k\log\left(\frac{ed^2}{\ell^2}\right) &\leq \sum_{\ell=1}^k\int_{u=\ell-1}^{\ell}\log\left(\frac{ed^2}{u^2}\right)\mathrm{d}u = \int_{u=0}^{k}\log\left(\frac{ed^2}{u^2}\right) =k\log(ed^2) - 2 \left[u\log(u)-u\right]_0^k \nonumber \\ 
&= k\log\left(\frac{ed^2}{k^2}\right) + 2k \leq 3 k\log\left(\frac{ed^2}{k^2}\right) .
\end{align}
Combining (\ref{eq:Gammak_bound_01}), (\ref{eq:Gammak_bound_02}) and (\ref{eq:Gammak_bound_03}), we finally get that $\Gamma_k^*\leq (9+3\sqrt{2})k\sqrt{\log\left(\frac{ed^2}{k^2}\right)}\leq 14 k\sqrt{\log\left(ed^2/k^2\right)}$. As a consequence, we have
\begin{align*}
10\Gamma_k^*\leq 140 k \sqrt{\log\left(\frac{ed^2}{k^2}\right)}\leq 140 k \sqrt{2\log\left(d\right)} \leq 140 k \sqrt{\frac{2d}{e^2k}} \leq d,
\end{align*}
since we assumed that $k\leq \min\left(d/(e^2\log(d)), (e/(140\sqrt{2}))^2d\right)$. We conclude that $\min\left(10\Gamma_k^*, d\right) = 10\Gamma_k^* \leq  140 k\sqrt{\log\left(ed^2/k^2\right)}$. Plugging this result into the expression of $r^*$ and $\rho^*$, we finally get that:
\begin{align*}
r^*\leq 140^2b^2A^2\frac{k^2}{N}\log\left(\frac{ed^2}{k^2}\right) \quad \mbox{and} \quad \rho^*\leq 140^2bA\frac{k^2}{\sqrt{N}}\log\left(\frac{ed^2}{k^2}\right)
\end{align*}
so that $r^*/\rho^* = bA/\sqrt{N}$. As a consequence, (\ref{hyp:lambda}) is satisfied as soon as:
\[
\frac{10b}{21\sqrt{N}} < \lambda < \frac{2b}{3\sqrt{N}}
\]
which is (\ref{eq:lambda_RMOM_SLOPE}). We are then in position to apply Theorem \ref{theo:main-penalized}, which allows us to conclude that, with probability at least $1-2k^2/ed^2$:
\begin{align*}
\norm{\hat{Z}^{RERM}_{\lambda}-Z^*}_{SLOPE} \leq \rho^* ~~,~~ G(\hat{Z}^{RERM}_{\lambda}-Z^*) \leq r^* ~~and ~~P\cL_{\hat{Z}^{RERM}_{\lambda}} \leq A^{-1}r^*.
\end{align*}
This concludes the proof.

\subsubsection{Proof of Corollary \ref{coro:main_SPCA_RSLOPE}}
The proof follows exactly the same lines as the one of Corollary \ref{coro:main_RERM_L1}, so we do not detail it here.

\subsubsection{Proof of Lemma \ref{Lemma:rho-spars-RMOM-SPCA}}\label{proof:Lemma_rho-spars-RMOM-SPCA}
Consider $A=2/\theta $ and $\gamma >0$. In the rest of the proof we write $r^*(\rho)$ for $r^*_{\mathrm{RMOM, G}}(A, \gamma, \rho)$. For any $J\subset [d]^2$, let $P_J$ be the coordinate projection on $J$. Consider $\rho >0$. Let $I:=\mathrm{supp}(Z^*)$ be the set of non-zero coefficients of $Z^*$. From Lemma \ref{coro:exactness1_spiked}, we have that $|I|\leq k^2$. Moreover, we know that for any $Z\neq 0$, $\partial \norm{.}_1(Z)=\left\{\Phi \in S_{\infty}:\inr{\Phi, Z}= \norm{Z}_1\right\}$, where $S_{\infty}$ is the unit-sphere for $\norm{.}_{\infty}$. Since $Z^*\in Z^*+\frac{\rho}{20}B_1$, we have that $\partial \norm{.}_1(Z^*)\subset \Gamma_{Z^*}(\rho)=\underset{Z\in Z^*+\frac{\rho}{20}B_1}{\cup}\partial \norm{.}_1(Z)$. Let then $\Phi^* \in \partial \norm{.}_1(Z^*)$. Consider $Z\in \bar{H}_{\rho, A} := \left\{Z\in\cC:\norm{Z-Z^*}_1=\rho\mbox{ and }\norm{Z-Z^*}_2\leq \sqrt{2/\theta} r^*(\rho)\right\}$. Since $Z^*$ and $P_I^c(Z)$ have disjoint supports, we can choose $\Phi^*$ so that it is also norming for $P_I^c(Z)$. Then, we have:
\begin{align}\label{proof:Lemma-rho-spars-RMOM-SPCA-01}
\inr{\Phi^*, Z-Z^*} &= \inr{\Phi^* , P_I(Z-Z^*)} +  \inr{\Phi^* , P_I^c(Z-Z^*)} \geq - \norm{\Phi^*}_{\infty}\norm{P_I(Z-Z^*)}_1 +  \inr{\Phi^* , P_I^c(Z)}\nonumber \\
& =  -(\norm{Z-Z^*}_1 - \norm{P_I(Z-Z^*)}_1) + \norm{P_I^c(Z-Z^*)}_1 \nonumber\\
& =  2 \norm{P_I^c(Z-Z^*)}_1 - \norm{Z-Z^*}_1  = \norm{Z-Z^*}_1 - 2 \norm{P_I(Z-Z^*)}_1 = \rho - 2 \norm{P_I(Z-Z^*)}_1
\end{align}
where we used the fact that $ \norm{\Phi^*}_{\infty}=1$. Then, since $Z\in\bar{H}_{\rho, A}$, we have:
\begin{align}\label{proof:Lemma-rho-spars-RMOM-SPCA-02}
\norm{P_I(Z-Z^*)}_1 \leq k\norm{P_I(Z-Z^*)}_2 \leq k\norm{Z-Z^*}_2\leq k  \sqrt{\frac{2}{\theta}}  r^*(\rho)
\end{align}
Combining (\ref{proof:Lemma-rho-spars-RMOM-SPCA-01}) and (\ref{proof:Lemma-rho-spars-RMOM-SPCA-02}), we finally get that:
\begin{align}\label{proof:Lemma-rho-spars-RMOM-SPCA-03}
\inr{\Phi^*, Z-Z^*} \geq \rho -2 k  \sqrt{\frac{2}{\theta}}  r^*(\rho).
\end{align}
As a consequence, $\sup_{\Phi\in\Gamma_{Z^*}(\rho)} \inr{\Phi, Z-Z^*} \geq \inr{\Phi^*, Z-Z^*} \geq  \rho -2 k  \sqrt{2/\theta}  r^*(\rho)$. This being true whatever $Z\in \bar{H}_{\rho, A}$, it follows that $\bar{\Delta}(\rho) \geq  \rho -2 k  \sqrt{2/\theta}  r^*(\rho)$. We conclude that any $\rho$ such that $\rho\geq 10k\sqrt{2/\theta}  r^*(\rho)$ satisfies $\bar{\Delta}(\rho) \geq (4/5)\rho$.

\subsubsection{Proof of Lemma \ref{Lemma:r2-MOM}}\label{proof:Lemma_r2-MOM}

Consider $\gamma > 0$. From Lemma \ref{coro:curvature_excess_risk}, we get that Assumption \ref{ass:curvature_RMOM_G-loc} holds with $G:Z\in\bR^{d\times d}\rightarrow (\theta/2) \norm{Z}_2^2$ and $A=1$, for any $\gamma>0$, in particular for the value of $\gamma$ we have just set. Moreover, Assumption \ref{assum:weak_moment} is granted for $t=1$ and $w\geq 0$. Let then $c_0>0$ be the constant provided by Theorem \ref{theo:main_processus_sto}, and consider $B := 3c_0w^2$ and $D:=1600w^2$. Let us define the following function:
\begin{align*}
r : (\gamma, \rho) \to \max\left(\sqrt{\frac{B\rho}{\gamma}}\left(\frac{6}{N}\log\left(\frac{2B(ed)^2}{\gamma\theta\rho}\sqrt{\frac{6}{N}}\right)\right)^{1/4} ; D\sqrt{\frac{K}{N\theta}}\right).
\end{align*}
We also consider 
\[
\rho^* := \max\left(400\sqrt{3}B\frac{k^2}{\gamma}\sqrt{\frac{1}{N\theta^2}\log\left(\frac{ed}{k}\right)}; 10Dk \sqrt{\frac{2K}{N\theta^2}}\right),
\]
as well as $r^*(\gamma) = r(\gamma, \rho^*)$. One can check that $\rho^*$ such defined satisfies both of the two conditions below:
\begin{align}\label{eq:proof_rstar_RMOM_l1_sparsePCA_002}
\mathrm{(1)} \quad \rho \geq  10k \sqrt{\frac{2B\rho}{\theta\gamma}}\left(\frac{6}{N}\log\left(\frac{2B(ed)^2}{\gamma\theta\rho}\sqrt{\frac{6}{N}}\right)\right)^{1/4} \quad \mbox{and} \quad
\mathrm{(2)} \quad \rho \geq  10k D\sqrt{\frac{2K}{N\theta^2}},
\end{align}
so that $\rho^*\geq 10k\sqrt{2/\theta}r^*$. Let us define $k^* = \sqrt{\theta/2}\rho^*/r^*$. We have $\log(ed/k^*)+1\leq \log(ed/10k)+1$, so that Assumption \ref{assum:weak_moment} still holds with $w$, $t=1$ and $k^*$. Then, since we assumed that $N\geq 2\log(ed/k) + 1\geq 2\log(ed/k^*) + 1$, Theorem \ref{theo:main_processus_sto} applies and allows us to affirm that
\begin{align}\label{eq:proof_rstar_RMOM_l1_sparsePCA_001}
\bE\left[\bignorm{\frac{1}{N}\sum_{i=1}^N\widetilde{X}_i\widetilde{X}_i^{\top}-\bE[\widetilde{X}_i\widetilde{X}_i^{\top}]}_{k^*}\right] \leq c_0 w^2 \sqrt{\frac{6(k^*)^2\log(ed/k^*)}{N}},
\end{align} 
where $\norm{\cdot}_{k^*}$ is the $\ell_1/\ell_2$ interpolation norm defined in (\ref{eq:norm_interpol_mat})  for $k=k^*$. 

For $r$ and $\rho>0$, define $\cC_{r, \rho} := \left\{Z\in\cC:\norm{Z-Z^*}_1=\rho\mbox{ and }\norm{Z-Z^*}_2\leq \sqrt{2/\theta} r^*(\rho)\right\}$. Let us now upper bound $E_G(r^*, \rho^*)$ and $V_{K, G}(r^*, \rho^*)$ from Definition~\ref{def:G-loc-fixed-point-RMOM}.

\paragraph{Bounding the complexity term $E_G(r^*, \rho^*)$.} Let $\sigma_1, \ldots, \sigma_N$ be $i.i.d.$ rademacher variables independent from the $\widetilde{X}_i$'s. We have $\cC_{r^*, \rho^*}\subset \sqrt{2/\theta}r^* \left(k^* B_1\cap B_2\right)$. As a consequence: 

\begin{align}\label{eq:Lemma-r2-MOM-proof-01}
\underset{Z\in\cC_{r^*,\rho^*}}{\mathrm{sup}}\left|\frac{1}{N}\sum_{i = 1}^N\sigma_i\cL_Z(\widetilde{X}_i) \right| & \leq \sqrt{\frac{2}{\theta}} r^* \underset{Z\in (k^*B_1 \cap B_2)\cap(\cC-Z^*)}{\mathrm{sup}} \left|\frac{1}{N} \sum_{i = 1}^N \sigma_i\cL_Z(\widetilde{X}_i)\right| \nonumber \\ 
& = \sqrt{\frac{2}{\theta}} r^* \underset{Z\in (k^*B_1 \cap B_2)\cap(\cC-Z^*)}{\mathrm{sup}} \left|  \biginr{\frac{1}{N}\sum_{i = 1}^N \sigma_i\widetilde{X}_i\widetilde{X}_i^\top, Z}\right|.  
\end{align}
Now, it follows from the desymmetrization inequality (see Theorem 2.1 in \cite{KoltchinskiiSaintFlour}) that:
\begin{align}\label{eq:proof_rstar_RMOM_l1_sparsePCA_01}
&\bE\left[\underset{Z\in (k^*B_1 \cap B_2)\cap(\cC-Z^*)}{\mathrm{sup}} \left|  \biginr{\frac{1}{N}\sum_{i = 1}^N \sigma_i\widetilde{X}_i\widetilde{X}_i^\top, Z}\right|\right]\nonumber \\  
&\leq 2 \bE\left[\underset{Z\in (k^*B_1 \cap B_2)\cap(\cC-Z^*)}{\mathrm{sup}} \left|  \biginr{\frac{1}{N}\sum_{i = 1}^N \widetilde{X}_i\widetilde{X}_i^\top - E[\widetilde{X}_i\widetilde{X}_i^\top], Z}\right|\right]
 + \frac{2}{\sqrt{N}} \underset{Z\in (k^*B_1 \cap B_2)\cap(\cC-Z^*)}{\mathrm{sup}}\left|\inr{\bE\left[\widetilde{X}\widetilde{X}^\top\right], Z}\right| \nonumber \\ 
 \leq &2 \bE\left[\bignorm{\frac{1}{N}\sum_{i = 1}^N \widetilde{X}_i\widetilde{X}_i^\top - E[\widetilde{X}_i\widetilde{X}_i^\top]}_{k^*}\right] + \frac{2}{\sqrt{N}} \underset{Z\in (k^*B_1 \cap B_2)\cap(\cC-Z^*)}{\mathrm{sup}}\left|\biginr{\bE\left[\widetilde{X}\widetilde{X}^\top\right], Z}\right| \nonumber \\ 
 &\leq 2 c_0 w^2 \sqrt{\frac{6(k^*)^2\log(ed/k^*)}{N}} + \frac{2}{\sqrt{N}} \underset{Z\in (k^*B_1 \cap B_2)\cap(\cC-Z^*)}{\mathrm{sup}}\left|\biginr{\bE\left[\widetilde{X}\widetilde{X}^\top\right], Z}\right|,
\end{align}
where we used (\ref{eq:proof_rstar_RMOM_l1_sparsePCA_001}) in the last inequality. \\ 

Concerning the second term in (\ref{eq:proof_rstar_RMOM_l1_sparsePCA_01}), we have for any $Z\in (k^*B_1 \cap B_2)\cap(\cC-Z^*)$:
\[
\inr{\bE\left[\widetilde{X}\widetilde{X}^\top, Z\right]} = \inr{\theta\beta^*(\beta^*)^\top + Id, Z} \overset{(i)}{=} \theta\inr{\beta^*(\beta^*)^\top, Z} \overset{(ii)}{\leq} \theta \norm{\beta^*}_2^2\norm{Z}_2 \leq \theta
\]
where we used the fact that $\inr{Id, Z} = \mathrm{Tr}(Z) = 0$ in $(i)$ and Cauchy-Schwarz in $(ii)$. as a consequence:
\begin{align}\label{eq:proof_rstar_RMOM_l1_sparsePCA_02}
\underset{Z\in (k^*B_1 \cap B_2)\cap(\cC-Z^*)}{\mathrm{sup}}\left|\inr{\bE\left[\widetilde{X}\widetilde{X}^\top\right], Z}\right| \leq \theta.
\end{align}

Combining (\ref{eq:Lemma-r2-MOM-proof-01}), (\ref{eq:proof_rstar_RMOM_l1_sparsePCA_01}) and (\ref{eq:proof_rstar_RMOM_l1_sparsePCA_02}), we finally get that:
\begin{align}\label{eq:last_equa_paper}
\nonumber E_G(r^*, \rho^*) &= \bE\left[\underset{\cC_{r^*,\rho^*}}{\mathrm{sup}}\left|\frac{1}{N}\sum_{i = 1}^N\sigma_i\cL_Z(\widetilde{X}_i) \right|\right] \leq \sqrt{\frac{2}{\theta}} r^* \left(2c_0 w^2 \sqrt{\frac{6(k^*)^2\log(ed/k^*)}{N}} + \frac{2\theta}{\sqrt{N}}\right)\nonumber \\ 
& \leq \sqrt{2\theta} r^* \left(3c_0 w^2 \sqrt{\frac{6(k^*)^2\log(ed/k^*)}{\theta^2N}}\right),
\end{align}
where we used the assumption that $\theta \leq k\leq k^*$.

\paragraph{Bounding the variance term $V_{K,G}(r^*, \rho^*)$.} 
Let us now  upper bound the variance term $$V_{K,G}(r^*, \rho^*) = \sqrt{\frac{K}{N}}\underset{Z\in\cC_{r^*,\rho^*}}{\mathrm{sup}} \sqrt{\bV ar(\cL_Z(\widetilde{X}_i))}.$$ 

For $\widetilde{X}$ distributed as the $\widetilde{X}_i$'s and $Z\in \cC_{r^*,\rho^*}$, one has:
\begin{align*}
\bV ar(\cL_Z(\widetilde{X})) & = \bE[((\cL_Z(\widetilde{X}) - P(\cL_Z(\widetilde{X}))^2] = \bE[\inr{\widetilde{X}\widetilde{X}^\top - \bE[\widetilde{X}\widetilde{X}^\top], Z-Z^*}^2] \nonumber\\
& = \sum_{p, q, s, t = 1}^d \bE\left[(\widetilde{X}^{(p)}\widetilde{X}^{(q)} - \bE[\widetilde{X}^{(p)}\widetilde{X}^{(q)}])(\widetilde{X}^{(s)}\widetilde{X}^{(t)} - \bE[\widetilde{X}^{(s)}\widetilde{X}^{(t)}])\right](Z-Z^*)_{pq}(Z-Z^*)_{st} \nonumber\\
& = \sum_{p, q, s, t = 1}^d T_{p, q, s, t}(Z-Z^*)_{pq}(Z-Z^*)_{st}
\end{align*}
where we defined $T_{p, q, s, t} := \bE\left[(\widetilde{X}^{(p)}\widetilde{X}^{(q)} - \bE[\widetilde{X}^{(p)}\widetilde{X}^{(q)}])(\widetilde{X}^{(s)}\widetilde{X}^{(t)} - \bE[\widetilde{X}^{(s)}\widetilde{X}^{(t)}])\right]$ for all $1\leq p,q,s,t\leq d$. Remembering that Assumption \ref{assum:weak_moment} is granted, we have:
\begin{align*}
T_{p, q, s, t}= \left\{ \begin{array}{ll}
\norm{(\widetilde{X}^{(p)})^2 - \bE[(\widetilde{X}^{(p)})^2]}_{L_2}^2 \leq (2w^2)^2 & \mbox{if } p=q=s=t \\
\norm{\widetilde{X}^{(p)}\widetilde{X}^{(q)} - \bE[\widetilde{X}^{(p)}\widetilde{X}^{(q)}]}_{L_2}^2 \leq (2w^2)2 & \mbox{if } (p, q) = (s, t), p\neq q \\
0 \mbox{ otherwise}.
\end{array}\right.
\end{align*}
Then:
\begin{align*}
\bV ar(\cL_Z(\widetilde{X})) & \leq  \sum_{p = 1}^d4w^4 (Z-Z^*)_{pp}^2+ \sum_{p \neq q} 4w^4 (Z-Z^*)_{pq}^2 \\ 
& =  \sum_{p, q=1}^q 4w^4 (Z-Z^*)_{pq}^2 = 4w^4 \norm{Z-Z^*}_2^2 \leq (8 w^4/\theta) (r^*)^2.
\end{align*}
This being true for any $Z\in\cC_{\rho^*, r^*}$ and any $\widetilde{X}$ distributed as the $\widetilde{X}_i$'s, we conclude that:
\begin{align}\label{eq:Lemma-r2-MOM-proof-03}
V_{K,G}(r, \rho) \leq 2w^2\sqrt{\frac{2K}{N\theta}} r^*.
\end{align}
Combining \eqref{eq:last_equa_paper} and \eqref{eq:Lemma-r2-MOM-proof-03}, we finally get that:
\begin{align*}
\max\left(\frac{E_G(r^*, \rho^*)}{\gamma}, 400\sqrt{2}V_{K,G}(r^*, \rho^*)\right) \leq \max\left(\frac{B}{\gamma}\sqrt{\frac{6(\rho^*)^2}{N}\log\left(\sqrt{\frac{2}{\theta}}\frac{edr^*}{\rho^*}\right)}, D \sqrt{\frac{K}{N\theta}}r^*\right)
\end{align*}
Now, one can  check that $r^*$ satisfies both of the two conditions below:
\begin{align*}
\mathrm{(3)}  \quad\frac{B}{\gamma} \sqrt{\frac{6(\rho^*)^2}{N}\log\left(\sqrt{\frac{2}{\theta}}\frac{edr^*}{\rho^*}\right)} \leq (r^*)^2 \quad \mbox{and}\quad \mathrm{(4)} \quad   D\sqrt{\frac{K}{N\theta}} r^* \leq (r^*)^2
\end{align*}
Then, we have:
\begin{align*}
\max\left(\frac{E_G(r^*, \rho^*)}{\gamma}, 400\sqrt{2}V_{K,G}(r^*, \rho^*)\right) \leq (r^*)^2
\end{align*}
which, according to Definition \ref{def:G-loc-fixed-point-RMOM}, allows us to conclude that $r^*_{\mathrm{RMOM, G}}(\gamma, \rho^*)\leq r^*$. Moreover, we have from \eqref{eq:proof_rstar_RMOM_l1_sparsePCA_002} that $\rho^*\geq \sqrt{2/\theta}10kr^*\geq \sqrt{2/\theta}10k r^*_{\mathrm{RMOM, G}}(\gamma, \rho^*) $, that is, $\rho^*$ satisfies the sparsity equation from Definition \ref{def:sparsity_RMOM_G-loc}. This concludes the proof.

\subsubsection{Proof of Theorem \ref{theo:main-RMOM-SPCA}}
The assumptions of Lemma \ref{Lemma:r2-MOM} are met, which gives us the existence of two positive constants $B$ and $D$ such that, defining $\rho^* := \max\left(400\sqrt{3}Bk^2\gamma^{-1}\sqrt{(N\theta^2)^{-1}\log\left(ed/k\right)}; 10Dk \sqrt{2K(N\theta^2)^{-1}}\right)$ and \\  $r^*(\gamma, \rho):= \max\left(\sqrt{B\rho\gamma^{-1}}\left((6/N)\log\left(2B(ed)^2(\gamma\theta\rho)^{-1}\sqrt{(6/N)}\right)\right)^{1/4} ; D\sqrt{K/(N\theta)}\right)$, one has $r^*_{\mathrm{RMOM, G}}(\gamma, \rho^*)\leq r^*(\gamma, \rho^*)$ and $\rho^*$ satisfies the sparsity equation from Definition \ref{def:sparsity_RMOM_G-loc}. From Lemma \ref{coro:curvature_excess_risk}, we get that Assumption \ref{ass:curvature_RMOM_G-loc} holds with $G:Z\in\bR^{d\times d}\rightarrow (\theta/2)\norm{Z}_2^2$ and $A=1$ for any $\gamma>0$, as a result of which the validity conditions of Theorem \ref{theo:main_RMOM_G-loc} are met. Then, fixing $\gamma = 1/32000$ and defining $\lambda = (11r^*(\gamma, 2\rho^*))/(40\rho^*)$, it is true that with probability at least $1-2\exp(-72K/625)$, 
\begin{align}\label{eq:proof_theo:main-RMOM-SPCA_01}
\norm{\hat{Z}^{\mathrm{RMOM}}_{K, \lambda}-Z^*}_1\leq 2\rho^*, \, P\cL_{\hat{Z}^{\mathrm{RMOM}}_{K, \lambda}} \leq \frac{93}{100}(r^*(\gamma, 2\rho^*))^2 \,\mbox{ and } \norm{\hat{Z}^{\mathrm{RMOM}}_{K, \lambda}-Z^*}_2\leq \sqrt{\frac{2}{\theta}} r^*(\gamma, 2\rho^*).
\end{align}
Now, we can write:
\begin{align}\label{eq:proof_theo:main-RMOM-SPCA_03}
\rho^*  \leq D_1 \frac{k}{\sqrt{N\theta^2}} \max\left(k\sqrt{\log\left(\frac{ed}{k}\right)} ; \sqrt{K}\right)
\end{align}
with $D_1:=\max(400\sqrt{3}B\gamma^{-1}, 10\sqrt{2}D)$. On the other hand, since $d\geq k$, we get that:
\begin{align*}
\rho^* \geq D_2k^2\sqrt{\frac{1}{N\theta^2}\log\left(\frac{ed}{k}\right)} \geq D_2\frac{k^2}{\sqrt{N\theta^2}}
\end{align*}
where $D_2 := 400\sqrt{3}B\gamma^{-1}$. As a consequence, we have:
\begin{align*}
\log\left(\frac{B(ed)^2}{\gamma\theta\rho*}\sqrt{\frac{6}{N}}\right) & \leq \log\left(\frac{B(ed)^2}{\gamma\theta}\sqrt{\frac{6}{N}}\frac{\sqrt{N\theta^2}}{D_2k^2}\right) = \log\left(\frac{1}{4\sqrt{5}}\left(\frac{ed}{k}\right)^2\right) \leq 2 \log\left(\frac{ed}{k}\right),
\end{align*}
so that:
\begin{align}\label{eq:proof_theo:main-RMOM-SPCA_04}
&r^*(\gamma, 2\rho^*)  \leq \max\left(\sqrt{\frac{2B\rho^*}{\gamma}}\left(\frac{12}{N}\log\left(\frac{ed}{k}\right)\right)^{1/4} ; D\sqrt{\frac{K}{N\theta}}\right) \nonumber \\ 
& \leq \max\left(\sqrt{\frac{2B}{\gamma}}\left(\frac{12}{N}\log\left(\frac{ed}{k}\right)\right)^{1/4}\frac{\sqrt{D_1k}}{(N\theta^2)^{1/4}}\max\left(\sqrt{k}\log\left(\frac{ed}{k}\right)^{1/4} ; K^{1/4}\right) ; D\sqrt{\frac{K}{N\theta}}\right) \nonumber \\ 
& \leq \max\left(\sqrt{\frac{2BD_1}{\gamma N\theta}}12^{1/4}\max\left(\sqrt{k}\log\left(\frac{ed}{k}\right)^{1/4} ; K^{1/4}\right)^2 ; D\sqrt{\frac{K}{N\theta}}\right) \leq \frac{C}{\sqrt{N\theta}}  \max\left(k\sqrt{\log\left(\frac{ed}{k}\right)} ; \sqrt{K}\right),
\end{align}
where $C:=\max\left(12^{1/4}\sqrt{2BD_1\gamma^{-1}} ; D\right)$. Combining (\ref{eq:proof_theo:main-RMOM-SPCA_01}), (\ref{eq:proof_theo:main-RMOM-SPCA_03}) and (\ref{eq:proof_theo:main-RMOM-SPCA_04}), we finally get that, with probability at least $1-2\exp(-72K/625)$:
\begin{align*}
&\norm{\hat{Z}^{\mathrm{RMOM}}_{K, \lambda}-Z^*}_1 \leq 2D_1 \frac{k}{\sqrt{N\theta^2}} \max\left(k\sqrt{\log\left(\frac{ed}{k}\right)} ; \sqrt{K}\right) \\
&\norm{\hat{Z}^{\mathrm{RMOM}}_{K, \lambda} - Z^*}_2 \leq \frac{\sqrt{2}C}{\sqrt{N\theta^2}}  \max\left(k\sqrt{\log\left(\frac{ed}{k}\right)} ; \sqrt{K}\right)
\end{align*} and
\begin{equation*}
 P\cL_{\hat{Z}^{\mathrm{RMOM}}_{K, \lambda}} \leq \frac{93C^2}{100N\theta} \max\left(k^2 \log\left(\frac{ed}{k}\right); K\right).
 \end{equation*} This concludes the proof.

\subsubsection{Proof of Corollary \ref{coro:main_SPCA_RMOM}}
 From Theorem \ref{theo:main-RMOM-SPCA}, we get the existence of a universal constant $C_2 > 0$ such that with probability at least $1-\exp(-72K/625)$, $\norm{\hat{Z}^{\mathrm{RMOM}}_{K, \lambda} - Z^*}_2 \leq C_2(N\theta^2)^{-1/2} \max\left(k\sqrt{\log\left(ed/k\right)} ; \sqrt{K}\right)$. Now, we can use Davis-Kahan sin-theta theorem (see Corollary~1 in \cite{yu2014useful}) to get the existence of a universal constant $c_0>0$ such that $\mathrm{sin}(\Theta(\hat{\beta},\beta^*)) = (1/\sqrt{2})\norm{\hat{\beta}\hat{\beta}^\top-\beta^*(\beta^*)^\top}_2  \leq (c_0/g) \norm{\hat{Z}^{\mathrm{RMOM}}_{K, \lambda} - Z^*}_2 $ where $g := \lambda_1 - \lambda_2$ ($\lambda_i$ being the $i^{\mathrm{th}}$ largest eigen value of $Z^*$) is the spectral gap of $Z^*$. Here, we know that $Z^*=\beta^*(\beta^*)^\top$ is rank one, with $1$ as order one eigen value and $0$ as order $d-1$ eigen value. Then we get $g=1$, which leads us to the desired result, with $D=\sqrt{2}c_0\times C_2$.

\section{Appendix}\label{sec:appendix}

\subsection{Distance metric learning: convexity of the constraint set}\label{app:distance_metric_learning_convexity}

Here we show that the constraint set $\cC$ of the ERM estimator of the distance metric learning problem presented in Section \ref{sec:intro} is convex. We recall the definition of this set:
\[
\cC := \left\{Z\in\bR^{d\times d} : Z\succeq 0, \sum_{i, j=1}^M\inr{\left(Y_i-Y_j\right)\left(Y_i-Y_j\right)^\top, Z   }^{1/2}\geq 1\right\}
\]
where $(Y_i)_{i=1}^N$ are $N$ given points in $\bR^d$. Fot the sake of simplicity, we define, for $(i,j)\in[d]^2$, $V_{ij} = \left(Y_i-Y_j\right)\in\bR^d$. Let $Z_1$ and $Z_2$ be two elements of $\cC$, and consider $t\in[0, 1]$. Let us show that $Z^\prime = tZ_1+(1-t)Z_2$ still belongs to $\cC$. We have:
\begin{align*}
\left(\sum_{i, j=1}^M\inr{V_{ij}V_{ij}^\top, Z^\prime}^{1/2}\right)^2 & = \sum_{(i,j)\in[N]^2}\inr{V_{ij}V_{ij}^\top, Z^\prime} + \sum_{(i,j) \neq (pq)} \inr{V_{ij}V_{ij}^\top, Z^\prime}^{1/2}\inr{V_{pq}V_{pq}^\top, Z^\prime}^{1/2} \\ 
& \geq \sum_{(i,j)\in[N]^2}\inr{V_{ij}V_{ij}^\top, Z^\prime} =  t \sum_{(i,j)\in[N]^2}\inr{V_{ij}V_{ij}^\top, Z_1} + (1-t)\sum_{(i,j)\in[N]^2}\inr{V_{ij}V_{ij}^\top, Z_2} \\ 
& = t\left(\sqrt{\sum_{(i,j)\in[N]^2}\inr{V_{ij}V_{ij}^\top, Z_1}}\right)^2 + (1-t)\left(\sqrt{\sum_{(i,j)\in[N]^2}\inr{V_{ij}V_{ij}^\top, Z_2}}\right)^2 \\ 
& \geq t \left(\sum_{(i,j)\in[N]^2}\inr{V_{ij}V_{ij}^\top, Z_1}^{1/2}\right)^2 + (1-t) \left(\sum_{(i,j)\in[N]^2}\inr{V_{ij}V_{ij}^\top, Z_2}^{1/2}\right)^2\\ 
& \geq t + (1-t) = 1
\end{align*}
since each $\sum_{(i,j)\in[N]^2}\inr{V_{ij}V_{ij}^\top, Z_{\ell}}^{1/2}$, $\ell\in\left\{1, 2\right\}$, is larger or equal to one, as $Z_{\ell}\in\cC$. Then, $Z^\prime \in \cC$. We conclude that $\cC$ is convex.

\subsection{A property of local complexity fixed points} 
\label{sub:a_property_of_local_complexity_fixed_points}
Let $H$ be a Hilbert space and $\cC\subset H$. We consider a linear loss function defined for all $Z\in\cC$ by $\ell_Z:X\in H \to -\inr{X,Z}$ and its associated oracle over $\cC$: $Z^*\in\argmin_{Z\in\cC}P \ell_Z$. The excess loss function of $Z\in\cC$ is defined as $\cL_{Z}=\ell_Z-\ell_{Z^*}$. Let $\norm{\cdot}$ be a norm defined (at least) over the span of $\cC$. Let $G:H\to \bR$ be a function. For all $\rho>0$ and $r>0$, we consider the localized model $\cC_{\rho, r} = \{Z\in\cC: \norm{Z-Z^*}\leq \rho, G(Z-Z^*)\leq r\}$ with respect to a $G$ localization and the associated Rademacher complexity 
\begin{equation*}
E(r, \rho) = \bE \left[\sup_{Z\in\cC_{\rho, r}} \Big|\frac{1}{N}\sum_{i=1}^N \sigma_i \cL_Z(X_i)\Big|\right]
\end{equation*}and variance term
\begin{equation*}
V(r,\rho) = \sup_{Z\in\cC_{\rho, r}} \sqrt{{\rm Var}(\cL_Z)}.
\end{equation*}Let $\theta$ and $\tau$ be two positive constants. We consider a local complexity fixed point: for all $\rho>0$,
\begin{equation*}
r^*(\rho) = \inf\left(r>0: \max\left(\theta E(r,\rho), \tau V(r,\rho)\right)\leq r^2\right).
\end{equation*}

\begin{Proposition}\label{prop:prop_local_fixed_point_appendix}We assume that $\cC$ is star-shaped in $Z^*$. We assume that $G$ is such that for all $\alpha\geq 1$ and all $Z\in\cC, G(\alpha(Z-Z^*))\geq \alpha G(Z-Z^*)$. Then, for all $\rho>0$ and $b\geq1$, we have $r^*(\rho)\leq r^*(b\rho)\leq \sqrt{b}r^*(\rho)$.
\end{Proposition}

\textit{Proof.} Let $\rho>0$ and $b\geq1$. For all $r>0$, $\cC_{\rho, r}\subset \cC_{b\rho, r}$ and so $r^*(\rho)\leq r^*(b\rho)$. Let us now prove the second inequality.

We start with some homogeneity property of the complexity and variance terms:
\begin{equation}\label{eq:prop_comp_et_var}
E(\sqrt{b}r, b\rho) \leq b E(r,\rho) \mbox{ and } V(\sqrt{b}r, b\rho) \leq b V(r,\rho). 
\end{equation}We prove \eqref{eq:prop_comp_et_var} for the complexity term, the proof for the variance term is identical. Let $Z\in\cC_{b\rho,\sqrt{b} r}$ and define $Z_0$ such that $Z=Z^*+b(Z_0-Z^*)$. Since $b\geq 1$ and $\cC$ is star-shaped in $Z^*$, $Z_0\in\cC$. Moreover, $b\norm{Z_0-Z^*}=\norm{Z-Z^*}\leq b \rho$ and, by the property of $G$, $bG(Z_0-Z^*)\leq G(Z-Z^*)\leq b r^2$. We conclude that $Z_0\in\cC_{\rho, r}$. Moreover, by linearity of the loss function, we have $\cL_{Z}=b \cL_{Z_0}$. We deduce that 
\begin{equation}
\sup_{Z\in \cC_{b\rho,\sqrt{b} r}} \Big|\frac{1}{N}\sum_{i=1}^N \sigma_i \cL_Z(X_i)\Big| \leq b \sup_{Z_0\in \cC_{\rho, r}} \Big|\frac{1}{N}\sum_{i=1}^N \sigma_i \cL_{Z_0}(X_i)\Big| 
\end{equation}and so \eqref{eq:prop_comp_et_var}  holds for the complexity term. It also holds for the variance using similar tools.

Next, it follows from \eqref{eq:prop_comp_et_var} that
\begin{align*}
r^*(b \rho) &= \inf\left(r>0: \max\left(\theta E(r,b\rho), \tau V(r,b\rho)\right)\leq r^2\right)=\inf\left(r>0: \max\left(\theta E\left(\sqrt{b}\frac{r}{\sqrt{b}},b\rho\right), \tau V\left(\sqrt{b}\frac{r}{\sqrt{b}},b\rho\right)\right)\leq r^2\right)\\
& \leq \inf\left(r>0: \max\left(\theta E\left(\frac{r}{\sqrt{b}},\rho\right), \tau V\left(\frac{r}{\sqrt{b}},\rho\right)\right)\leq \left(\frac{r}{\sqrt{b}}\right)^2\right)\leq \sqrt{b}r^*(\rho).
\end{align*}\endproof

\subsection{A property of the sparsity equation} 
\label{sub:a_property_of_the_sparsity_equation}
We consider the same setup as in Section~\ref{sub:a_property_of_local_complexity_fixed_points} and define for all $\rho>0$,
\begin{equation*}
H_\rho = \left\{Z\in\cC: \norm{Z-Z^*}=\rho, G(Z-Z^*)\leq (r^*(\rho))^2\right\}, \Gamma_{Z^*}(\rho) = \bigcup_{Z:\norm{Z-Z^*}\leq \rho/20}\partial \norm{\cdot}(Z)
\end{equation*}and  $\Delta(\rho) = \inf_{Z\in H_\rho}\sup_{\Phi\in\Gamma_{Z^*}(\rho)}\inr{\Phi, Z-Z^*}$. In the previous section we said that $\rho$ satisfies the sparsity equation when $\Delta(\rho)\geq c_0 \rho$ where $0<c_0<1$ is some absolute constant. In the following result we show that if $\rho$ satisfies the sparsity equation then any number larger than $\rho$ also satisfies this equation.

\begin{Proposition}\label{prop:prop_sparsity_equation}We assume that $\cC$ is star-shaped in $Z^*$. We assume that $G$ is such that for all $\alpha\geq 1$ and all $Z\in\cC, G(\alpha(Z-Z^*))\geq \alpha G(Z-Z^*)$. Let $0<c_0<1$. Then, for all $\rho>0$ and $b\geq1$, if $\rho$ is such that $\Delta(\rho)\geq c_0\rho$ then  $\Delta(b\rho)\geq c_0b\rho$. 
\end{Proposition}

\textit{Proof.} Let $\rho>0$ be such that $\Delta(\rho)\geq c_0\rho$ and let $b\geq1$. Let $Z\in H_{b\rho}$. Let us show that there exists $\Phi\in \Gamma_{Z^*}(b\rho)$ such that $\inr{\Phi, Z-Z^*}\geq c_0 b\rho$. 

Let $Z_0$ be such that $Z=Z^*+b(Z_0-Z^*)$. Since $b\geq 1$ and $\cC$ is star-shaped in $Z^*$, $Z_0\in \cC$. Moreover, $b\norm{Z_0-Z^*}=\norm{Z-Z^*}=b\rho$ and, using the property of $G$ and Proposition~\ref{prop:prop_local_fixed_point_appendix}, $b G(Z_0-Z^*)\leq G(Z-Z^*)\leq  (r^*(b\rho))^2\leq b (r^*(\rho))^2$. Therefore, we have $Z_0\in H_{\rho}$. But, since we assumed that $\Delta(\rho)\geq c_0\rho$, there exists $\Phi\in  \Gamma_{Z^*}(\rho)$ such that $\inr{\Phi,Z_0-Z^*}\geq c_0 \rho$ and so $\inr{\Phi,Z-Z^*}\geq c_0 b\rho$. We conclude the proof by noting that $\Gamma_{Z^*}(\rho)\subset \Gamma_{Z^*}(b\rho)$ and so $\Phi\in  \Gamma_{Z^*}(b\rho)$. 
\endproof





\end{document}